\documentclass[11pt]{report}
\usepackage[latin1]{inputenc}
\usepackage{amsmath,amsthm,amssymb}
\newtheorem{teo}{Theorem}
\newtheorem{defin}{Definition}
\newtheorem{prop}{Proposition}
\newtheorem{cor}{Corollary}
\newtheorem{ex}{Example}
\newtheorem{lemma}{Lemma}
\newtheorem{rem}{Remark}

\def\eq#1{(\ref{#1})}
\def\neweq#1{\begin{equation}\label{#1}}
\def\endeq{\end{equation}}
\newcommand{\ii}{\^{\i}}
\newcommand{\ep}{\varepsilon}
\newcommand{\ab}{\^a}
\newcommand{\ai}{such that }
\newcommand{\ri}{\rightarrow }

\newcommand{\fla}{the functional }

\newcommand{\RR}{\mathbb R}

\newcommand{\pth}{\mbox{$p(t_h)$}}
\newcommand{\funu}{\mbox{$f(\pth -hw(t_h))$} }
\newcommand{\ptn}{\mbox{$p(t_{h_{\scriptstyle n}})$} }
\newcommand{\fdoi}{\mbox{$f(\ptn  -h_{n}w(t_{h_{\scriptstyle n}}))$} }
\newcommand{\limn}{\mbox{$\displaystyle\lim_{n\rightarrow\infty}$} }

\newcommand{\pslie}{\mbox{({\rm s-PS})$_a$} }
\newcommand{\psg}{\mbox{\rm(s-PS)} }
\newcommand{\fxr}{\mbox{$f:X\rightarrow \RR$} }
\newcommand{\catx}{\mbox{{\rm Cat}$_{X}$} }
\newcommand{\catpi}{\mbox{{\rm Cat}$_{\pi (X)}$} }
\newcommand{\psz}{\mbox{({\rm PS})$_{Z}$} }
\newcommand{\pszc}{\mbox{({\rm PS})$_{Z,c}$} }
\newcommand{\crfc}{\mbox{{\rm Cr}$(f,c)$} }
\newcommand{\ginf}{\mbox{$\underline{g}$} }
\newcommand{\gsup}{\mbox{$\overline{g}$} }
\newcommand{\pslc}{\mbox{({\rm PS})$_c$} }
\newcommand{\tmn}{\mbox{$t-{1\over n}$} }
\newcommand{\tpn}{\mbox{$t+{1\over n}$} }
\newcommand{\limeps}{\mbox{$\displaystyle\lim_{\varepsilon\searrow 0}$} }
\newcommand{\lpu}{\mbox{$L^{p+1}(\Omega )$} }
\newcommand{\huo}{\mbox{$H^1_0(\Omega )$} }
\newcommand{\intom}{\mbox{$\displaystyle\int_{\Omega}$} }
\newcommand{\ldo}{\mbox{$L^{2}(\Omega )$} }
\newcommand{\umo}{\mbox{$[u\leq 0]$} }
\newcommand{\upo}{\mbox{$[u\geq 0]$} }

\newcommand{\mfin}{\nolinebreak\makebox[3cm]{\rule{2mm}{2mm}}}
\newcommand{\finf}{\mbox{$\underline{f}$} }
\newcommand{\fsup}{\mbox{$\overline{f}$} }
\newcommand{\spsc}{\mbox{({\rm s-PS})$_c$} }
\newcommand{\Uinf}{\mbox{$\underline{U}$} }
\newcommand{\Usup}{\mbox{$\overline{U}$} }
\newcommand{\uinf}{\mbox{$\underline{u}$} }
\newcommand{\usup}{\mbox{$\overline{u}$} }
\newcommand{\di}{\displaystyle}
\newcommand{\itab}{\left\{\begin{tabular}{ll} }
\newcommand{\ttab}{\end{tabular}\right. }
\begin{document}
\baselineskip16pt
\title{\Large\sc Nonlinear Partial Differential Equations of Elliptic Type}

\author{\sc Vicen\c tiu D. R\u adulescu\\ Department of Mathematics, University of Craiova, 200585 Craiova,
Romania\\ E-mail: {\tt radulescu@inf.ucv.ro}\qquad {\tt
http://inf.ucv.ro/$^\sim$radulescu} }

\date{}
\maketitle
\tableofcontents
\chapter*{Introduction}
\addcontentsline{toc}{chapter}{Introduction}
 Nonlinear Analysis is one of the fields of Mathematics
with the most spectacular development in the last decades of this end of century.
The impressive number of results in this area
is also a consequence of various problems
raised by Physics, Optimization or Economy. In the modelling of natural
phenomena a crucial role is played by the study of partial differential equations
of elliptic type; they arise in every field of science. Consequently, the desire
to understand the solutions of these equations has always a prominent place in
the efforts of mathematicians; it has inspired such diverse fields as
Functional Analysis, Variational Calculus or Algebraic Topology.

The present book is based on a one semester course at the University of Craiova.
The goal of this textbook is to provide the background which is necessary to
initiate   work   on a Ph.D. thesis in Applied Nonlinear Analysis. My purpose
is to provide for the student a broad perspective in the subject, to illustrate
the rich variety of phenomena encompassed by it and to impart a working
knowledge of the most important techniques of analysis of the solutions of
the equations. The level
of this book is aimed at beginning graduate students. Prerequisites include
a truly advanced Calculus course, basic knowledge on Functional Analysis and
PDE, as well as the necessary tools on Sobolev spaces.

Throughout this work we have used intensively two classical
results: the Mountain-Pass Lemma (in its $C^1$ statement!) of
Ambrosetti and Rabinowitz (1973, \cite{AR}) and Ekeland's
Variational Principle (1974, \cite{E1}). We recall in what follows
these celebrated results.

\medskip
{\bf  Mountain Pass Lemma.} {\sl Let $X$ be a real Banach space and let
 $F:X\rightarrow\RR$ be a $C^1$-functional  which satisfies the following
 assumptions:

i) $\ F(0)<0$ and there exists $e\in X\setminus\{ 0\}$ such that $F(e)<0$.

ii) there is some $0<R<\| e\|$ \ai  $F(u)\geq 0$, for all $u\in X$ with
$\| u\| =R$.

Put
$$c=\inf_{p\in{\cal P}}\max_{t\in [0,1]}F(p(t)),$$
where ${\cal P}$ denotes the set of all continuous paths joining
 0 and $e$.

Then the number
 $c$ is an ``almost critical value" of
 $F$, in the sense that there exists a sequence
  $(x_n)$ in $X$ \ai
$$\lim_{n\rightarrow\infty}F(x_n)=c\quad\hbox{and}\quad
\lim_{n\rightarrow\infty}\| F'(x_n)\|_{X^*}=0\, .$$}\rm

\medskip
{\bf  Ekeland's Variational Principle.} {\sl Let $(M,d)$ be a complete metric
space and let
 $\psi :M\rightarrow (-\infty ,+\infty ]$, $\psi\not\equiv +\infty$ be a lower
 semicontinuous functional bounded from below.

Then the following hold:

i) Let $\varepsilon >0$ be arbitrary and let  $z_0\in M$ be \ai
$$\psi (z_0)\leq \inf_{x\in M}\psi (x)+\varepsilon \, .$$

Then, for every $\lambda >0$, there exists $z_{\lambda}\in M$ \ai\, ,
for all $x\in M\setminus\{z_{\lambda}\}$,
$$\psi (z_\lambda )\leq\psi (z_0)\, ,$$
$$\psi (x)>\psi (z_\lambda )-{\varepsilon\over\lambda}\,
d(x,z_\lambda )\, ,$$
$$d(z_\lambda ,z_0)\leq\lambda\, .$$

ii) For each $\varepsilon >0$ and $z_0\in M$, there is some
$z\in M$ such that, for any $x\in M$,
$$\psi (x)\geq\psi (z)-\varepsilon d(x,z)\, ,$$
$$\psi (z)\leq \psi (z_0)-\varepsilon d(z_0,z)\, .$$}\rm

\medskip
In the first part of this book we present the method of
sub and super solutions, which is one of the main tools in Nonlinear Analysis for
finding solutions to a boundary value problem. The proofs are simple
and we give several examples to illustrate the theory. We
continue with another elementary method for finding solutions,
namely the Implicit Function Theorem. The main application is
a celebrated theorem due to H.~Amann  which is related to a bifurcation
problem associated to a convex and positive function.
This kind of equations arises frequently in physics,
biology, combustion-diffusion etc. We mention only
 Brusselator type reactions, the combustion theory, dynamics of population, the
 Fitzhugh-Nagumo system, morphogenese, superconductivity, super-fluids etc.
 We give
complete details in the case where the functional is
asymptotically linear at infinity.

In the third chapter we give some basic results related
to the Clarke generalized gradient of a locally Lipschitz functional
 (see Clarke \cite{Cl1}, \cite{Cl2}). Then we develop a nonsmooth critical point theory which
 enables us to deduce the
 Brezis-Coron-Nirenberg Theorem \cite{BCN},
the ``Saddle Point" Theorem of Rabinowitz \cite{Ra1} or the
Ghoussoub-Preiss Theorem \cite{GP}. The motivation of this study is the following:
some of the strongest tools for proving existence results in PDE are
the ``Mountain Pass" Lemma of Ambrosetti-Rabinowitz and the
Lusternik-Schnirelmann Theorem. These results apply when the solutions of the
given problem are critical points of a suitable ``energetic" functional
 $f$, which is assumed to be of class $C^1$ and defined on a real Banach space.
 A natural question is what happens if the energy functional, associated to a
 given problem in a natural way, fails to be
 differentiable. The results we give here are based on the notion of Clarke
 generalized gradient of a locally Lipschitz functional which is very useful
 for the treatment of many problems arising in the calculus of variations,
 optimal control, hemivariational inequalities etc. Clarke's generalized
 gradient coincides with the usual one if
 $f$ is differentiable or convex. In the classical framework the Fr\'echet
 differential of a $C^1$-functional is a linear and continuous operator. For
 the case of locally Lipschitz maps the property of linearity of the gradient
does not remain valid. Thus, for fixed
 $x\in X$, the directional derivative
$f^0(x,\cdot )$ is subadditive and positive homogeneous and its generalized
gradient
$\partial f(x)$ is a nonempty closed subset of the dual space.

In Chapter 4 the main results are two
 Lusternik-Schnirelmann type theorems. The first one uses the notion of
 critical point for a pairing of operators
(see, e.g., Fucik-Necas-Soucek-Soucek \cite{FNSS}). The second theorem is related to locally Lipschitz
functionals which are periodic with respect to a discrete subgroup and also
uses the notion of Clarke subdifferential.

In the following part of this work we give several applications of the abstract results which
appear in the first  chapters. Here we recall a classical result of Chang
\cite{Ch} and prove multivalued variants of the Brezis-Nirenberg problem, as well as
a solution of the forced pendulum problem, which was studied in the smooth case
in Mawhin-Willem \cite{MW2}. We also study multivalued problems at resonance of
 Landesman-Lazer type. The methods we develop here enable us to study several
 classes of discontinuous problems and all these techniques are based on
 Clarke's generalized gradient theory.
This tool is very useful in the study of critical periodic orbits
of hamiltonian systems (Clarke-Ekeland), the mathematical
programming (Hiriart-Urruty), the duality theory (Rockafellar),
optimal control (V.~Barbu and F.~Clarke), nonsmooth analysis (A.D.
Ioffe and F. Clarke), hemivariational inequalities (P.D.
Panagiotopoulos) etc.
\begin{flushright}{November 2004}\end{flushright}

\chapter{Method of sub and super solutions}
Let $\Omega$ be a smooth bounded domain in $\RR^N$ and consider a
Carath\'eodory function $f(x,u):\overline\Omega\times\RR\ri\RR$
such that $f$ is of class $C^1$ with respect to the variable $u$.
Consider the problem
\neweq{P1}
\itab
$\di-\Delta u=f(x,u)\,,$ & \quad $\mbox{in}\ \Omega$\\
$\di u=0\,,$ & \quad $\mbox{on}\ \partial\Omega\, .$\\
\ttab
\endeq
By solution of the problem \eq{P1} we mean a function
$u\in C^2(\Omega )\cap C(\overline\Omega )$ which satisfies
\eq{P1}.

\begin{defin}\label{subsol}
A function $\Uinf\in C^2(\Omega )\cap C(\overline\Omega )$ is
said to be subsolution of the problem \eq{P1} provided that
$$
\itab
$\di -\Delta \Uinf\leq f(x,\Uinf )\,,$ & \quad $\mbox{in}\ \Omega$\\
$\di \Uinf\leq 0\,,$ & \quad $\mbox{on}\ \partial\Omega\, .$\\
\ttab
$$
Accordingly, if the signs are reversed in the above inequalities,
we obtain the definition of  a supersolution
$\Usup$
for the
problem \eq{P1}.
\end{defin}

\begin{teo}\label{ssol}
Let $\Uinf$ (resp., $\Usup$) be a subsolution (resp., a supersolution)
to the problem \eq{P1} such that $\Uinf\leq \Usup$ in $\Omega$.
The following hold:

(i) there exists a solution $u$ of \eq{P1} which, moreover,
satisfies $\Uinf\leq u\leq\Usup$;

(ii) there exist a minimal and a maximal solution $\uinf$ and
$\usup$ of the problem \eq{P1} with respect to the interval
$[\Uinf ,\Usup ]$.
\end{teo}

\begin{rem}\label{maxi}
The existence of the solution in this theorem, as well as
the maximality (resp., minimality) of  solutions given by (ii)
have to be
understood with respect to the given pairing of ordered sub
and supersolutions.
It is very possible that \eq{P1} has solutions which are {\bf not}
in the interval $[\Uinf ,\Usup ]$. It may also happen that \eq{P1}
has {\bf no} maximal or minimal solution. Give such an exemple!
\end{rem}

\begin{rem}\label{ordo}
The hypothesis $\Uinf\leq \Usup$ is {\bf not} automatically
fulfilled for arbitrary sub and supersolution of \eq{P1}.
Moreover, it may occur that $\Uinf > \Usup$ on the {\bf whole} of
$\Omega$. An elementary example is the following: consider the
eigenvalue problem
$$
 \itab
 $ \di -\Delta u=\lambda_1u\,,$ & \quad $\mbox{in}\ \Omega $\\
 $ \di u=0\,,$ & \quad $\mbox{on}\ \partial\Omega\, .$\\
 \ttab
 $$
 We know that all solutions of this problem are of the form
 $u=Ce_1$, where $C$ is a real constant and $e_1$ is not vanishing
 in $\Omega$, say $e_1(x)>0$, for any $x\in\Omega$.
 Choose $\Uinf =e_1$ and $\Usup =-e_1$. Then $\Uinf$ (resp.,
 $\Usup$) is subsolution (resp., supersolution) to the problem
 \eq{P1}, but $\Uinf >\Usup$.
 \end{rem}
{\bf Proof of Theorem \ref{ssol}.} (i)
 Let $g(x,u):=f(x,u)+au$, where $a$ is a real constant.
 We can choose $a\geq 0$ sufficiently large so that the map
 $\RR\ni u\longmapsto  g(x,u)$ is increasing on
 $[\Uinf (x),\Usup (x)]$, for every $x\in\Omega$. For this aim, it
 is enough to have $a\geq 0$ and
 $$a\geq\max\left\{ -f_u(x,u); x\in\overline\Omega \ \mbox{and}\
 u\in [\Uinf (x),\Usup (x)]\right\}\, .$$

For this choice of $a$ we define the sequence of
functions $u_n\in C^2(\Omega )\cap
C(\overline\Omega )$ as follows: $u_0=\Usup$ and, for every
$n\geq 1$, $u_n$ is the unique solution of the linear problem
\neweq{udefin}
 \itab
 $\di -\Delta u_n+au_n=g(x,u_{n-1})\,,$ &\quad $\mbox{in}\ \Omega $\\
 $ \di u_n=0\,,$ & \quad $\mbox{on}\ \partial\Omega\, .$\\
 \ttab
 \endeq
 {\bf Claim}: $\Uinf\leq\cdots \leq u_{n+1}\leq u_n\leq\cdots
 \leq u_0=\Usup$.\\
 {\it Proof of Claim}. \rm Our arguments use in an essential manner
 the Weak Maximum Principle. So, in order to prove that
 $u_1\leq \Usup$ we have, by the definition of $u_1$,
 $$
 \itab
 $\di -\Delta (\Usup -u_1)+a(\Usup -u_1)
 \geq g(x,\Usup )-g(x,\Usup )=0\, ,$ &\quad $\mbox{in}\ \Omega $\\
 $\di \Usup -u_1\geq 0\,,$ & \quad $\mbox{on}\ \partial\Omega\, .$\\
 \ttab
 $$
Since the operator $-\Delta +aI$ is coercive,
it follows that $\Usup \geq u_1$ in $\Omega$.
 For the proof of $\Uinf\leq u_1$ we observe that $\Uinf\leq 0=
 u_1$ on $\partial\Omega$ and, for every $x\in\Omega$,
 $$-\Delta (\Uinf -u_1)+a(\Uinf -u_1)\leq f(x,\Uinf )+
 a\Uinf -g(x,\Usup )\leq 0\, ,$$
 by the monotonicity of $g$. The Maximum Principle implies
 $\Uinf\leq u_1$.

 Let us now assume that
 $$\Uinf\leq\cdots \leq u_{n}\leq u_{n-1}\leq\cdots
 \leq u_0=\Usup\, .$$
It remains to prove that
 $$
 \Uinf\leq u_{n+1}\leq u_n\, .
 $$
Taking into account the equations satisfied by $u_n$ and $u_{n+1}$
we obtain
$$
 \itab
 $\di -\Delta (u_n -u_{n+1})+a(u_n -u_{n+1})
 = g(x,u_{n-1} )-g(x,u_n )\geq 0 \, ,$ & \quad $\mbox{in}\ \Omega$ \\
 $ \di u_n -u_{n+1}\geq 0\,,$ & \quad $\mbox{on}\ \partial\Omega\, ,$\\
 \ttab
$$
 which implies $u_n\geq u_{n+1}$ in $\Omega$.

 On the other hand, by
 $$
 \itab
 $\di -\Delta \Uinf +a\Uinf
 \leq g(x,\Uinf )\,, $ & \quad $\mbox{in}\ \Omega $\\
 $\di \Uinf \leq 0\,,$ &\quad $\mbox{on}\ \partial\Omega$\\
 \ttab
 $$
 and the definition of $u_{n+1}$ we have
 $$
 \itab
 $\di-\Delta (u_{n+1}-\Uinf )+a(u_{n+1} -\Uinf )
 \geq g(x,u_{n} )-g(x,\Uinf )\geq 0\,, $ & \quad $\mbox{in}\ \Omega$ \\
 $\di u_{n+1} -\Uinf \geq 0\,,$ &\quad $\mbox{on}\ \partial\Omega\, .$\\
 \ttab
 $$
 Again, by the Maximum Principle, we deduce that $\Uinf \leq u_{n+1}$
 in $\Omega$, which completes the proof of the Claim.

It follows that there exists a function $u$ such that,
  for every fixed $x\in\Omega$,
  $$u_n(x)\searrow u(x)\qquad\mbox{as}\ n\ri\infty\, .$$
 Our aim is to show that we can pass to the limit in \eq{udefin}.
 For this aim we use a standard bootstrap argument. Let
 $g_n(x):=g(x,u_n(x))$. We first observe
 that the sequence $(g_n)$ is bounded in $L^\infty (\Omega )$, so
 in every $L^p(\Omega )$ with $1<p<\infty$. It follows by
 \eq{udefin} and standard Schauder estimates that the sequence
 $(u_n)$ is bounded in $W^{2,p}(\Omega )$, for any $1<p<\infty$. But
 the space $W^{2,p}(\Omega )$ is continuously embedded in
 $C^{1,\alpha }(\overline\Omega )$, for $\alpha =1-\frac{N}{2p}$,
 provided that $p>\frac{N}{2}$. This implies that $(u_n)$ is
 bounded in $C^{1,\alpha }(\overline\Omega )$. Now, by
 standard estimates in H\"older spaces we deduce that $(u_n)$ is
 bounded in $C^{2,\alpha }(\overline\Omega )$. Since
 $C^{2,\alpha }(\overline\Omega )$ is compactly embedded in
 $C^{2}(\overline\Omega )$, it follows that, passing eventually
 at a subsequence,
 $$u_n\ri u\qquad\mbox{in}\ C^{2}(\overline\Omega )\, .$$
 Since the sequence is monotone we obtain that the whole sequence
 converges to $u$ in $C^2$. Now, passing at the limit in
 \eq{udefin} as $n\ri\infty$ we deduce that $u$ is solution of
 the problem \eq{P1}.

(ii) Let us denote by $\usup$ the solution obtained with the above
technique and choosing $u_0=\Usup$. We justify that $\usup$ is a
maximal solution with respect to the given pairing $(\Uinf , \Usup
)$. Indeed, let $u\in [\Uinf ,\Usup ]$ be an arbitrary solution.
With an argument similar to that given in the proof of (i) but
with respect to the pairing of ordered sub-supersolutions
$(u,\Usup )$ we obtain that $u\leq u_n$, for any $n\geq 0$, which
implies $u\leq \usup$. \qed

\medskip
Let us now  assume that $f$ is a continuous function. We give in what follows
an elementary variational
proof for the existence of a solution to the problem
(\ref{P1}), provided that  sub and supersolution $\Uinf$
and $\Usup$ exist, with $\Uinf\leq\Usup$. For this aim, let
$$E(u)=\frac{1}{2}\intom \mid \nabla u\mid^2-\intom F(x,u)$$
be the energy functional associated to the problem (\ref{P1}).
Here, $F(x,u)=\int_0^uf(x,t)dt$.

Set
$$f_0(x,t)=
\itab
$ \di f(x,t)\,,$ & \qquad $\mbox{if}\ \, \Uinf (x)<t<\Usup (x) $\\
$ \di f(x,\Usup (x))\,,$ &\qquad $\mbox{if}\ \, t\geq\Usup (x) $\\
$\di f(x,\Uinf (x))\,,$ &\qquad $\mbox{if}\ \, t\leq\Uinf (x)\, .$\\
\ttab
$$
The associated energy functional is
$$E_0(u)=\frac{1}{2}\intom \mid \nabla u\mid^2-\intom F_0(x,u)\, ,$$
with an appropriate definition for $F_0$.

We observe the following:

\noindent - $E_0$ is well defined on $H^1(\Omega )$, since
$f_0$ is uniformly bounded, so $F_0$ has a sublinear growth;

\noindent - $E_0$ is weak lower semicontinuous;

\noindent - the first term of $E_0$ is the dominating one at $+\infty$
and, moreover,
$$\lim_{\| u\|\ri\infty}E_0(u)=+\infty\, .$$

Let
$$\alpha =\inf_{u\in H^1_0(\Omega )} E_0(u)\, .$$
We show in what follows that $\alpha$ is attained. Indeed, since $E_0$
is coercive, there exists a minimizing sequence
$(u_n)\subset H^1_0(\Omega )$. We may assume without loss of generality that
$$u_n\rightharpoonup u\, ,\qquad\mbox{weakly in}\ \, \huo\, .$$
So, by the lower semicontinuity of $E_0$ with respect to the weak
topology,
$$\frac{1}{2}\intom \mid\nabla u_n\mid^2-\intom F_0(x,u_n)\leq
\alpha +o(1)\, .$$
This implies $E_0(u)=\alpha$.
Now, since $u$ is a critical point of $E_0$, it follows that it satisfies
$$-\Delta u=f_0(x,u)\, ,\qquad\mbox{in}\ \, {\cal D}'(\Omega )\, .$$
The same bootstrap argument as in the
above proof shows that $u$ is smooth.

We prove in what follows that $\Uinf\leq u\leq\Usup$. Indeed, we have
 $$-\Delta\Uinf\leq f(x,\Uinf)\, ,\qquad
 \mbox{in}\ \, \Omega\, .$$
Therefore
$$-\Delta (\Uinf -u)\leq f(x,\Uinf)-f_0(x,u)\, .$$
After multiplication by $(\Uinf-u)^+$ in this inequality and integration
over $\Omega$ we find
$$\intom\mid \nabla (\Uinf -u)^+\mid^2\leq\intom (\Uinf -u)^+
\left( f(x,\Uinf )-f_0(x,u)\right)\, .$$
Taking into account the definition of $f_0$ we obtain
 $$\intom\mid \nabla (\Uinf-u)^+\mid^2=0$$
 which implies $\nabla (\Uinf -u)^+=0$ in $\Omega$. Therefore,
 $\Uinf\leq u$ in $\Omega$. \qed

\smallskip
We can interpret a solution $u$ of the problem \eq{P1} as an
equilibrium solution of the associated parabolic problem
\neweq{parabb}
\left\{\begin{array}{lll}
 v_t-\Delta v=f(x,v)\, ,\quad\mbox{in}\ \Omega\times (0,\infty )\\
 v(x,t)=0\, ,\quad\mbox{on}\ \partial\Omega\times (0,\infty )\\
 v(x,0)=u_0(x)\, ,\quad\mbox{in}\ \Omega\, .
\end{array}\right.
\endeq
Suppose that the initial data $u_0(x)$ does not deviate too much
from a stationary state $u(x)$. Does the solution of \eq{parabb}
return to $u(x)$ as $t\ri\infty$? If this is the case then the
solution $u$ of the problem \eq{P1} is said to be stable. More
precisely, a solution $u$ of \eq{P1} is called {\bf stable} if for
every $\ep >0$ there exists $\delta >0$ such that $\|
u(x)-v(x,t)\|_{L^\infty (\Omega\times (0,\infty ))}<\ep$, provided
that $\| u(x)-u_0(x)\|_{L^\infty (\Omega )}<\delta$. Here,
$v(x,t)$ is a solution of problem \eq{parabb}. We establish in
what follows that the solutions given by the method of sub and
supersolutions are, generally, stable, in the following sense.

\begin{defin}\label{stable}
A solution $u$ of the problem \eq{P1} is said to be stable
provided that the first eigenvalue of the linearized
operator at $u$ is positive, that is,
$\lambda_1\left(-\Delta -f_u(x,u)\right)>0.$
The solution $u$ is called semistable if
$\lambda_1\left(-\Delta -f_u(x,u)\right)\geq 0.$
\end{defin}
In the above definition we understand the first eigenvalue of the
linearized operator with respect to homogeneous Dirichlet boundary condition.

\begin{teo}\label{stabil}
Let $\Uinf$ (resp., $\Usup$) be subsolution (resp.,
 supersolution) of the problem
\eq{P1} such that $\Uinf\leq\Usup$ and let $\uinf$
(resp., $\usup$) be the corresponding minimal (resp., maximal)
solution of \eq{P1}. Assume that $\Uinf$ is not a solution
of \eq{P1}. Then $\uinf$ is semistable. Furthermore, if
$f$ is concave, then $\uinf$ is stable.

Similarly, if
$\Usup$ is not a solution
then $\usup$ is semistable and, if
$f$ is convex, then $\usup$ is stable.
\end{teo}
{\bf Proof.}
Let $\lambda_1=
\lambda_1\left(-\Delta -f_u(x,\uinf )\right)$ and let
$\varphi_1$ be the corresponding eigenfunction, that is,
$$
\itab
$\di -\Delta\varphi_1-f_u(x,\uinf )\varphi_1=\lambda_1\varphi_1
\, ,$& \quad $\mbox{in}\ \Omega $\\
$\di \varphi_1=0\,, $ & \quad $\mbox{on}\ \partial\Omega\, .$\\
\ttab
$$
We can suppose, without loss of generality, that $\varphi_1>0$
in $\Omega$. Assume, by contradiction, that $\lambda_1<0$.
Let us consider the function $v:=\uinf -\ep\varphi_1$, with
$\ep >0$. We prove in what follows that the following hold:

(i) $v$ is a supersolution to the problem \eq{P1}, for
$\ep$ small enough;

(ii) $v\geq\Uinf$.

By (i), (ii) and Theorem \ref{ssol} it follows that there
exists a solution $u$ such that $\Uinf\leq u\leq v<
\uinf$ in $\Omega$, which contradicts the minimality of
$\uinf$ and the hypothesis that $\Uinf$ is not a solution.

In order to prove (i), it is enough to show that
$$-\Delta v\geq f(x,v)\quad\mbox{in}\ \Omega\, .$$
But
$$\begin{array}{lll}
\Delta v+ f(x,v)=\Delta\uinf -\ep\Delta\varphi_1+f(x,\uinf -
\ep\varphi_1)= \\
-f(x,\uinf )+\ep\lambda_1\varphi_1+
\ep f_u(x,\uinf )\varphi_1
+f(x,\uinf -\ep\varphi_1)= \\
-f(x,\uinf )+\ep\left(\lambda_1\varphi_1+
f_u(x,\uinf )\varphi_1\right) +f(x,\uinf )-
\ep f_u(x,\uinf )\varphi_1+o(\ep\varphi_1)= \\
\ep\lambda_1\varphi_1+o(\ep )\varphi_1=
\varphi_1\left( \ep\lambda_1+o(\ep )\right)=\varphi_1\ep
\left(\lambda_1+o(1)\right)\leq 0\, ,
\end{array}
$$
provided that $\ep >0$ is sufficiently small.

Let us now prove (ii). We observe that $v\geq \Uinf$ is equivalent to
$\uinf -\Uinf\geq\ep\varphi_1$, for small
$\ep$. But $\uinf -\Uinf\geq 0$ in $\Omega$. Moreover,
$\uinf -\Uinf\not\equiv 0$, since $\Uinf$ is not solution.
Now we are in position to apply the Hopf Strong Maximum
Principle in the following variant: assume $v$ satisfies
$$
\itab
$\di -\Delta w+aw=f(x)\geq 0\, ,$ &\quad $\mbox{in}\ \Omega $\\
$\di w\geq 0 \, ,$ & \quad $\mbox{on}\ \partial\Omega\, ,$\\
\ttab
$$
where $a$ is a nonnegative number. Then $w\geq 0$ in
$\Omega$ and the following alternative holds: either

\noindent (i) $w\equiv 0$ in $\Omega$

\noindent or

\noindent (ii) $w>0$ in $\Omega$ and $\frac{\partial w}
{\partial\nu}<0$ on the set $\{ x\in\partial\Omega ;
w(x)=0\}.$

Let $w=\uinf -\Uinf\geq 0$. We have
$$
-\Delta w+aw=f(x,\uinf )+\Delta\Uinf +a(\uinf -\Uinf )\geq
f(x,\uinf )-f(x,\Uinf )+a(\uinf -\Uinf )\, .
$$
So, in order to have $-\Delta w+aw\geq 0$ in $\Omega$, it is
sufficient to choose $a\geq 0$ so that the mapping $\RR\ni
u\longmapsto  f(x,u)+au$ is increasing on $[\underline
U(x),\overline U(x)]$, as already done in the proof of Theorem
\ref{ssol}. Observing that $w\geq 0$ on $\partial\Omega$ and
$w\not\equiv 0$ in $\Omega$ we deduce by the Strong Maximum
Principle that
$$w>0\ \, \mbox{in}\ \Omega\ \, \mbox{and}\ \,
\frac{\partial w}{\partial\nu}<0\ \, \mbox{on}\
\{ x\in\partial\Omega ; \uinf (x)=\Uinf (x)=0\}\, .$$

We prove in what follows that we can choose $\ep >0$
sufficiently small so that $\ep\varphi_1\leq w$. This
is an interesting consequence of the fact that the normal derivative
is negative in the points of the boundary where the
function vanishes. Arguing by contradiction, there exist
a sequence $\ep_n\ri 0$ and $x_n\in\Omega$ such that
\neweq{normal}
(w-\ep_n\varphi_1)(x_n)<0\, .\endeq
Moreover, we can choose the points $x_n$ with the
additional property
\neweq{gradientw}
\nabla (w-\ep_n\varphi_1)(x_n)=0\, .\endeq
But, passing eventually at a subsequence, we can assume that $x_n\ri x_0
\in\overline\Omega$. It follows now by \eq{normal} that
$w(x_0)\leq 0$ which implies $w(x_0)=0$, that is,
$x_0\in\partial\Omega$. Furthermore, by \eq{gradientw},
$\nabla w(x_0)=0$, a contradiction, since
$\frac{\partial w}{\partial\nu}(x_0)<0$, by the Strong
Maximum Principle.

Let us now assume that $f$ is concave. We have to show
that $\lambda_1>0$. Arguing again by contradiction, let
us suppose that $\lambda_1=0$. With the same arguments as
above we can show that $v\geq\Uinf$. If we prove that $v$
is a supersolution then we contradicts the minimality of
$\uinf$. The above arguments do not apply  since, in order to
find a contradiction, the estimate
$$\Delta v+f(x,v)=\ep\varphi_1\left( \lambda_1+o(1)\right) $$
is not relevant in the case
where $\lambda_1=0$. However
$$\begin{array}{lll}
\Delta v+f(x,v)=-f(x,\uinf )+
\ep\left( \lambda_1\varphi_1+f_u(x,\uinf )\varphi_1\right)
+f(x,\uinf -\ep\varphi_1)\leq \\
\ep f_u(x,\uinf )\varphi_1+f_u(x,\uinf )(-\ep\varphi_1)
=0\, .
\end{array}$$
\qed

\medskip
If neither $\Uinf$ nor $\Usup$ are solutions to the
problem \eq{P1} it is natural to ask if
there exists a solution $u$ such that $\Uinf <u<\Usup$
and $\lambda_1(-\Delta -f_u(x,u))>0$. In general such
a situation does not occur, as showed by the following
example: consider the problem
$$
\itab
$\di-\Delta u=\lambda_1u-u^3\, ,$ & \quad $\mbox{in}\ \Omega $\\
$\di u=0\, ,$ & \quad $\mbox{on}\ \partial\Omega\, ,$\\
\ttab
$$
where $\lambda_1$ denotes the first eigenvalue of
$-\Delta$ in $H_0^1(\Omega )$. We remark that we can
choose $\Usup =a$ and $\Uinf =-a$, for every $a>\sqrt{\lambda_1}$.
On the other hand, by Poincar\'e's Inequality,
$$\lambda_1\intom u^2\leq\intom\mid\nabla u\mid^2
=\lambda_1\intom u^2-\intom u^4\, ,$$
which shows that the unique solution is $u=0$.
However this solution is not stable, since $\lambda_1(-\Delta
-f_u(0))=0$.

Another question which arises is under what hypotheses
there exists a global maximal (resp., minimal) solution
 of \eq{P1}, not only with respect to a
prescribed pairing of sub and supersolutions. The
following result shows that a sufficient condition is
that the nonlinearity has a kind of sublinear growth.
More precisely, let us consider the problem
\neweq{P2}
\itab
$\di -\Delta u=f(x,u)+g(x)\, ,$ & \quad $\mbox{in}\ \Omega $\\
$\di u=0\, ,$ & \quad $\mbox{on}\ \partial\Omega\, .$\\
\ttab
\endeq

\medskip
\begin{teo}\label{global}
Assume $g\in C^\alpha (\overline\Omega )$, for some
$\alpha\in (0,1)$ and, for every $(x,u)\in\Omega
\times\RR$,
\neweq{sign}
f(x,u)\, {\rm sign}\, u\leq a\mid u\mid +C\qquad
\mbox{with}\ a<\lambda_1\, .
\endeq
Then there exists a global minimal (resp., maximal)
solution $\uinf$ (resp., $\usup$) to the problem
\eq{P2}.
\end{teo}
{\bf Proof.} Assume without loss of generality that $C>0$.
We choose as supersolution of  \eq{P2}
the unique solution $\Usup$ of the problem
$$
\itab
$\di-\Delta \Usup -a\Usup =C'\, ,$ &\quad $\mbox{in}\ \Omega $\\
$\di \Usup =0\,,$ & \quad $\mbox{on}\ \partial\Omega\, ,$\\
\ttab
$$
where $C'$ is taken such that $C'\geq C+\sup_{\overline \Omega}\,\mid g\mid$.
Since $a<\lambda_1$ it follows by the Maximum Principle
that $\Usup\geq 0$.

Let $\Uinf =-\Usup$ be a subsolution of \eq{P2}. Thus,
by Theorem \ref{ssol}, there exists $\uinf$ (resp.,
$\usup$) minimal (resp., maximal) with respect to
$(\Uinf$, $\Usup)$. We prove in what follows that $\uinf
\leq u\leq \usup$, for {\bf every} solution $u$ of the
problem \eq{P2}. For this aim, it is enough to show
that $\Uinf
\leq u\leq \Usup$. Let us prove that $u\leq\Usup$.
Denote
$$\Omega_0=\{ x\in\Omega ;\ u(x)>0\}\, .$$
Consequently, it is sufficient to show that $u\leq \Usup$ in
$\Omega_0$. The idea is to prove that
$$
\itab
$\di -\Delta (\Usup -u)-a(\Usup -u)\geq 0\, ,$
& \quad $\mbox{in}\ \Omega_0 $\\
$\di \Usup -u\geq 0\,,$ & \quad $\mbox{on}\ \partial\Omega_0\, $\\
\ttab
$$
and then to apply the Maximum Principle. On the one hand,
it is obvious that
$$\Usup -u=\Usup\geq 0\, ,\quad\mbox{on}\ \partial\Omega_0\,.$$
On the other hand,
$$\begin{array}{lll}
-\Delta (\Usup -u)-a(\Usup -u)=-\Delta\Usup -
a\Usup -(-\Delta u-au)\geq\\
C'-f(x,u)+au\geq 0\, ,\qquad\mbox{in}\ \Omega_0\, ,
\end{array}
$$
which ends our proof. \qed

\medskip
Let us now consider the problem
\neweq{P3}
\itab
$\di-\Delta u=f(u)\, ,$ & \quad $\mbox{in}\ \Omega $\\
$\di u=0\, ,$ & \quad $\mbox{on}\ \partial\Omega $\\
$\di u>0\, , $ & \quad $\mbox{in}\ \Omega\, ,$\\
\ttab
\endeq
where
\neweq{x0}
f(0)=0
\endeq
\neweq{lsup}
\limsup_{u\ri +\infty}\frac{f(u)}{u}<\lambda_1\, .
\endeq
Observe that \eq{lsup} implies
$$f(u)\leq au+C\, ,\qquad\mbox{for every}\ u\geq 0\, ,$$
with $a<\lambda_1$ and $C>0$.

Clearly $\Uinf =0$ is a subsolution of \eq{P3}. Choose
$\Usup$ the unique solution of the problem
$$ \itab
$\di-\Delta \Usup -a\Usup =C\, ,$ &\quad $\mbox{in}\ \Omega $\\
$\di \Usup =0\, , $ & \quad $\mbox{on}\ \partial\Omega\,. $\\
\ttab$$
We then obtain a minimal solution $\uinf$
and a maximal solution $\usup\geq 0$. However we can
not state that $\usup >0$ (give an example!). A positive
answer is given by

\begin{teo}\label{sufic}
Assume $f$ satisfies hypotheses \eq{x0}, \eq{lsup} and
\neweq{fprim}
f'(0)>\lambda_1\, .
\endeq
Then there exists a maximal solution $u$ to the problem
\eq{P3} such that $u>0$ in $\Omega$.\end{teo}
{\bf Proof.} The idea is to find another subsolution.
Let $\Uinf =\ep\varphi_1$, where $\varphi_1>0$ is the
first eigenfunction of $-\Delta$ in $H^1_0(\Omega )$.
To obtain our conclusion it is sufficient to
 verify  that for $\ep >0$ small enough we have

 \noindent (i) $\ep\varphi_1$ is a subsolution;

 \noindent (ii) $\ep\varphi_1\leq \Usup$.

Let us verify (i). We observe that
$$f(\ep\varphi_1)=f(0)+\ep\varphi_1f'(0)+o(\ep\varphi_1)
=\ep\varphi_1 f'(0) +o(\ep\varphi_1)\, .$$
So the inequality $-\Delta (\ep\varphi_1)\leq f(\ep\varphi_1)$
is equivalent to
$$\ep\lambda_1\varphi_1\leq
\ep\varphi_1 f'(0) +o(\ep\varphi_1)\, ,$$
that is
$$\lambda_1\leq f'(0)+o(1)\, .$$
This is true, by our hypothesis \eq{fprim}.

Let us now verify (ii). Recall that $\Usup$ satisfies
$$
\itab
$\di -\Delta \Usup =a\Usup +C\,,$ &\quad $\mbox{in}\ \Omega$\\
$\di \Usup =0\,,$ & \quad $\mbox{on}\ \partial\Omega\, .$\\
\ttab
$$
Thus, by the Maximum Principle, $\Usup >0$ in $\Omega$ and
$\frac{\partial \overline{U}}{\partial\nu}<0$ on $\partial\Omega$.
We observe that the other variant, namely $\Usup\equiv 0$, becomes
impossible, since $C>0$. Using the same trick as in the proof of
Theorem \ref{ssol} (more precisely, the fact that $\frac{\partial
\overline{U}}{\partial\nu}<0$ on $\partial\Omega$) we find $\ep
>0$ small enough so that $\ep\varphi_1\leq\Usup$ in $\Omega$. \qed

\begin{rem}\label{exist}
We observe that a necessary condition for the existence of a
solution to the problem \eq{P3} is that the line $\lambda_1u$
intersects the graph of the function $f=f(u)$ on the positive
semi-axis.
\end{rem}
Indeed, if $f(u)<\lambda_1u$ for any $u>0$
then the unique solution is $u=0$. After multiplication
with $\varphi_1$ in \eq{P3} and integration we find
$$
\intom (-\Delta u)\varphi_1=-\intom u\Delta\varphi_1=
\lambda_1\intom u\varphi_1= \intom f(u)\varphi_1<\lambda_1\intom
u\varphi_1\, ,
$$
a contradiction. \qed

\begin{rem}\label{ghergu}
We can require instead of \eq{fprim} that $f\in C^1(0,\infty)$ and
$f'(0+)=+\infty$.
\end{rem}
Indeed, since $f(0)=0$, there exists $c>0$ such that, for all
$0<\ep<c$,
$$\frac{f(\ep\varphi_1)}{\ep\varphi_1}=\frac{f(\ep\varphi_1)-f(0)}{\ep\varphi_1-0}>\lambda_1\,.$$
It follows that $f(\ep\varphi_1)>\ep\lambda_1\varphi_1=-\Delta
(\ep\varphi_1)$ and so, $\underline U=\ep\varphi_1$ is a
sub-solution. It is easy to check that $\overline U\geq \underline
U$ in $\Omega$ and the proof continues with the same ideas as
above. \qed

The following result gives a sufficient condition for
that the solution of \eq{P3} is unique.

\begin{teo}\label{unic}
 Under hypotheses \eq{x0}, \eq{lsup}, \eq{fprim}
 assume furthermore that
 \neweq{bo}
 \mbox{the mapping}\ \ (0,+\infty )\ni u\longmapsto
 \frac{f(u)}{u}\ \ \mbox{is decreasing}\, .
 \endeq
 Then the problem \eq{P3} has a unique solution.
 \end{teo}
{\bf Example}. If $f$ is concave then the mapping
$ \frac{f(u)}{u}$ is decreasing. Hence by the previous
results the solution is unique and stable. For instance
the problem
$$
\itab
$\di -\Delta u=\lambda u-u^p\, , $ &\quad $\mbox{in}\ \Omega $\\
$\di u=0\,,$ &\quad $\mbox{on}\ \partial\Omega$\\
$\di u>0\,,$ &\quad $\mbox{in}\ \Omega\, ,$\\
\ttab
$$
with $p>1$ and $\lambda >\lambda_1$ has a unique solution
which is also stable.\vspace*{0.2cm}\\
{\bf Proof of Theorem \ref{unic}.} Let $u_1, u_2$ be
arbitrary solutions of \eq{P3}. We may assume that
$u_1\leq u_2$; indeed, if not, we choose $u_1$ as the
minimal solution.
Multiplying the equalities
$$-\Delta u_1=f(u_1)\, ,\quad\mbox{in}\ \Omega$$
and
$$-\Delta u_2=f(u_2)\, ,\quad\mbox{in}\ \Omega$$
by $u_2$, resp. $u_1$, and integrating on $\Omega$
we find
$$\intom\left(f(u_1)u_2-f(u_2)u_1\right)=0\, ,$$
or, equivalently,
$$\intom u_1u_2\left( \frac{f(u_1)}{u_1}-
\frac{f(u_2)}{u_2}\right)=0\, .$$ So, by $0<u_1\leq u_2$  we
deduce $\frac{f(u_1)}{u_1}= \frac{f(u_2)}{u_2}$ in $\Omega$. Now,
by \eq{bo} we conclude that $u_1=u_2$. \qed

Let us now consider the problem
\neweq{P4}
\itab
$\di -\Delta u=f(u)\,,$ & \quad $\mbox{in}\ \Omega $\\
$\di u=0\,,$ & \quad $\mbox{on}\ \partial\Omega$\\
$\di u>0\,,$ &\quad $\mbox{in}\ \Omega\, .$\\
\ttab
\endeq
Our aim is to establish in what follows a corresponding result
in the case where the nonlinearity does not satisfy any
growth assumptions, like \eq{lsup} or \eq{fprim}. The
following deeper result in this direction is due to
Krasnoselski.

\begin{teo}\label{krasn}
Assume that the nonlinearity $f$ satisfies \eq{bo}. Then
the problem \eq{P4} has a unique solution.
\end{teo}
{\bf Proof.} In order to probe the uniqueness it is
enough to show that for any arbitrary solutions $u_1$ and
$u_2$, we
can suppose that $u_1\leq u_2$. Let
$$A =\{ t\in [0,1]; tu_1\leq u_2\}\, .$$
We observe that $0\in A$. Now we show that there exists
$\ep_0>0$ such that for every $\ep\in (0,\ep_0)$ we
have $\ep\in A$. This follows easily with arguments
which are already done and using the crucial observations
that
$$ \itab
$\di u_2>0\, ,$ &\quad $\mbox{in}\ \Omega$\\
$\di \frac{\partial u_2}{\partial\nu}<0\, ,$ & \quad $\mbox{on}
\ \partial\Omega\, .$\\
\ttab
$$
Let $t_0=\max A$. Assume, by contradiction, that
$t_0<1$. Hence $t_0u_1\leq u_2$ in $\Omega$. The idea
is to show the existence of some $\ep >0$ such that
$(t_0+\ep )u_1\leq u_2$, which contradicts the choice
of $t_0$. For this aim we use the Maximum Principle.
We have
$$-\Delta (u_2-t_0u_1)+a(u_2-t_0u_1)=f(u_2)+au_2
-t_0\left( f(u_1)+au_1\right)\, .$$
Now we choose $a>0$ so that the mapping $u\longmapsto
f(u)+au$ is increasing. Therefore
$$
\begin{array}{lll}
-\Delta (u_2-t_0u_1)+a(u_2-t_0u_1)\geq\\
f(t_0u_1)+at_0u_1-t_0f(u_1)-at_0u_1=f(t_0u_1)-t_0f(u_1)
\geq 0\, .
\end{array}$$
This implies either

\noindent (i) $u_2-t_0u_1\equiv 0$

\noindent or

\noindent (ii) $u_2-t_0u_1>0$ in $\Omega$ and
$\frac{\partial }{\partial\nu}\, (u_2-t_0u_1)<0$
on $\partial\Omega$.

The first case is impossible since it would imply
$t_0f(u_1)=f(t_0u_1)$, a contradiction. This reasoning is based on
the elementary fact that if $f$ is continuous and $f(Cx)=Cf(x)$
for any  $x$ in a nonempty interval then $f$ is linear. \qed

\chapter{Implicit Function Theorem and Applications to
Boundary Value Problems}

\section{Abstract theorems}
Let $X,Y$ be  Banach spaces. Our aim is to develop a general
method which will enable us to solve equations of the type
$$F(u,\lambda )=v\, ,$$
where $F:X\times\RR\ri Y$ is a prescribed sufficiently smooth
function and $v\in Y$ is given.

\begin{teo}\label{IFT}
Let $X,Y$ be real Banach spaces and let $(u_0,\lambda_0)
\in X\times\RR$. Consider a $C^1$-mapping $F=F(u,\lambda ):
X\times\RR\ri Y$ such that the following conditions hold:

(i) $F(u_0,\lambda_0)=0$;

(ii) the linear mapping $F_u(u_0,\lambda_0):X\ri Y$ is bijective.

Then there exists a neighbourhood $U_0$ of $u_0$ and a
neighbourhood $V_0$ of $\lambda_0$ such that for every
$\lambda\in V_0$ there is a unique element $u(\lambda )
\in U_0$ so that $F(u(\lambda ),\lambda )=0$.

Moreover, the mapping $V_0\ni\lambda\longmapsto u(\lambda )$
is of class $C^1$.
\end{teo}

{\bf Proof}. Consider the mapping $\Phi (u,\lambda ):
X\times\RR\ri Y\times\RR$ defined by
$\Phi (u,\lambda )=\left( F(u,\lambda ),\lambda\right)$.
It is obvious that $\Phi\in C^1$. We apply to $\Phi$
the Inverse Function Theorem. For this aim, it remains to
verify that the mapping $\Phi '(u_0,\lambda_0):
X\times\RR\ri Y\times\RR$ is bijective. Indeed, we have
$$\begin{array}{lll}
\Phi (u_0+tu,\lambda_0+t\lambda )=
\left( F(u_0+tu,\lambda_0+t\lambda ),
\lambda_0+t\lambda\right)= \\
\left( F(u_0,\lambda_0)+F_u(u_0,\lambda_0)\cdot (tu)+
F_\lambda (u_0,\lambda_0)\cdot (t\lambda )+o(1),\lambda_0+
t\lambda\right)\, .
\end{array}
$$
It follows that
$$F'(u_0,\lambda_0)= \left(
\begin{array}{cc}
F_u(u_0,\lambda_0) &  F_\lambda (u_0,\lambda_0) \\
       0            &           \mbox{I}
       \end{array}
       \right)
$$
which is a bijective operator, by our hypotheses. Thus, by the
Inverse Function Theorem, there exist a neighbourhood ${\cal U}$
of the point $(u_0,\lambda_0)$ and a neighbourhood ${\cal V}$ of
$(0,\lambda )$ such that the equation
$$\Phi (u,\lambda )=(f,\lambda_0 )$$
has a unique solution, for every $(f,\lambda )\in {\cal V}.$ Now
it is sufficient to take here $f=0$ and our conclusion follows.
\qed

 With a similar proof one can justify the following {\bf global}
version of the Implicit Function Theorem.

\begin{teo}\label{globally}
Assume $F:X\times\RR\ri Y$ is a $C^1$-function on $X\times\RR$
satisfying

(i) $F(0,0)=0$;

(ii) the linear mapping $F_u(0,0):X\ri Y$ is bijective.

Then there exist an open neighbourhood $I$ of 0 and a $C^1$
mapping $I\ni\lambda\longmapsto u(\lambda )$ such that
$u(0)=0$ and $F(u(\lambda ),\lambda )=0$.
\end{teo}

The following result will be of particular importance in the
next applications.

\begin{teo}\label{teobif}
Assume the same hypotheses on $F$ as in Theorem \ref{globally}.
Then there exists an open and maximal interval $I$ containing
the origin and there exists a unique $C^1$-mapping
$I\ni\lambda\longmapsto u(\lambda )$ such that the
following hold:

a) $F(u(\lambda ),\lambda )=0$, for every $\lambda\in I$;

b) the linear mapping $F_u(u(\lambda ),\lambda )$ is
bijective, for any $\lambda\in I$;

c) $u(0)=0$.
\end{teo}

{\bf Proof}. Let $u_1$, $u_2$ be solutions and consider the
corresponding open intervals $I_1$ and $I_2$ on which these
solutions exist, respectively. It follows that $u_1(0)=
u_2(0)=0$ and
$$F(u_1(\lambda ),\lambda )=0\, ,\quad\mbox
{for every $\lambda\in I_1$}\, ,$$
$$F(u_2(\lambda ),\lambda )=0\, ,\quad\mbox
{for every $\lambda\in I_2$}\, .$$
Moreover, the mappings $F_u(u_1(\lambda ),\lambda )$ and
$F_u(u_2(\lambda ),\lambda )$ are one-to-one and onto on
$I_1$, resp. $I_2$. But, for $\lambda$ sufficiently close
to 0 we have $u_1(\lambda )=u_2(\lambda )$. We wish to
show that we have global uniqueness. For this aim, let
$$I=\{\lambda\in I_1\cap I_2;\ u_1(\lambda )=u_2(\lambda)\}\, .$$
Our aim is to show that $I=I_1\cap I_2$. We first observe that
$0\in I$, so $I\not=\emptyset$. A standard argument then shows
that $I$ is closed in $I_1\cap I_2$. In order to show that
$I=I_1\cap I_2$, it is sufficient now to prove that $I$ is an open
set in $I_1\cap I_2$. The proof of this statement follows by
applying Theorem \ref{IFT} for $\lambda$ instead of 0. Thus,
$I=I_1\cap I_2$.

Now, in order to justify the existence of a maximal interval $I$,
we consider the $C^1$-curves $u_n(\lambda )$ defined on the
corresponding open intervals $I_n$, such that $0\in I_n$,
$u_n(0)=u_0$, $F(u_n(\lambda ), \lambda )=0$ and $F_u(u_n(\lambda
),\lambda )$ is an isomorphism, for any $\lambda\in I_n$. Now a
standard argument enables us to construct a maximal solution on
the set $\cup_n I_n$. \qed

\begin{cor}\label{surj}
Let $X,Y$ be Banach spaces and let $F:X\ri Y$ be a $C^1$-function.
Assume that the linear mapping $F_u(u):X\ri Y$ is bijective, for every
$u\in X$ and there exists $C>0$ such that $\| \left(
F_u(u)\right)^{-1}\|\leq C$, for any $u\in X$.

Then $F$ is onto.
\end{cor}
{\bf Proof}. Assume, without loss of generality, that $F(0)=0$ and
fix arbitrarily $f\in Y$. Consider the operator $F(u,\lambda
)=F(u)-\lambda f$, defined on $X\times\RR$. Then, by Theorem
\ref{teobif}, there exists a $C^1$-function $u(\lambda )$ which is
defined on a maximal interval $I$ such that $F(u(\lambda
))=\lambda f$. In particular, $u:=u(1)$ is a solution of the
equation $F(u)=f$. We assert that $I=\RR$. Indeed, we have
$$u_\lambda (\lambda )= \left(
F_u(u)\right)^{-1} f\, ,$$ so $u$ is a Lipschitz map on $I$, which
implies $I=\RR$.

The Implicit Function Theorem is used to solve equations
of the type $F(u)=f$, where $F\in C^1(X,Y).$ A simple
method for proving that $F$ is onto, that is,
Im$\, F\, =\, Y$ is to prove the following:

(i) Im$\, F$ is open;

(ii) Im$\, F$ is closed.

For showing (i), usually we use the Inverse Function
Theorem, more exactly, if $F_u(u)$ is one-to-one,
for every $u\in X$, then (i) holds. A sufficient condition
for that (ii) holds is that $F$ is a proper map.

Another variant of the Implicit Function Theorem is
given by

\begin{teo}\label{variant}
Let $F(u,\lambda )$ be a $C^1$-mapping in a neighbourhood
of $(0,0)$ and such that $F(0,0)=0$. assume that

(i) {\rm Im}$\, F_u(0,0)=Y$;

(ii) the space $X_1:={\rm Ker}\, F_u(0,0)$ has a
closed complement $X_2$.

Then there exist $B_1=\{ u_1\in X;\| u_1\|<\delta\}$,
$B_2=\{ \lambda\in\RR ;\mid\lambda\mid <r\}$,
$B_3=\{ g\in Y;\| g\| <R\}$ and a neighbourhood
$U$ of the origin in $X_2$ such that, for any $u_1\in B_1$,
$\lambda\in B_2$ and $g\in B_3$, there exists a unique
 solution $u_2=\varphi (u_1,\lambda ,g)\in U$ of the equation
 $$F(u_1+\varphi (u_1,\lambda ,g),\lambda )=g\, .$$
 \end{teo}

 {\bf Proof}. Let $\Gamma =X\times\RR\times Y$, that is,
 every element $\nu\in\Gamma$ has the form
 $\nu =(u_1,\lambda ,g)$. It remains to apply then
 Implicit Function Theorem to the mapping
 $G:X\times\Gamma\ri Y$ which is defined by
 $G(u_2,\nu )=F(u_1+u_2,\lambda )-g\, .$
 \qed

 \medskip
We conclude this paragraph with the following elementary
example: let $\Omega$ be a smooth bounded domain in
$\RR^N$ and let $g$ be a $C^1$ real function defined on
a neighbourhood of 0 and such that $g(0)=0$. Consider
the problem
\neweq{exeift}
\itab
$-\Delta u=g(u)+f(x)\, ,$ &\quad $\mbox{in}\ \Omega$\\
$ u=0\, ,$ & \quad $\mbox{on}\ \partial\Omega\, .$\\
\ttab
\endeq
Assume that $g'(0)$ is not a real number of $-\Delta$
in $H^1_0(\Omega )$, say $g'(0)\leq 0$. If $f$ is
sufficiently small, then the problem \eq{exeift} has
a unique solution, by the  Implicit Function
Theorem. Indeed,
it is enough to apply Theorem \ref {IFT} to the operator
$F(u)=-\Delta u-g(u)$ and after observing that
$F_u(0)=-\Delta -g'(0)$. There are at least two distinct
possibilities for defining $F$:

\noindent (i) $F:C_0^{2,\alpha}(\overline\Omega )\ri
C_0^{\alpha}(\overline\Omega )$, for some $\alpha\in
(0,1)$.

\noindent or

\noindent (ii) $F:W^{2,p}(\Omega )\cap W_0^{1,p}
(\Omega )\ri L^p(\Omega )$. In order to obtain classical
solutions (by a standard bootstrap argument that we will describe
later), it is sufficient to choose $p>\frac{N}{2}$.

\section{A basic bifurcation theorem}
  Consider a $C^2$ map $f:\RR\rightarrow\RR$ which is
  convex, positive and such that
  $f'(0)>0$.
  Our aim is to study the problem
\neweq{P}
\itab
$ -\Delta u=\lambda f(u),$ & \quad $\hbox{in}\ \ \Omega $\\
$u=0,$ &\quad $\hbox{on}\ \ \partial\Omega\, ,$\\
\ttab
\endeq
where $\lambda$ is a positive parameter. We are looking for classical solutions of this
problem, that is,
$u\in C^2(\Omega )\cap C(\overline{\Omega})$.

\smallskip
Trying to apply the Implicit Function Theorem to our problem
\eq{P}, set
$$X=\{ u\in C^{2,\alpha}(\overline\Omega );u=0\ \mbox{on}\ \partial\Omega\}$$
and $Y=C^{0,\alpha}(\overline\Omega )$, for some $0<\alpha <1$. Define $F(u,\lambda
)=-\Delta u-\lambda f(u)$. It is clear that $F$ verifies all the assumptions of the
Implicit Function Theorem. Hence, there exist a {\rm maximal} neighbourhood of the
origin $I$ and a {\rm unique} map $u=u(\lambda )$ which is solution of the problem (P)
and such that the linearized operator $-\Delta -\lambda f'(u(\lambda ))$ is bijective.
In other words, for every $\lambda\in I$, the problem \eq{P} admits a stable solution which
is given by the Implicit Function Theorem. Let $\lambda^\star :=\sup I\leq +\infty$. We
shall denote from now on by $\lambda_1(-\Delta -a)$ the first eigenvalue in
$H^1_0(\Omega )$ of the operator $-\Delta -a$, where $a\in L^\infty (\Omega )$.

Our aim is to prove in what follows the following celebrated
result.

\medskip
{\bf Amann's Theorem}. {\sl Assume that $f:\RR\ri\RR$
 is a $C^2$ which is
  convex, positive and such that
  $f'(0)>0$. Then the following hold:

  i) $\lambda^\star <+\infty$;

  ii) $\lambda_1(-\Delta -\lambda f'(u(\lambda )))>0$;

  iii) the mapping $I\ni\lambda \longmapsto u(x,\lambda )$ is increasing, for every
  $x\in \Omega$;

iv) for every $\lambda\in I$ and $x\in\Omega$, we have $u(x,\lambda )>0$;

v) there is no solution of the problem \eq{P}, provided that $\lambda >\lambda^\star$;

vi) $u(\lambda )$ is a minimal solution of the problem \eq{P};

vii) $u(\lambda )$ is the unique stable solution of the problem
\eq{P}}.\rm

\medskip
{\bf Proof}. i) It is obvious, by the variational characterization of the first
eigenvalue, that if $a,b\in L^\infty (\Omega )$ and $a\leq b$  then
$$\lambda_1(-\Delta -a(x))\geq\lambda_1(-\Delta -b(x))\, .$$
 Assume iv) is already proved.
Thus, by the convexity of $f$, $f'(u(\lambda ))\geq f'(0)>0$, which implies that
$$\lambda_1(-\Delta -\lambda f'(0))\geq\lambda_1(-\Delta -\lambda f'(u(\lambda )))>0,$$
for every $\lambda\in I$. This implies $\lambda_1-\lambda f'(0)>0$, for any $\lambda
<\lambda^\star$, that is, $\lambda^\star\leq\lambda_1/f'(0)<+\infty$.

ii) Set $\varphi (\lambda )=\lambda_1(-\Delta -\lambda f'(u(\lambda )))$. So, $\varphi
(0)=\lambda_1>0$ and, for every $\lambda <\lambda^\star$, $\varphi (\lambda )\not= 0$,
since the linearized operator $-\Delta -\lambda f'(u(\lambda ))$ is bijective, by the
Implicit Function Theorem. Now, by the continuity of the mapping $\lambda\longmapsto
\lambda f'(u(\lambda ))$, it follows that $\varphi$ is continuous, which implies, by the
above remarks, $\varphi >0$ on $[0,\lambda^\star )$.

iii) We differentiate in \eq{P} with respect to $\lambda$. Thus
$$-\Delta u_\lambda =f(u(\lambda ))+\lambda f'(u(\lambda ))\cdot
u_\lambda\qquad\mbox{in}\ \ \Omega$$
and $u_\lambda =0$ on $\partial\Omega$. Hence
$$(-\Delta -\lambda f'(u(\lambda )))u_\lambda =f(u(\lambda ))\qquad\mbox{in}\ \ \Omega
.$$
But the operator $(-\Delta -\lambda f'(u(\lambda )))$ is coercive. So, by Stampacchia's
Maximum Principle, either $u_\lambda \equiv 0$ in $\Omega$, or $
u_\lambda > 0$ in $\Omega$. The first variant is not convenient, since it would imply
that $f(u(\lambda ))=0$, which is impossible, by our initial hypotheses. It remains that
$u_\lambda > 0$ in $\Omega$.

iv) follows from iii).

v) Assume that there exists some $\nu >\lambda^\star$ and there exists a corresponding
solution $v$ to our problem \eq{P}.

\medskip
{\bf Claim}. $u(\lambda )<v$ in $\Omega$, for every $\lambda <\lambda^\star$.

\medskip
{\bf Proof of the claim}. By the convexity of $f$ it follows that
$$\begin{array}{cl}
& -\Delta (v-u(\lambda ))=\nu f(v)-\lambda f(u(\lambda ))\geq\\
& \lambda (f(v)-f(u(\lambda )))\geq\lambda f'(u(\lambda ))(v-u(\lambda )).
\end{array}
$$
Hence
$$\begin{array}{cl}
& -\Delta (v-u(\lambda ))-\lambda f'(u(\lambda ))(v
-u(\lambda ))=\\
&(-\Delta -\lambda f'(u(\lambda )))
(v-u(\lambda ))\geq 0\quad\mbox{in}\ \ \Omega ,
\end{array}
$$
since the operator $-\Delta -\lambda f'(u(\lambda ))$ is coercive. Thus, by
Stampacchia's Maximum Principle, $v\geq u(\lambda )$ in $\Omega$, for every $\lambda
<\lambda^\star$.

\medskip
Hence, $u(\lambda )$ is bounded in $L^\infty$ by $v$. Passing to the limit as
$\lambda\ri\lambda^\star$ we find
that $u(\lambda )\ri u^\star<+\infty$ and $u^\star =0$ on $\partial\Omega$.

We prove in what follows that
$$\lambda_1(-\Delta -\lambda^\star f'(u^\star )))=0.$$
 We already know that
$\lambda_1(-\Delta -\lambda^\star f'(u^\star )))\geq 0$. Assume
that $\lambda_1(-\Delta -\lambda^\star f'(u^\star )))> 0$, so,
this operator is coercive. We apply the Implicit Function Theorem
to $F(u,\lambda )=-\Delta u-\lambda f(u)$ at the point $(u^\star
,\lambda^\star )$. We obtain that there is a curve of solutions of
the problem \eq{P} passing through $(u^\star ,\lambda^\star )$,
which contradicts the maximality of $\lambda^\star$. We have
obtained that $\lambda_1(-\Delta -\lambda^\star f'(u^\star )))=0$.
So, there exists $\varphi_1>0$ in $\Omega$, $\varphi_1=0$ on
$\partial\Omega$, so that
\neweq{mgherg}-\Delta\varphi_1-\lambda^\star f'(u^\star
)\varphi_1=0\qquad\mbox{ in}\ \, \Omega .\endeq Passing to the
limit as $\lambda\ri\lambda^\star$ in the relation
$$(-\Delta -\lambda f'(u(\lambda )))(v-u(\lambda ))\geq 0$$
we find
$$(-\Delta -\lambda^\star f'(u^\star )))(v-u^\star)\geq 0.$$
Multiplying this inequality by $\varphi_1$ and integrating, we obtain
$$-\intom (v-u^\star )\Delta\varphi_1dx-\lambda^\star\intom f'(u^\star)(v-u^\star
)\varphi_1dx\geq 0.$$ In fact, by \eq{mgherg}, the above relation
is an equality, which implies that
$$-\Delta (v-u^\star )=\lambda^\star f'(u^\star )(v-u^\star )\qquad\mbox{in}\ \,\Omega
.$$ It follows that $\nu f(v)=\lambda^\star f(u^\star)$ in
$\Omega$. But $\nu >\lambda^\star$ and $f(v)\geq f(u^\star )$. So,
$f(v)=0$ which is impossible.

vi) Fix an arbitrary $\lambda <\lambda^\star$.
Assume that $v$ is another solution of the problem \eq{P}. We have
$$-\Delta (v-u(\lambda ))=\lambda f(v)-\lambda f(u(\lambda ))\geq\lambda f'(u(\lambda
))(v-u(\lambda))\quad\mbox{in}\ \Omega .$$
Again, by Stampacchia's Maximum Principle applied to the coercive operator
$-\Delta -\lambda f'(u(\lambda ))$ we find that $v\geq u(\lambda )$ in $\Omega$.

vii) Let $v$ be another stable solution, for some $\lambda
<\lambda^\star$. With the same reasoning as in vi), but applied to
the coercive operator $-\Delta -\lambda f'(v)$, we get that
$u(\lambda )\geq v$. Finally, $u(\lambda )=v$. \qed

\section{Qualitative properties of the minimal solution in a neighbourhood of
the bifurcation point}

We assume, throughout this section, that
\neweq{4.1}
\lim_{t\rightarrow\infty}{{f(t)}\over t}=a>0\, .
\endeq

We impose in this section a supplementary hypothesis to
 $f$, which essentially means that as
 $t\rightarrow\infty$, then $f$ increases faster than $t^a$, for some $a>1$.
In this case we shall prove that the problem \eq{P} has a weak solution
$u^*\in\huo$ provided
$\lambda =\lambda^*$. However, we can not obtain in this case a supplementary regularity
of $u^*$.

\medskip
\begin{teo}\label
{4.1.1}  Assume that $f$ satisfies the following additional condition:
there is some $a>0$ and $\mu >1$ such that, for every $t\geq a$,
$$tf'(t)\geq\mu f(t).$$
The following hold

i) the problem {P} has a solution
 $u^*$, provided that $\lambda =\lambda^*$;

ii) $u^*$ is the weak limit in \huo
of stable solutions $u(\lambda )$,
if $\lambda\nearrow\lambda^*$;

iii) for every $v\in\huo$,
$$f'(u^*)v^2\in L^1(\Omega )$$
and
$$\lambda^*\intom f'(u^*)v^2dx\leq\intom |\nabla v|^2dx\, .$$
\end{teo}

\medskip
{\bf Proof.} For every $v\in\huo$ and $\lambda\in [0,\lambda^* )$
we have, by Amann's Theorem,
\neweq{PP}
\lambda^*\intom f'(u(\lambda ))v^2dx\leq\intom
|\nabla v|^2dx\, .
\endeq

Choosing here $v=u(\lambda )$, we find
$$\lambda^*\intom f'(u(\lambda ))u^2(\lambda )dx\leq\intom |\nabla u(\lambda )|^2dx=$$
$$=\lambda\intom f(u(\lambda ))u(\lambda )dx\, .$$
So, if for $a$ as in our hypothesis we define
$$\Omega (\lambda )=\{ x\in\Omega ;\ u(\lambda )(x)>a\} ,$$
then
$$\intom f'(u(\lambda ))u^2(\lambda )dx
\leq \frac{\lambda}{\lambda^*}\int_\Omega f(u(\lambda ))u(\lambda
)dx\leq \frac{\lambda}{\lambda^*}\, \left [\ \int_{\Omega
 (\lambda )}f(u(\lambda ))u(\lambda )dx
+a\cdot |\Omega |\cdot f(a)\ \right ].$$ By our hypotheses on $f$,
we have
$$f'(u(\lambda ))u^2(\lambda )\geq\mu f(u(\lambda ))u(\lambda )\quad\mbox{in}\ \
\Omega (\lambda )\, .
$$
Therefore
\neweq{4.4}
(\mu -1)\int_{\Omega (\lambda )}f(u(\lambda ))u(\lambda )
dx\leq C,
\endeq
where the constant $C$ depends only on
 $\lambda$. This constant can be chosen sufficiently large so that
\neweq{4.5}
\int_{\Omega\setminus\Omega (\lambda )}f(u(\lambda ))u(\lambda )dx\leq C
\, .
\endeq
By \eq{4.4}, \eq{4.5} and $\mu >1$ it follows that there exists $C>0$ independent of
$\lambda$ such that
\neweq{4.6}
\intom f(u(\lambda ))u(\lambda )dx\leq C\, .
\endeq
From here and from
$$\intom |\nabla u(\lambda )|^2dx=\lambda\intom f(u(\lambda ))u(\lambda )dx,$$
it follows that $u(\lambda )$ is bounded in \huo, independently
with respect to $\lambda$. Consequently, up to a subsequence, we may suppose that there
exists  $u^*\in\huo$ such that
\neweq{4.7}
u(\lambda )\rightharpoonup u^*\quad\hbox{weakly in}\ \ \huo\ \ \hbox{if}\
\ \lambda\rightarrow\lambda^*
\endeq
\neweq{4.8}
u(\lambda )\rightarrow u^*\quad\hbox{a.e. in}\ \ \Omega\
\ \hbox{if}\ \ \lambda\rightarrow\lambda^*\, .
\endeq
Hence
$$f(u(\lambda ))\rightarrow f(u^*)\quad\hbox{a.e. in}\ \
 \Omega\ \ \hbox{if}\ \ \lambda
\rightarrow\lambda^*\, .$$
By \eq{4.6} we get
$$\intom f(u(\lambda ))dx\leq C\, .$$

Since the mapping $\lambda\longmapsto f(u(\lambda ))$ is increasing and the integral is
bounded, we find, by the Monotone Convergence Theorem,
$f(u^*)\in L^1(\Omega )$ and
$$f(u(\lambda ))\rightarrow f(u^*)\ \ \ \hbox{in}\
\ L^1(\Omega )\quad\hbox{if}\ \ \lambda\rightarrow\lambda^*\, .$$

Let us now choose $v\in\huo$, $v\geq 0$. So,
$$\intom \nabla u(\lambda )\cdot\nabla vdx=\lambda\intom
f(u(\lambda ))vdx\, .$$
On the other hand, we have already remarked that
$$f(u(\lambda ))v\rightarrow f(u^*)v\quad\hbox{a.e. in}\ \ \Omega ,\ \ \hbox{if}\
\lambda\rightarrow \lambda^*\, .$$
and
$$\lambda\longmapsto f(u(\lambda ))v\quad\hbox{is increasing}\, .$$

By \eq{4.7}, it follows that
$$\intom f(u(\lambda ))v\leq C\, .$$
Now, again by the Monotone Convergence Theorem,
$$f(u(\lambda ))v\rightarrow f(u^*)v\quad\hbox{in}\ \ L^1(\Omega )\ \ \hbox{if}\
 \ \lambda\rightarrow\lambda^*$$
and
$$\intom \nabla u^*\cdot\nabla vdx=\lambda^*\intom f(u^*)vdx\, .$$

If $v\in\huo$ is arbitrary, we find the same conclusion if we consider
$v=v^+-v^-$.

So, for every $v\in\huo$,
$$f(u(\lambda ))v\rightarrow f(u^*)v\quad\hbox{in}\ \ L^1(\Omega )$$
$$\intom \nabla u^*\cdot\nabla vdx=\lambda^*\intom f(u^*)vdx$$
and $u^*\in\huo$. Consequently, $u^*$ is a weak solution
of the problem \eq{P}. Moreover, for every $v\in\huo$,
$$\lambda^*\intom f'(u(\lambda ))v^2dx\leq \intom |\nabla v|^2dx\, .$$

Applying again the Monotone Convergence Theorem we find that
$$f'(u^*)v^2\in L^1(\Omega )$$
and
$$\lambda^*\intom f'(u^*)v^2dx\leq\intom |\nabla v|^2dx\, .$$
\qed

\medskip
The facts that
$$f(u^*)\in L^1(\Omega ),\ \ f(u^*)u^*\in L^1(\Omega )$$
and
$$f'(u^*)v^2\quad\hbox{in}\quad L^1(\Omega ),\quad\hbox{ for every}
 \ \ v\in\huo ,$$
imply a supplementary regularity of $u^*$. In many concrete situations one may show
that there exists
 $n_0$ such that, if $N\leq n_0$,
then $u^*\in L^\infty (\Omega )$ and $f(u^*)\in L^\infty (\Omega )$.

We shall deduce in what follows some known results, for the special case
$f(t)=e^t$.

\medskip
\begin{teo}\label
{4.1.2}  Let $f(t)=e^t$. Then
$$f(u(\lambda ))\rightarrow f(u^*)\quad\hbox{in}\ \
L^p(\Omega ),\ \ \hbox{if}\ \ \lambda\rightarrow\lambda^* ,$$
for every $p\in [1,5)$.

Consequently
$$u(\lambda )\rightarrow u^*\quad\hbox{in}\ \ W^{2,p}(\Omega ),\ \ \hbox{if}
\ \ \lambda\rightarrow\lambda^* ,$$
for every $p\in [1,5)$.

Moreover, if $N\leq 9$, then
$$u^*\in L^\infty (\Omega )\ \ \ \hbox{and}\ \ \
 f(u^*)\in L^\infty (\Omega )\, .$$
 \end{teo}

\medskip
{\bf Proof.} As we have already remarked, for every $v\in\huo$,
\neweq{4.10}
\intom \nabla u(\lambda )\cdot\nabla vdx=\lambda\intom
e^{u(\lambda )}vdx\endeq
\neweq{4.11}
\lambda\intom e^{u(\lambda )}v^2dx\leq\intom
|\nabla v|^2dx\, .
\endeq

In \eq{4.10} we choose $v=e^{(p-1)u(\lambda )}-1\in\huo$, for $p>1$ arbitrary. We find
\neweq{4.12}
(p-1)\intom e^{(p-1)u(\lambda )}|\nabla u(\lambda )|^2dx=\lambda\intom e^{u(\lambda )}
[e^{(p-1)u(\lambda )}-1]dx\, .
\endeq

In \eq{4.11} we put $$v=e^{{{p-1}\over\scriptstyle 2}u(\lambda )}
-1\in\huo\, .$$
Hence
\neweq{4.13}
\left\{\begin{array}{lll}
\lambda\intom e^{u(\lambda )}[e^{{{p-1}\over 2}u(\lambda
)}-1]^2dx\leq  \\
\leq {{(p-1)^2}\over 4}\intom e^{(p-1)u(\lambda )}
|\nabla u(\lambda )|^2dx\, .
\end{array}\right.
\endeq

Taking into account the relation \eq{4.12}, our relation
\eq{4.13} becomes
\neweq{4.14}
\left\{\begin{array}{lll}
\lambda\intom e^{pu(\lambda)}dx-2\lambda\intom e^{{{p+1}\over 2}
u(\lambda )}dx+\lambda\intom e^{u(\lambda )}dx
\leq \\
\leq {{p-1}\over 4}[\lambda\intom e^{pu(\lambda )}dx-\lambda\intom
 e^{u(\lambda )}dx]\, .
\end{array}\right.
\endeq

By H\"older's Inequality, our relation \eq{4.14} yields
$$\lambda \left(1-{{p-1}\over 4}\right)\intom e^{pu(\lambda )}dx+
\lambda \left(1+{{p-1}\over 4}\right) \intom e^{u(\lambda )}dx
\leq 2\lambda\intom e^{{{p+1}\over 2}u(\lambda )}dx \leq
C\left(\intom e^{pu(\lambda )}dx\right)^2,$$ where $C$ is a
constant which does not depend on $\lambda$.

So, if $1-{{p-1}\over 4}>0$, that is $p<5$, then the mapping $e^{u(\lambda )}$
is bounded  in $L^p(\Omega )$.

We have already proved that if $\lambda\rightarrow\lambda^*$, then
$$e^{u(\lambda )}\rightarrow e^{u^*}\quad\hbox{a.e. in}\ \ \Omega\, .$$
Moreover
$$e^{u(\lambda )}\leq e^{u^*}\quad\hbox{in}\ \ \Omega\, .$$

By the Dominated Convergence Theorem, we find that, for every
$1\leq p<5$,
$$e^{u^*}\in L^p(\Omega )$$
and
$$e^{u(\lambda )}\rightarrow e^{u^*}\quad\hbox{in}\ \ L^p(\Omega )$$

Taking into account the relation \eq{P} which is fulfilled by $u(\lambda )$, and
using a standard regularity theorem for elliptic equations,
it follows that
$$u(\lambda )\rightarrow u^*\quad\hbox{in}\ \
W^{2,p}(\Omega )\ \ \hbox{if}\ \ \lambda\rightarrow\lambda^*\, .$$
On the other hand, by Sobolev inclusions,
$$W^{2,p}\subset L^\infty (\Omega ),\quad\hbox{if}\ \ N<2p\, .$$
So, if $N\leq 9$,
$$u(\lambda )\rightarrow u^*\quad\hbox{in}\ \ L^\infty (\Omega ),\ \ \hbox{if}
\ \ \lambda\rightarrow\lambda^*\, .$$
Hence
$$u^*\in L^\infty (\Omega )\quad\hbox{and}\quad e^{u^*}\in
L^\infty (\Omega )\, .$$ \qed

\medskip
\begin{rem}\label
{4.1.3} The above result is optimal, since D. Joseph and
T. Lundgren showed in \cite{JL} that if $N=10$ and
 $\Omega $ is an open ball, then $u^*\notin L^\infty (\Omega )$.
 \end{rem}

\section{Bifurcation problems associated to convex positive asymptotic
linear functions}

Throughout this section we assume that $f$ satisfies the same general assumptions as
in the preceding ones and, moreover, $f$ is asymptotic linear at infinity, in the
sense that
\neweq{4.15}
\lim_{t\rightarrow\infty}{{f(t)}\over t}=a\in
(0,+\infty )\, .
\endeq

Under these hypotheses, we shall study in detail the problem
\eq{P}, more exactly, we
shall try to ask the following questions:

i) what happens if $\lambda =\lambda^*$?

ii) what about the behaviour of the minimal solution $u(\lambda )$
 provided $\lambda$
approaches $\lambda^*$?

iii) the existence of other solutions except the minimal one
and, in this case,
 their behaviour.

First of all we make the following notations:

a) if $\alpha\in L^\infty (\Omega )$, let $\lambda_j(\alpha )$ and
$\varphi_j(\alpha )$ be the $j$th-eigenvalue (resp., eigenfunction) of the operator
$-\Delta -\alpha$. Moreover, we can assume that $\varphi_1(\alpha )>0$
and
$$\intom \varphi_j(\alpha )\varphi_k(\alpha )dx=\delta_{jk}\, .$$

If $\alpha =0$ we denote $\lambda_j$ (resp., $\varphi_j$).

b) a solution $u$ of the problem (P) is said to be {\it stable} if
$\lambda_1(\lambda f'(u))>0$ and {\it unstable}, in the opposite case.

\medskip

 \begin{teo}\label
 {4.2.1}  Under the preceding hypotheses, we also assume that
$$\lim_{t\rightarrow\infty}(f(t)-at)=l\geq 0\, .$$
Then

i)$\ \ \ \lambda^*={{\lambda_1}\over a}$.

ii)$\ \ \displaystyle\lim_{\lambda\nearrow\lambda^*}u(\lambda )
=\infty$, uniformly on compact subsets of $\Omega$.

iii) for every $\lambda\in (0,\lambda^*)$, the problem \eq{P}
 has only the minimal
solution.

iv) the problem \eq{P} has no solution if $\lambda =\lambda^*$.\rm
\end{teo}

\medskip
In the proof of this theorem we shall make use of the following auxiliary results:

\medskip
\begin{lemma}\label{4.2.2}  Let $\alpha\in L^\infty (\Omega )$,
 $w\in\huo\setminus\{ 0\},\ w\geq 0$,
such that $\lambda_1(\alpha )\leq 0$ and
\neweq{4.16}
-\Delta w\geq\alpha w\, ,\endeq

Then

i)$\ \ \ \lambda_1(\alpha )=0$.

ii) $\ \ -\Delta w=\alpha w$.

iii)$\ \ w>0$ in $\Omega$.\rm
\end{lemma}

{\bf Proof of the lemma.}
Multiplying \eq{4.16} by $\varphi_1 (\alpha )$ and integrating by parts,
we find
$$\intom\alpha\varphi_1(\alpha )w+\lambda_1(\alpha )\intom\varphi_1(\alpha )w\geq
\intom\alpha\varphi_1(\alpha )w\, .$$

Since $\lambda_1(\alpha )\leq 0$, it follows that
$\lambda_1(\alpha )=0$ and $-\Delta w=\alpha w$. But $w\geq 0$ and
$w\not\equiv 0$. Therefore there exists $C>0$ such that
$w=C\varphi_1(\alpha )$, that is $w>0$ in $\Omega$. \qed

\medskip
\begin{lemma}\label{4.2.3}  If there exist $a,b>0$ such that $f(t)=at+b$ for all
$t\geq 0$, then

i)$\ \ \lambda^*={{\lambda_1}\over\scriptstyle a}\, .$

ii) the problem \eq{P} has no solution if $\lambda =\lambda^*$.
\end{lemma}

\medskip
{\bf Proof of the lemma.} For every $0<\lambda <
{{\lambda_1}\over\scriptstyle a}$,
the linear problem
\neweq{4.17}\itab
$ -\Delta u-\lambda au=\lambda b,$ &\quad $\hbox{in}\ \ \Omega $\\
$u=0,$ & \quad $\hbox{on}\ \ \partial\Omega $\\
\ttab
 \endeq
has a unique solution in \huo, which, moreover, is positive,
by Stampacchia's Maximum Principle.

Since $\Omega$ is smooth and $-\Delta u=\lambda au+\lambda b\in\huo$, we find
$u\in H^3(\Omega )$. From now on, with a standard
 bootstrap regularization, it follows that $u\in H^\infty (\Omega )$, that is
$u\in C^\infty (\overline{\Omega})$. So, we have proved the existence of a smooth
solution to the problem \eq{P}, for every
 $0<\lambda <{{\lambda_1}\over\scriptstyle a}$.

For concluding the proof, it is sufficient to show that the
problem \eq{P}
has no solution if $\lambda^*={{\lambda_1}\over\scriptstyle a}$. Indeed,
if $u$ would be a solution, by multiplication in \eq{P} with $\varphi_1$ and
integration by parts, it follows that
 $\intom \varphi_1=0$, which contradicts $\varphi_1 >0$ in $\Omega$.
\qed

\medskip

\begin{lemma}\label{4.2.4}
The following hold:

i)$\ \ \lambda^*\geq{{\lambda_1}\over\scriptstyle a}\, .$

ii) if \eq{P} has a solution for $\lambda =\lambda^*$, then it is necessarily unstable.

iii) the problem \eq{P} has at most one solution for $\lambda =\lambda^*$.

iv) $\ u(\lambda )$ is the unique solution $u$ of the problem \eq{P} such that
$\lambda_1(\lambda f'(u))\geq 0$.\rm
\end{lemma}

\medskip
{\bf Proof.} i) By the theorem of sub and super solutions, it is sufficient to show that, for every
$0<\lambda <{{\lambda_1}\over\scriptstyle a}$, the problem has a sub and a super
solution. More precisely, we show that there exists $\Uinf ,\Usup\in
C^2(\Omega )\cap C(\overline{\Omega})$ with $\Uinf\leq\Usup$ and such that
$$ \itab
$-\Delta\Usup \geq\lambda f(\Usup ),$ &\quad $\hbox{in}\ \Omega$\\
$\Usup\geq 0,$ & \quad $\hbox{on}\ \partial\Omega $\\
\ttab $$
and $\Uinf$ verifies a similar inequality, but with reversed signs.

Let \Usup be the solution of the problem \eq{4.17} for $b=f(0)$
and $\Uinf \equiv 0$. From $f(t)\leq at+b$, for every $t>0$ and
$\Usup >0$ in $\Omega$ it follows that $f(\Usup )\leq a\Usup +b$,
which implies
 $-\Delta \Usup\geq \lambda f(\Usup )$,
in $\Omega$.

The fact that \Uinf is subsolution is obvious.

ii) Assume that the problem \eq{P} has a stable solution $u^*$ for
$\lambda =\lambda^*$. Consider the operator
$$G:\{ u\in C^{2,{1\over 2}}(\overline{\Omega});\ u=0\ \ \hbox{on}
\ \ \partial\Omega\}
\times\RR\rightarrow C^{0,{1\over 2}}(\overline{\Omega}),$$
defined by
$$G(u,\lambda )=-\Delta u-\lambda f(u)\, .$$

Applying the Implicit Function Theorem to $G$
it follows that the problem \eq{P} has a solution for
$\lambda$ in a neighbourhood of $\lambda^*$, which contradicts the maximality of
 $\lambda^*$.

iii) Let $u$ be a solution corresponding to $\lambda =\lambda^*$. Then $u$
is a supersolution for the problem \eq{P}, for every $\lambda\in (0,\lambda^*)$,
that is $u\geq u(\lambda )$, for every $\lambda\in (0,\lambda^*)$.
Since the map $\lambda\longmapsto u(\lambda )$ is increasing,
we obtain, by the monotone convergence theorem, that
 there exists $u^*\leq u$ such that
$$u(\lambda )\rightarrow u^*\quad\hbox{in}\ \ L^1(\Omega )\, .$$

Since $-\Delta u(\lambda )=\lambda f(u(\lambda ))$, for every
$\lambda\in (0,\lambda^*)$, it follows that $-\Delta
u^*=\lambda^*f(u^*)$. In order to show that $u^*$ is solution of
the problem \eq{P} for $\lambda =\lambda^*$, it is enough to show
that $u\in\huo$. Indeed, it follows then, by bootstrap, that if
$N>2$, then
$$-\Delta u^*\in L^{2^*}(\Omega )$$
and, consequently,
$$u^*\in W^{2,2^*}(\Omega )\, .$$
(Here we have denoted by $2^*$ the critical Sobolev exponent
of 2, that is, $2^*={{2N}\over{N-2}}$).

If $N=1,2$, since $u^*\in\huo$, it follows that  $-\Delta u^*\in L^4(\Omega )$
and, by Theorems 8.34 and 9.15 in \cite{GT}, we find that
$$u^*\in C^{0,{1\over 2}}(\overline{\Omega}).$$
Applying now Theorem 4.3 in Gilbarg-Trudinger
\cite{GT}, it follows that $u^*$
is a solution for the problem \eq{P}.

Let us show now that $u(\lambda )$ is bounded in \huo. Indeed, by
multiplication in \eq{P} with $u(\lambda )$ and integration by
parts, we find
$$\intom |\nabla u(\lambda )|^2=\lambda\intom
 f(u(\lambda ))u(\lambda )\leq \lambda^*\intom uf(u)\, .$$
Hence there exists $v\in\huo$ such that, up to a subsequence,
$$u(\lambda )\rightharpoonup v\quad\hbox{weakly in}\ \ \huo,\ \hbox{if}\ \
\lambda\rightarrow\lambda^*$$
and
$$u(\lambda )\rightarrow v\quad\hbox{a.e. in}\ \ \Omega,\ \hbox{if}
\ \ \lambda\rightarrow\lambda^*\, .$$

But $u(\lambda )\rightarrow u^*$ a.e. in $\Omega$. Consequently, $v=u^*$,
that is $u^*\in\huo$ and
$$u(\lambda )\rightharpoonup u^*\quad\hbox{weakly in}\ \ \huo,\ \hbox{if}\ \
\lambda\rightarrow\lambda^*\, .$$

For concluding the proof, it remains to show that
$u=u^*$. Let $w=u-u^*\geq 0$. Then
\neweq{4.18}
-\Delta w=\lambda^*(f(u)-f(u^*))\geq\lambda^*
f'(u^*)w\, .
\endeq

We also have
$$\lambda_1(\lambda^*f'(u^*))\leq 0\, .$$

By Lemma \ref{4.2.2} it follows that, either $w=0$, or $w>0$ and,
in both cases, $-\Delta w=\lambda^*f'(u^*)w$. If $w>0$, then, by the last equality
and by
 \eq{4.18} it follows that $f$ is linear in all intervals
$[u^*(x),u(x)]$, for every $x\in\Omega$. This implies easily that
 $f$ is linear in the interval $\displaystyle [0,\max_{\Omega}u]$,
which contradicts Lemma \ref{4.2.3}.

iv) Assume that the problem \eq{P}
 has a solution  $u\not= u(\lambda )$ with
$\lambda_1(\lambda f'(u))\geq 0$. Then, by Hopf's strong maximum principle,
(Theorem 3.5. in \cite{GT}), $u>u(\lambda )$. Put $w=u-u(\lambda )>0$. It follows that
\neweq{4.19}
-\Delta w=\lambda (f(u)-f(u(\lambda )))\leq \lambda
f'(u)w\, .\endeq
Multiplying \eq{4.19} by $\varphi =\varphi_1(\lambda f'(u))$ and integrating by parts,
 we find
$$\lambda\intom f'(u)\varphi w+\lambda_1(\lambda f'(u))\intom \varphi w\leq
\lambda\intom f'(u)\varphi w\, .$$
So, $\lambda_1(\lambda f'(u))=0$ and in \eq{4.19} the equality holds, which means that
$f$ is linear in the interval $\displaystyle [0,\max_{\Omega}u]$. Therefore
$$0=\lambda_1(\lambda f'(u))=\lambda_1(\lambda f'(u(\lambda ))),$$
contradiction. \qed

\medskip
The following result is a reformulation of Theorem 4.1.9.
in H\"ormander \cite{Ho}.

\medskip
\begin{lemma}\label{4.2.5}
 Let $(u_n)$ be a sequence of superharmonic nonnegative functions
defined on $\Omega$. The following alternative holds:

either

i) $\limn u_n=\infty$, uniformly on compact subsets of $\Omega$

or

ii) $(u_n)$ contains a subsequence which converges in
 $L^1_{\rm{loc}}(\Omega )$
to some $u^*$.\rm
\end{lemma}

\medskip
\begin{lemma}\label{4.2.6}  The following conditions are equivalent:

i) $\ \ \lambda^*={{\lambda_1}\over \scriptstyle a}$.

ii) the problem \eq{P} has no solution if  $\lambda =\lambda^*$.

iii) $\ \displaystyle \lim_{\lambda\rightarrow\lambda^*}u(\lambda )=\infty$,
uniformly on every compact subset of
 $\Omega$.\rm\end{lemma}

\medskip
{\bf Proof.} i)$\Longrightarrow$ii) Let us first assume that there exists a solution
$u$ provided $\lambda =\lambda^*$. As we have already observed in
Lemma \ref{4.2.4}, $u$ is
necessarily unstable. On the other hand,
$$\lambda_1(\lambda^* f'(u))\geq \lambda_1(\lambda^* a)=0$$

Hence $\lambda_1(\lambda^*f'(u))=0$, that is $f'(u)=a$, which contradicts Lemma
 \ref{4.2.3}.

ii)$\Longrightarrow$iii) Assume the contrary. We first prove that the sequence
$u(\lambda )$, for $0<\lambda <\lambda^*$, is bounded in $L^2(\Omega )$.
Indeed, if not, passing eventually to a subsequence, we may assume that
 $u(\lambda )=k(\lambda )w(\lambda )$, with
$$\intom w^2(\lambda )=1\quad\hbox{and}\quad
\lim_{\lambda\nearrow\lambda^*}k(\lambda )=\infty\, .$$
Thus, by Lemma \ref{4.2.5}, going again to a subsequence,
$$u(\lambda )\rightarrow u^*\ \ \ \hbox{in}\ \ L^1_{\rm{loc}}(\Omega ),
\quad\hbox{if}\ \lambda\rightarrow\lambda^*\, .$$
Therefore
$${{\lambda}\over{k(\lambda )}}f(u(\lambda ))\rightarrow 0\ \ \ \hbox{in}
\ \
L^1_{\rm{loc}}(\Omega ),$$
that is,
\neweq{4.20}-\Delta w(\lambda )\rightarrow 0\quad\hbox{in}
\ \ L^1_{\rm{loc}}(\Omega )\, .\endeq

We prove in what follows that
 $(w(\lambda ))$ is bounded in \huo . Indeed,
$$\intom |\nabla w(\lambda )|^2=\intom-\Delta w(\lambda )w(\lambda )=
\intom {{\lambda}\over{k(\lambda )}}f(u(\lambda ))w(\lambda )\leq$$
$$\leq \lambda^*\intom (aw^2(\lambda )+{{f(0)}\over{k(\lambda )}}
w(\lambda ))\leq\lambda^*a+c\intom w(\lambda )\leq$$
$$\leq \lambda^*a+C\sqrt{|\Omega |}\, ,$$
where $C>0$ is a constant independent on $\lambda$.

Let $w\in\huo$ be such that, passing again to a subsequence,
\neweq{4.21}
w(\lambda )\rightharpoonup w\quad\hbox{weakly in }\ \ \huo ,\ \ \hbox{if}\ \
\lambda\rightarrow\lambda^*\endeq
$$w(\lambda )\rightarrow w\quad\hbox{in}\ \ L^2(\Omega )\ \ \hbox{if}
\ \lambda\rightarrow\lambda^*\, .$$

By \eq{4.20} and \eq{4.21} it follows that
$$-\Delta w=0,\ \ w\in\huo,\ \ \intom w^2=1,$$
which yields a contradiction. Thus
$(u(\lambda ))$ is bounded in $L^2(\Omega )$. With the same arguments as above,
 $(u(\lambda ))$
is bounded in \huo. Let $u\in\huo$ be such that, up to a subsequence,
$$u(\lambda )\rightharpoonup u\quad\hbox{weakly in}\ \ \huo,\ \hbox{if}
\ \ \lambda\rightarrow\lambda^*$$
$$u(\lambda )\rightarrow u\quad\hbox{in}\ \ L^2(\Omega ),\ \hbox{if}
\ \ \lambda\rightarrow\lambda^*\, .$$

It follows that $-\Delta u=\lambda^*f(u)$, that is $u$ is a solution of the problem
\eq{P} for
$\lambda =\lambda^*$, contradiction.

iii)$\Longrightarrow$ii) As observed, if \eq{P}
 has a solution provided
$\lambda =\lambda^*$, then it is necessarily equal to
$\displaystyle\lim_{\lambda\rightarrow\lambda^*}u(\lambda )$, which is not possible
in our case.

[iii) and ii)]$\Longrightarrow$i) Let $u(\lambda )=k(\lambda )w(\lambda )$, with
 $k(\lambda )$
and $w(\lambda )$ as above. With the same proof one can show that, $(w(\lambda ))$
is bounded in \huo . Let $w\in\huo$ such that , up to a subsequence,
$$w(\lambda )\rightharpoonup w\quad\hbox{weakly in}\ \ \huo,\ \ \hbox{if}\ \
\lambda\rightarrow
\lambda^*$$
$$w(\lambda )\rightarrow w\quad\hbox{in}\ \ L^2(\Omega ),\ \ \hbox{if}\ \
 \lambda\rightarrow
\lambda^*\, .$$

Then
$$-\Delta w(\lambda )\rightarrow -\Delta w,\quad\hbox{in}\ \ {\cal D}'(\Omega )
$$
and
$${\lambda\over{{k(\lambda )}}}f(u(\lambda ))\rightarrow\lambda^*aw,\quad\hbox{in}\
\ L^2(\Omega )\, .$$
It follows that
$$-\Delta w=\lambda^*aw\quad\hbox{in}\ \ \Omega ,\ w\in\huo ,\ w\geq 0\ \ \hbox{and}
\ \intom w^2=1\, .$$ This means that
$\lambda^*={{\lambda_1}\scriptstyle\over a}$ (and $w=\varphi_1$).
\qed

\medskip
\begin{lemma}\label{4.2.7}  The following conditions are equivalent:

i)$\ \ \lambda^*>{{\lambda_1}\over\scriptstyle a}\, .$

ii) the problem \eq{P} has exactly one solution, say  $u^*$, corresponding to
 $\lambda =\lambda^*$.

iii) $\ u(\lambda )$ converges uniformly in $\overline{\Omega}$ to $u^*$, which is
the unique solution of the problem \eq{P}
 provided $\lambda =\lambda^*$.\rm\end{lemma}

\medskip
{\bf Proof.} We have already remarked that $\lambda^*\geq{{\lambda_1}\over
\scriptstyle a}$. So, this lemma becomes a reformulation of the preceding one, but
with the difference that the limit appearing in
 iii) is uniform in
$\overline{\Omega}$. Since $(u(\lambda ))$ converges a.e.
in $\Omega$ to $u^*$, it is sufficient to show that $u(\lambda )$
has a limit in $C(\overline{\Omega})$ if $\lambda\rightarrow\lambda^*$.
Even less, it is enough to show that  $u(\lambda )$ is relatively compact
 in $C(\overline{\Omega})$. This follows easily by the Arzela-Ascoli theorem,
if we show that $u(\lambda )$ is bounded
 in $C^{0,{1\over 2}}(\overline{\Omega})$. From $0<u(\lambda )<u^*$ we get
 $0<f(u(\lambda ))<f(u^*)$, which yields a uniform bound for
  $-\Delta u(\lambda )$ in
$L^{2N}(\Omega )$. The desired bound for $u(\lambda )$ is now a consequence of the
theorem
 8.34 and of the remark at the page 212 in \cite{GT}, as well as
 of the Closed Graph
 Theorem.
\qed

\medskip
{\bf Proof of Theorem \ref{4.2.1}} By Lemma \ref{4.2.6}, the assertions i),
ii) and iv) are equivalent. We shall deduce that $\lambda^*={{\lambda_1}\over
\scriptstyle a}$, by showing that the problem \eq{P} has no solution
 if $\lambda ={{\lambda_1}\over\scriptstyle a}$.
Indeed, if $u$ would be such a solution then
\neweq{4.24}
-\Delta u=\lambda f(u)\geq\lambda_1u\, .\endeq
Multiplying \eq{4.24} by $\varphi_1$ and integrating by parts,
it follows that $\lambda f(u)=\lambda_1u$, which contradicts $f(0)>0$.

iii) By Lemma \ref{4.2.4} iv), it is enough to show that  if
$0<\lambda <{{\lambda_1}\over \scriptstyle a}$, then every solution $u$
verifies $\lambda_1(\lambda f'(u))\geq 0$. On the other hand,
$$-\Delta -\lambda f'(u)\geq -\Delta -\lambda a,$$
which means that
$$\lambda_1(\lambda f'(u))\geq\lambda_1(\lambda a)=
\lambda_1 -\lambda a>0\, .$$ \qed

\medskip
\begin{teo}\label{4.2.8}  Assume that
$$\lim_{t\rightarrow\infty}(f(t)-at)=l<0\, .$$

Then

i) $\ \ \lambda^*\in ({{\lambda_1}\over\scriptstyle a},{{\lambda_1}\over
{\lambda_0}})$, where $\lambda_0=\min\{ {{f(t)}\over t} \ ;\ \ t>0\}$.

ii) for every $\lambda =\lambda^*$, the problem \eq{P} has exactly one solution,
say
$u^*$.

iii) $\ \ \displaystyle\lim_{\lambda\rightarrow\lambda^*}u(\lambda )=u^*$,
 uniformly in $\Omega$.

iv) if $\lambda\in (0,{{\lambda_1}\over\scriptstyle a}]$, then $u(\lambda )$ is
the unique solution of the problem \eq{P}.

v) if $\lambda\in ({{\lambda_1}\over\scriptstyle a},\lambda^*)$, the the
problem \eq{P} has at least an unstable solution, say $v(\lambda )$. Moreover,
for every such a solution $v(\lambda )$,

vi) $\ \ \displaystyle\lim_{\lambda\rightarrow{{\lambda_1}
\over a}}v(\lambda )=\infty$, uniformly on compact subsets of $\Omega$.

vii)$\ \ \ \displaystyle\lim_{\lambda\rightarrow\lambda^*}v(\lambda )=u^*$,
 uniformly in $\Omega$.\rm
 \end{teo}

\medskip
{\bf Proof.} i) In order to show that $\lambda^*\leq{{\lambda_1}\over
\scriptstyle{\lambda_0}}$, we shall prove that the problem \eq{P} has no solution
if
  $\lambda =
{{\lambda_1}\over\scriptstyle{\lambda_0}}$. Contrary, let $u$ be such a
solution. Multiplying \eq{P} by  $\varphi_1$ and integrating by parts, we find
\neweq{4.25}
\lambda_1\intom\varphi_1u=\lambda_1\intom\varphi_1f(u)\, .\endeq
Since $\lambda ={{\lambda_1}\over\scriptstyle{\lambda_0}}$, the relation
\eq{4.25} becomes
$$\lambda_1\intom\varphi_1u={{\lambda_1}\over
{\lambda_0}}\intom\varphi_1f(u)\geq\lambda_1\intom\varphi_1u,$$
which implies $f(u)=\lambda_0u$. As above, this equality contradicts $f(0)>0$.

The inequality $\lambda^*>{{\lambda_1}\over\scriptstyle a}$, as well as
  ii), iii)
are equivalent, by Lemmas \ref{4.2.4} and \ref{4.2.7}. We shall prove by contradiction  that
$\lambda^*>{{\lambda_1}\over\scriptstyle a}$. In this case,
$$\lambda^*={{\lambda_1}\over a}$$
and
$$\lim_{\lambda\rightarrow\lambda^*}u(\lambda )=
\infty,\quad\hbox{uniformly on compact subsets of }\ \ \Omega\, .$$

By \eq{4.25},
$$0=\intom\varphi_1\cdot [\lambda_1u(\lambda )-
\lambda f(u(\lambda ))]=$$
$$=\intom\varphi_1\cdot [(\lambda_1-a\lambda )u(\lambda )-
\lambda (f(u(\lambda ))-
au(\lambda ))]\geq$$
$$\geq -\lambda\intom \varphi_1\cdot [f(u(\lambda ))-au(\lambda ))]\, .$$
Passing to the limit as $\lambda\nearrow\lambda^*$ and taking
into account that $l<0$, we find
$$0\geq -l\lambda\intom\varphi_1>0,$$
a contradiction.

Since $\lambda^*\leq{{\lambda_1}\over\scriptstyle{\lambda_0}}$, and the problem
\eq{P} has a solution if $\lambda =\lambda^*$, it follows that
$\lambda^*<{{\lambda_1}\over\scriptstyle{\lambda_0}}$.

In order to prove iv) we follow the same technique as in the proof of
Theorem \ref{4.2.1}, iii).

v) Consider the functional
$$J:\huo\rightarrow\RR,\quad J(u)={1\over 2}\intom |\nabla u|^2-\intom F(u),$$
where
$$F(t)=\lambda \int_0^tf(s)ds\, .$$

For $\lambda\in ({{\lambda_1}\over\scriptstyle a},\lambda^*)$ we shall obtain
the unstable solution as a critical point of $J$, which is different from
 $u(\lambda )$.

We shall make use the following result, which can be found
 in Brezis-Nirenberg \cite{BN1}:

\medskip
\begin{lemma}\label{4.2.9}
  The functional $J$ has the following properties:

i) $\ J$ is of class $C^1$.

ii) for every $u,v\in\huo$,
$$\langle J'(u),v\rangle =\intom \nabla u\cdot\nabla v-\lambda\intom f(u)v$$

iii) $\ u_0=u(\lambda )$ is a local minimum point of $J$.\rm
\end{lemma}

\medskip
The idea is to apply the Mountain-Pass Lemma. For this aim we shall modify the
functional
$J$, such that $u_0$ would become a strict minimum point. Therefore,
for every $\varepsilon >0$, let us define
$$J_{\varepsilon}:\huo\rightarrow\RR,\quad J_{\varepsilon}(u)=J(u)
+{\varepsilon\over 2}\,
\intom |\nabla (u-u_0)|^2\, .$$

By the preceding Lemma,

i) $\ J$ is of class $C^1$.

ii) $\ \langle J'_{\varepsilon}(u),v\rangle =\intom \nabla u\cdot\nabla v-
\lambda\intom
f(u)v+\varepsilon\intom\nabla (u-u_0)\cdot\nabla v\, .$

iii) for every $\varepsilon >0$, $u_0$ is a strict local minimum point of
 $J_{\varepsilon}$.

\medskip
\begin{lemma}\label{4.2.10}  Let $\varepsilon_0 ={{\lambda a-\lambda_1}\over{2
\lambda _1}}$. Then there exists
 $v_0\in\huo$ such that $J_{\varepsilon}(v_0)<J_{\varepsilon}(u_0)$, for every
  $\varepsilon\in [0,\varepsilon_0]$.\rm
  \end{lemma}

\medskip
{\bf Proof of the lemma.} We first remark that $J_{\varepsilon}(u)$
is bounded by $J_0(u)$ and $J_{\varepsilon_{\scriptstyle {0}}}(u)$. Therefore,
for concluding the proof it is enough to show that
$$\lim_{t\rightarrow\infty}J_{\varepsilon_0}(t\varphi_1)=-\infty\, .$$
But
\neweq{4.27}
J_{\varepsilon}(t\varphi_1)={{\lambda_1}\over 2}t^2+{{\varepsilon_0}\over 2}
\lambda_1t^2-
\varepsilon_0\lambda_1t^2\intom\varphi_1u_0+{{\varepsilon_0}\over 2}
\intom |\nabla u_0|^2-\intom F(t\varphi_1)\, .\endeq

Set $\alpha ={{3a\lambda +\lambda_1}\over{4\lambda}}<a$. It follows that there
exists a real number $\beta$ such that $f(s)\geq\alpha s+\beta$, for every
$s\in\RR$.
Consequently, for every $t\geq 0$,
$$F(t)\geq{{\alpha\lambda}\over 2}\, t^2+\beta\lambda t\, .$$
>From here and from \eq{4.27} we find
$$\limsup_{t\rightarrow\infty}{1\over{t^2}}\,
 J_{\varepsilon_0}(t\varphi_1)\leq{{\lambda_1+\varepsilon_0
 \lambda_1-\lambda\alpha}\over 2}<0,$$
by our choice of $\alpha$. \qed

\medskip
\begin{lemma}\label{4.2.11}  The Palais-Smale condition
is satisfied uniformly with respect to
 $\varepsilon$. More precisely, if
\neweq{4.28}(J_{\varepsilon_n}(u_n))\ \ \hbox{is bounded in}\ \
\RR,\ \hbox{for}\ \ \varepsilon_n\in [0,\varepsilon_0]
\endeq
and
\neweq{4.29}J'_{\varepsilon_n}(u_n)\rightarrow 0\quad\hbox{in}\ \ \ H^{-1}(\Omega ),
\endeq
then $(u_n)$ is relatively compact  in \huo .\rm
\end{lemma}

\medskip
{\bf Proof of the lemma.} It is enough to show that $(u_n)$
has a subsequence converging in $H^1_0(\Omega )$. Indeed, in this case, up to another
subsequence, there exist $u\in\huo$ and $\varepsilon\geq 0$ such that
$$u_n\rightharpoonup u\quad\hbox{weakly in}\ \ \huo$$
$$u_n\rightarrow u\quad\hbox{in}\ \ L^2(\Omega )$$
$$u_n\rightarrow u\quad\hbox{a.e. in}\ \ \Omega$$
$$\varepsilon_n\rightarrow\varepsilon\, .$$
By \eq{4.29} it follows that
\neweq{4.30}-\Delta u_n-\lambda f(u_n)-\varepsilon_n\Delta
(u_n-u_0)\rightarrow 0\quad\hbox{in}\ \ {\cal D}'(\Omega )
\, .\endeq

By $|f(u_n)-f(u)|\leq a|u_n-u)|$ it follows that
$$f(u_n)\rightarrow f(u)\quad\hbox{in}\quad L^2(\Omega )\, .$$
Using now \eq{4.30} we obtain
$$-(1+\varepsilon_n)\Delta u_n\rightarrow\lambda f(u)
-\varepsilon\Delta u_0\quad\hbox{in}\ \ {\cal D}'(\Omega ),$$
that is,
$$-\Delta u-\lambda f(u)-\varepsilon\Delta (u-u_0)=0\, .$$
Multiplying this equality by $u$ and integrating we find
\neweq{4.32}
(1+\varepsilon)\intom |\nabla u|^2-\lambda\intom uf(u)-\varepsilon\lambda
\intom uf(u_0)=0\, .
\endeq
By multiplication in \eq{4.29} with $(u_n)$ and integration we obtain
\neweq{4.33}
(1+\varepsilon_n)\intom |\nabla u_n|^2-\lambda\intom u_nf(u_n)-
\varepsilon_n
\lambda\intom u_nf(u_0)\rightarrow 0,
\endeq
by the boundedness of $(u_n)$. The second term appearing in
\eq{4.33} tends to
$-\lambda\intom uf(u)$, while the last tends to
 $-\varepsilon\lambda\intom uf(u_0)$,
by the convergence in $L^2(\Omega )$ of the sequences $(u_n)$ and
$(f(u_n))$. So, by comparing the first terms of
 \eq{4.32} and \eq{4.33}, it follows that
$$u_n\rightarrow u\quad\hbox{in}\quad\huo\, .$$

At this stage it is sufficient to prove that $(u_n)$ contains a subsequence
which is bounded in $L^2(\Omega )$. Indeed, the boundedness
in $L^2(\Omega )$ of the sequence $(u_n)$ implies the boundedness
in \huo,
as follows by \eq{4.28}.
Arguing by contradiction, let us assume that
$$\| u_n\|_{L^2(\Omega )}\rightarrow\infty\, .$$
Put $u_n=k_nw_n$, with
$$k_n>0,\ k_n\rightarrow\infty\ \ \hbox{and}\ \ \intom w_n^2=1\, .$$
We can assume that $\varepsilon_n\rightarrow\varepsilon$. So,
\neweq{4.34}
0=\limn {{J_{\varepsilon_n(u_n)}}\over{k_n^2}}=\endeq
$$=\limn [{1\over 2}\intom |\nabla w_n|^2-{1\over{k_n^2}}\intom F(u_n)+
{{\varepsilon_n}\over 2}\intom |\nabla (w_n-{{u_0}\over{k_n}})|^2]\, .$$
On the other hand,
$$\intom |\nabla (w_n-{{u_0}\over{k_n}})|^2=\intom |\nabla w_n|^2+
{1\over{k_n^2}}\intom |\nabla u_0|^2-
{{2\lambda}\over{k_n}}\intom w_nf(u_0)\, .$$
Relation \eq{4.34} becomes
$$\limn \ [{{1+\varepsilon_n}\over 2}\intom |\nabla w_n|^2
-{1\over{k_n^2}}\intom F(u_n)]=0\, .$$
But
$$|F(u_n)|=|F(k_nw_n)|\leq{{\lambda a}\over 2}k_n^2w_n^2+\lambda b|k_nw_n|,$$
because $\ |f(t)|\leq a|t|+b$, where $b=f(0)$.
Thus the sequence
$$\left({1\over{k_n^2}}\intom F(u_n)\right)$$
is bounded, which implies also the boundedness of $(w_n)$  in \huo.
Let $w\in\huo$ such that, up to a subsequence,
$$w_n\rightharpoonup w\quad\hbox{weakly in}\ \ \huo\, ,$$
$$w_n\rightarrow w\quad\hbox{strongly in}\ \ L^2(\Omega )\, ,$$
$$w_n\rightarrow w\quad\hbox{a.e. in}\ \ \Omega\, .$$
We also remark that $\intom w^2=1$.

We prove in what follows that
\neweq{4.35}
-(1+\varepsilon )\Delta w=\lambda aw^+\, .
\endeq
Indeed, dividing in \eq{4.29} by $k_n$ we find
\neweq{4.36}
(1+\varepsilon_n)\intom\nabla w_n\cdot\nabla v-\lambda\intom{{f(u_n)}
\over{k_n}}v-
{{\varepsilon_n\lambda}\over{k_n}}\intom f(u_0)v\rightarrow 0\, ,
\endeq
for all $v\in\huo$. We remark that
$$(1+\varepsilon_n)\intom\nabla w_n\cdot\nabla v\rightarrow
(1+\varepsilon )\intom
\nabla w\cdot\nabla v\, .$$
Relation \eq{4.35} follows from \eq{4.36}
 if we show that the sequence
$({1\over{k_n}}f(u_n))$ contains a subsequence which converges
 in $L^2(\Omega )$ to
$aw^+$.

Since
$${1\over{k_n}}f(u_n)={1\over{k_n}}f(k_nw_n),$$
it is obvious that the needed limit $aw^+$ is in the set
$$\{ x\in\Omega ;\ \ w_n(x)\rightarrow w(x)\not= 0\}\, .$$

If $w(x)=0$ and $w_n(x)\rightarrow w(x)$, let $\varepsilon >0$
and $n_0$ such that $|w_n(x)|<\varepsilon$, for every $n\geq n_0$. So,
$${{f(k_nw_n)}\over{k_n}}\leq\varepsilon a+{b\over{k_n}},\quad\hbox{for every}
\ \ n\geq n_0,$$
that is the asked limit is 0. Hence
$${{f(u_n)}\over{k_n}}\rightarrow aw^+,\quad\hbox{a.e. in}\ \
\Omega\, .$$

Since $w_n\rightarrow w$ in $L^2(\Omega )$, it follows that
(see Theorem IV.9 in Brezis
\cite{Br}), up to a subsequence, $(w_n)$
is dominated in $L^2(\Omega )$. From
$${1\over{k_n}}f(u_n)\leq a|w_n|+{1\over{k_n}}b\, ,$$
it follows that $({1\over{k_n}}f(u_n))$ is also dominated in
$L^2(\Omega )$. It follows that the relation \eq{4.35} is true.

By the Maximum Principle applied  in \eq{4.35}, we find that
 $w\geq 0$ and
 \neweq{4.37}
\left\{\begin{array}{lll}
\displaystyle -\Delta w = \frac{\lambda a}{1+\varepsilon}w \\
\displaystyle w \geq 0 \\
\displaystyle\int_{\Omega}\ w^2 = 1\, .
\end{array}\right.
\endeq
So, ${{\lambda a}\over{1+\varepsilon}}=\lambda_1$ and
$w=\varphi_1$,
which contradicts $\varepsilon\in [0,\varepsilon_0]$ and the
 choice of
 $\varepsilon_0$.
\qed

\medskip
For $u_0=u(\lambda )$ and $v_0$ found in Lemma \ref{4.2.10}, let
$${\cal P}=\{ p\in C([0,1],\huo );\ p(0)=u_0,\ p(1)=v_0\}$$
and
$$c_{\varepsilon}=\inf_{p\in{\cal P}}\max_{t\in [0,1]}
J_{\varepsilon}(p(t))\, .$$

\medskip
\begin{lemma}\label{4.2.12}
 $\ \ c_\varepsilon$ is uniformly bounded.
\end{lemma}

\medskip
{\bf Proof.} The fact that $J_{\varepsilon}$ increases with
$\varepsilon$ implies that $c_0\leq c_{\varepsilon}\leq
c_{\varepsilon_{\scriptstyle 0}}$, for every
$0\leq\varepsilon\leq\varepsilon_0$. \qed

\medskip
{\bf Proof of Theorem \ref{4.2.8} v) continued} For every $\varepsilon\in
(0,\varepsilon_0]$,
let $v_\varepsilon\in\huo$ such that
\neweq{4.38}
-\Delta v_{\varepsilon} ={{\lambda}\over{1+\varepsilon}}f(v_\varepsilon )+
{{\lambda\varepsilon}\over{1+\varepsilon}}f(u_0)
\endeq
and
\neweq{4.39}
J_{\varepsilon}(v_\varepsilon )=c_\varepsilon\, .
\endeq

By Lemmas \ref{4.2.11}, \ref{4.2.12} and from \eq{4.39}
 it follows that there exists
$v\in\huo$ such that
$$v_\varepsilon \rightarrow v\quad\hbox{in}\ \ \huo ,\quad\hbox{if}\ \
\varepsilon\searrow 0\, .$$
Taking into account the relation \eq{4.38} we obtain
$$-\Delta v=\lambda f(v)\, .$$

It remains to prove that $v\not= u_0$ and $v\in C^2(\Omega )
\cap C(\overline{\Omega})$.

We first remark that $v_\varepsilon$ is a solution of the equation
\eq{4.38}, which is different from  $u_0$, so it is unstable, that is
$$\lambda_1({{\lambda}\over{1+\varepsilon}}
f'(v_\varepsilon ))\leq 0\, .$$
Indeed, \eq{4.38} is an equation of the form
$$-\Delta u=g(u)+h(x),$$
where $g$ as convex and positive, while $h$ is positive.
Therefore, by the results
established in
Brezis-Nirenberg \cite{BN1}, if this equation has a solution,
then it has also a minimal solution $u$, with
$\lambda_1(g'(u))\geq 0$.
By the proof of Lemma \ref{4.2.4} iv) it follows that if
 $v$ is another solution,
then $\lambda_1(g'(v))<0$. In our case, $u_0$ plays the role of $u$, while
$v_\varepsilon$ plays the role of $v$. It remains to prove now that the limit of a
sequence of unstable solutions is unstable, too.

\medskip
\begin{lemma}\label{4.2.13}
 Let $u_n\rightarrow u$ in \huo and
 $\mu_n\rightarrow \mu$
such that $\lambda_1(\mu_nf'(u_n))\leq 0$.

Then $\lambda_1(\mu f'(u))\leq 0$.
\end{lemma}

\medskip
{\bf Proof of the lemma.} The fact that $\lambda_1(\alpha )\leq 0$ is equivalent
to the existence of some
 $\varphi\in\huo$ such that
$$\intom |\nabla\varphi |^2\leq \intom \alpha\varphi^2\quad\hbox{and}\quad
\intom\varphi^2=1\, .$$
This assertion follows from the Courant-Hilbert Principle.

Let $\varphi_n\in\huo$ be such that
\neweq{4.40}
\intom |\nabla \varphi_n|^2\leq\intom\mu_nf'(u_n)\varphi_n^2
\endeq
and
\neweq{4.41}
\intom\varphi_n^2=1\, .
\endeq
Since $f'\leq a$, we obtain by the relation \eq{4.40}
 that $(\varphi_n)$ is bounded
 in \huo. Let $\varphi\in\huo$ such that, up to a subsequence,
$$\varphi_n\rightharpoonup \varphi \quad\hbox{weakly in}\ \ \huo$$
$$\varphi_n\rightarrow\varphi\quad\hbox{strongly
in}\ \ L^2(\Omega )\, .$$
So, again up to a subsequence,
the right hand side of \eq{4.40}
 converges to $\mu\intom f'(u)\varphi^2$. This
assertion can be proved by extracting from
 $(\varphi_n)$ a subsequence dominated in
$L^2(\Omega )$, as in the proof of Theorem IV.9 in Brezis
\cite{Br}. By the relations
$$\intom\varphi^2=1\quad\hbox{and}\quad \intom
|\nabla\varphi |^2\leq\liminf_{n\rightarrow\infty} \intom
|\nabla\varphi_n|^2$$ our conclusion follows. \qed

\medskip
The fact that $v\in C^2(\Omega )\cap C(\overline{\Omega })$ follows by a standard
bootstrap argument:
$$v\in\huo\Longrightarrow f(v)\in L^{2^*}(\Omega )\Longrightarrow
v\in W^{2,2^*}(\Omega )\Longrightarrow ...$$

We point out that the main tools are:

a) if $v\in L^p(\Omega )$, then $f(v)\in L^p(\Omega )$.

b) a standard regularity result for elliptic equations (Theorem
9.15 in Gilbarg-Trudinger \cite{GT}).

c) the Sobolev embeddings.

\medskip
vi) Assume the contrary. Thus, there exists $\mu_n\rightarrow
{{\lambda_1}\over a}$ and a sequence $(v_n)$ of unstable solutions of the problem
\eq{P} corresponding to  $\lambda =\mu_n$, as well as
$v\in L^1_{\rm{loc}}(\Omega )$ such that $v_n\rightarrow v$ in
$L^1_{\rm{loc}}(\Omega )$.

We first prove that $(v_n)$ is unbounded in \huo. Indeed, if not,
Let $w\in\huo$ such that, up to a subsequence,
$$v_n\rightharpoonup w\quad\hbox{weakly in}\ \ \huo\, ,$$
$$v_n\rightarrow w\quad\hbox{in}\ \ L^2(\Omega )\, .$$
Therefore
$$-\Delta v_n\rightarrow -\Delta w\quad\hbox{in}\ \ {\cal D}'(\Omega )$$
and
$$f(v_n)\rightarrow f(w)\quad\hbox{in}\ \ L^2(\Omega ),$$
which means that
$$-\Delta w={{\lambda_1}\over a}f(w)\, .$$
It follows that $w\in C^2(\Omega )\cap C(\overline{\Omega})$, that is
$w$ is solution of the problem \eq{P}.

We prove now that
\neweq{4.42}
\lambda_1({{\lambda_1}\over a}f'(w))\leq 0\, .
\endeq
Indeed, since $v_n$ is an unstable solution, it is sufficient to show that
 $v_n\rightarrow w$ in \huo and to apply then Lemma \ref{4.2.13}. But
$$\mu_nf(v_n)\rightarrow{{\lambda_1}\over a}f(w)\quad\hbox{in}\
 \ H^{-1}(\Omega )$$
(because the convergence holds in $L^2(\Omega )$).

The operator $-\Delta :\huo\rightarrow H^{-1}(\Omega )$ is bicontinuous, so
 \eq{4.42} holds. Hence $w\not= u({{\lambda_1}\over\scriptstyle a})$,
which contradicts iv).

The fact that $(v_n)$ is unbounded in \huo implies its unboundedness in
$L^2(\Omega )$. Indeed, we have observed that the boundedness in $L^2(\Omega )$
implies the boundedness in \huo . Let $v_n=k_nw_n$ with
$k_n>0,\ \intom w_n^2=1$ and, up to a subsequence,
$k_n\rightarrow\infty$. Hence
$$-\Delta w_n={{\mu_n}\over{k_n}}f(u_n)\rightarrow 0\quad\hbox{in}\ \
L^1_{\rm{loc}}(\Omega ),$$
so also in ${\cal D}'(\Omega )$. The sequence $(w_n)$ may be assumed to be
bounded, with an argument already done.
 Let $w$ be its weak limit  in \huo . It follows that
$$-\Delta w=0\quad\hbox{and}\quad\intom w^2=1,$$
which is the desired contradiction.

\medskip
vii) As above, it is enough to prove the $L^2(\Omega )$-boundedness of
 $v(\lambda )$ for $\lambda$ in a neighbourhood of $\lambda^*$ and
to apply the uniqueness property of $u^*$. Assuming the contrary, let
 $\mu_n\rightarrow\lambda^*$ and $v_n$ be a solution of the
 problem \eq{P} for
 $\lambda =\mu_n$, with $\| v_n\|_{L^2(\Omega )}\rightarrow\infty$. Writing, as
 above, $v_n=k_nw_n$, it follows that
$$-\Delta w_n={{\mu_n}\over{k_n}}f(u_n)\, .$$
The right hand side of this relation is bounded in $L^2(\Omega )$, so
$(w_n)$ is bounded in \huo. Let $w\in\huo$ such that,
up to a subsequence,
$$w_n\rightharpoonup w\quad\hbox{weakly in}\ \ \huo$$
$$w_n\rightarrow w\quad\hbox{strongly in}\ \ L^2(\Omega )\, .$$
With an argument already done it follows that
$$-\Delta w=\lambda^*aw,\ \ w\geq 0,\ \ \intom w^2=1\, .$$
This forces the equality $\lambda^*={{\lambda_1}\over\scriptstyle
a}$, which yields a contradiction. \qed

\chapter{Nonsmooth Mountain-Pass Theory }

\section{Basic properties of locally Lipschitz functionals}

Throughout this chapter, $X$ will denote a real Banach space.
Let $X^*$ be its dual and, for every $x\in X$  and $x^* \in X^*$,
 let
$\langle x^* ,x\rangle$ be the duality pairing between
 $X^*$ and $X$.

\medskip
\begin{defin}\label
{1.1.1} \sl A functional $f:X\rightarrow\RR$ is said to be
locally Lipschitz provided that, for every  $x\in X$,
there exists a neighbourhood
$V$ of $x$ and a positive constant $k=k(V)$ depending on $V$ such that
$$|f(y)-f(z)|\leq k||y-z||,$$
for each $y,z\in V$.\rm
\end{defin}
\vspace*{0.2cm}
The set of all locally Lipschitz mappings
defined on $X$ with real values will be denoted by
$\hbox{Lip}_{loc}(X,\RR)$.

\medskip
\begin{defin}\label
{1.1.2} \sl Let $f:X\rightarrow \RR$ be a locally Lipschitz
functional and $x,v\in X$. We call the generalized directional derivative
of $f$ in  $x$ with respect to the direction $v$ the number
$$f^0(x,v)=\limsup_{y\rightarrow x \atop \lambda \searrow 0}
{{f(y+\lambda v)-f(y)}\over{\lambda}}\, .$$
\rm
\end{defin}

\medskip
It is obvious that if
$f$ is a locally Lipschitz functional, then
 $f^0(x,v)$ is a finite number and
 \neweq{1.1}
|f^0(x,v)|\leq k||v||\, .
\endeq
Moreover, if
 $x\in X$ is fixed, then the mapping $v\longmapsto f^0(x,v)$ is positive
 homogeneous and subadditive, so convex continuous. By the
 Hahn-Banach theorem, there exists a linear map $x^*:X\rightarrow\RR$
\ai for every $v\in X$,
$$x^*(v)\leq f^0(x,v)\, .$$
The continuity of
 $x^*$ is an immediate consequence of the above inequality and of
  \eq{1.1}.

\medskip
\begin{defin}\label
{1.1.3} \sl Let $f:X\rightarrow\RR$ be a locally Lipschitz
functional and
 $x\in X$. The generalized gradient (Clarke subdifferential) of
 $f$ at the point  $x$ is the nonempty subset $\partial f(x)$ of $X^*$
which is defined by
$$\partial f(x)=\{x^*\in X^*;\ f^0(x,v)\geq
\langle x^*,v\rangle ,\ \ \hbox{for all}\ \
v\in X\}\, .$$\rm
\end{defin}

\medskip
We also point out that if
 $f$ is convex then $\partial f(x)$ coincides with the subdifferential of
 $f$ in $x$ in the sense of the convex analysis, that is
$$\partial f(x)=\{x^*\in X^*;\ f(y)-f(x)\geq \langle x^*,y-x\rangle ,\ \
\hbox{for all}\ \ y\in X\}\, .$$

\medskip
We list in what follows the main properties of the Clarke gradient of a locally
Lipschitz functional. We refer to
 \cite{Cl1}, \cite{Cl2}, \cite{Ch} for further details and proofs.

a) For every $x\in X$, $\partial f(x)$ is a convex  and
$\sigma (X^*,X)$-compact set.

b) For every $x,v\in X$ the following holds
$$f^0(x,v)=\max\{\langle x^*,v\rangle ;\ x^*\in\partial f(x)\}\, .$$

c) The multivalued mapping $x\longmapsto \partial f(x)$ is upper
semicontinuous, in the sense that for every
  $x_0\in X,\ \varepsilon >0$ and $v\in X$,
there exists $\delta >0$ such that, for any $x^*\in\partial f(x)$ satisfying
$||x-x_0||<\delta$, there is some $x_0^*\in \partial f(x_0)$
satisfying
$|\langle x^*-x^*_0,v\rangle |<\varepsilon$.

d) The functional $f^0(\cdot ,\cdot )$ is upper semicontinuous.

e) If $x$ is an extremum point of $f$, then $0\in\partial f(x)$.

f) The mapping
$$\lambda (x)=\min_{x^*\in\partial f(x)}||x^*||$$
exists and is lower semicontinuous.

g) $\quad \partial (-f)(x)=-\partial f(x)$.

h) Lebourg's Mean Value Theorem: if $x$ and $y$ are two distinct point in
 $X$ then there exists a point $z$ situated in the open segment joining
$x$ and $y$ such that
$$f(y)-f(x)\in\langle \partial f(z),y-x\rangle$$

i) If $f$ has a G\^{a}teaux derivative $f'$ which is continuous in a neighbourhood
of
 $x$, then $\partial f(x)=\{f'(x)\}$. If $X$ is finite dimensional, then
$\partial f(x)$  reduces at one point if and only if
$f$ is Fr\'echet-differentiable at $x$.

\medskip
\begin{defin}\label{1.1.4}
\sl A point $x\in X$ is said to be a critical point of
the locally Lipschitz functional
 $f:X\rightarrow\RR$ if $0\in\partial f(x)$,
that is $f^0(x,v)\geq 0$, for every $v\in X$. A number $c$ is a critical value
of
 $f$ provided that there exists a critical point
 $x\in X$ such that $f(x)=c$.\end{defin}

\medskip \rm
Remark that a minimum point is a critical point. Indeed, if
 $x$ is a local minimum point, then for every  $v\in X$,
$$0\leq \limsup_{\lambda\searrow 0}
{{f(x+\lambda v)-f(x)}\over{\lambda}}\leq f^0(x,v)\, .$$

We now introduce a compactness condition for locally Lipschitz functionals.
This condition was used for the first time, in the case of
 $C^1$-functionals, by R. Palais and S. Smale (in the global variant) and by
H. Brezis, J.M. Coron and L. Nirenberg (in its local variant).

\medskip
\begin{defin}\label{1.1.5}
 If $f:X\rightarrow\RR$ is a locally Lipschitz
functional and $c$ is a real number, we say that $f$ satisfies the
 Palais-Smale condition at the level  $c$ (in short, $({\rm PS})_c$ ) if any
 sequence $(x_n)$ in $X$ satisfying $f(x_n)\longrightarrow c$
and $\lambda (x_n)\longrightarrow 0$, contains a convergent
subsequence.
The mapping $f$ satisfies the Palais-Smale condition (in short, $({\rm PS})$)
if every sequence $(x_n)$ which satisfies $(f(x_n))$ is bounded and
$\lambda (x_n)\longrightarrow 0$, has a
convergent subsequence.
\end{defin}

\section{Ekeland's Variational Principle}

\begin{teo}\label{ekeland} (Ekeland)
Let $(M,d)$ be a complete metric space and let
 $\psi :M\ri (-\infty ,+\infty ]$, $\psi\not\equiv +\infty$,
be a lower semicontinuous function which is bounded from below.
Then the following hold: for every $\ep >0$ and for any $z_0\in M$ there exists
$z\in M$ such that

(i) $\psi (z)\leq \psi (z_0)-\ep\, d(z,z_0)$;

(ii) $\psi (x)\geq \psi (z)-\ep\, d(x,z)$, for any $x\in M$.
\end{teo}

{\bf Proof}. We may assume without loss of generality that $\ep =1$. Define
the following binary relation on $M$:

$$y\leq x\qquad\mbox{if and only if}
\qquad \psi (y)-\psi (x) +d(x,y)\leq 0\, .$$
We verify obviously that ``$\leq$" is an order relation.
For arbitrary $x\in M$, set
$$S(x)=\{ y\in M; y\leq x\}\, .$$
Let $(\ep_n)$ be a sequence of positive numbers such that $\ep_n\ri 0$ and fix
$z_0\in M$. For any $n\geq 0$, let $z_{n+1}\in S(z_n)$ be such that
$$\psi (z_{n+1})\leq\inf_{S(z_n)}\psi +\ep_{n+1}\, .$$
The existence of $z_{n+1}$ follows by the definition of the set $S(x)$.
We will prove that the sequence $(z_n)$ converges to some element $z$
which satisfies (i) and (ii).

Let us first remark that $S(y)\subset S(x)$, provided that $y\leq x$.
Hence, $S(z_{n+1})\subset S(z_n)$. It follows that, for any $n\geq 0$,
$$\psi (z_{n+1})-\psi (z_n)+d(z_n,z_{n+1})\leq 0\, ,$$
which implies $\psi (z_{n+1})\leq \psi (z_n)$. Since $\psi$ is bounded
from below, it follows that the sequence $\{\psi (z_n)\}$ converges.

We prove in what follows that $(z_n)$ is a Cauchy sequence. Indeed, for
any $n$ and $p$ we have
\neweq{np}
\psi (z_{n+p})-\psi (z_n)+d(z_{n+p},z_n)\leq 0\, .
\endeq
Therefore
$$d(z_{n+p},z_n)\leq \psi (z_n)-\psi (z_{n+p})\ri 0\, ,\quad\mbox{as}\ \,
n\ri\infty\, ,$$
which shows that $(z_n)$ is a Cauchy sequence, so it converges to some $z\in M$.
Now, taking $n=0$ in \eq{np}, we find
$$\psi (z_{p})-\psi (z_0)+d(z_{p},z_0)\leq 0\, .$$
So, as $p\ri\infty$, we find (i).

In order to prove (ii), let us choose an arbitrary $x\in M$. There are two cases:

- $x\in S(z_n)$, for any $n\geq 0$. It follows that $\psi (z_{n+1})
\leq\psi (x)+\ep_{n+1}$ which implies that $\psi (z)\leq\psi (x)$.

- there exists an integer $N\geq 1$ such that $x\not\in S(z_n)$, for
any $n\geq N$ or, equivalently,
$$\psi (x)-\psi (z_n)+d(x,z_n)> 0\, ,\qquad\mbox{for every}\ \, n\geq N\, .$$
Passing at the limit in this inequality
 as $n\ri\infty$ we find (ii).
 \qed

 \medskip
 \begin{cor}\label{ek1}
Assume the same hypotheses on $M$ and $\psi$. Then, for any $\ep >0$,
there exists $z\in M$ such that
$$\psi (z)<\inf_M\psi +\ep$$
and
$$\psi (x)\geq\psi (z)-\ep\, d(x,z)\, ,\qquad\mbox{for any}\ \,
x\in M\, .$$
\end{cor}

\medskip
The conclusion follows directly from Theorem \ref{ekeland}.

\medskip
The following will be of particular interest in our next arguments.

\medskip
\begin{cor}\label{ek2}
Let $E$ be a Banach space and let $\psi :E\ri\RR$ be a $C^1$ function
which is bounded from below. Then, for any $\ep >0$ there exists
$z\in E$ such that
$$\psi (z)\leq\inf_E\psi +\ep\qquad\mbox{and}\qquad\|\psi '(z)\|_{E^\star}
\leq\ep\, . $$
\end{cor}

\medskip
{\bf Proof}. The first part of the conclusion follows directly from Theorem
\ref{ekeland}. For the second part we have
$$\|\psi '(z)\|_{E^\star}=\sup_{\| u\| =1}
\langle \psi '(z),u\rangle\, .$$
But
$$\langle \psi '(z),u\rangle =\lim_{\delta\ri 0}
{{\psi (z+\delta u)-\psi (z)}
\over{\delta \| u\|}}\, .$$
So, by Theorem \ref{ekeland},
$$\langle \psi '(z),u\rangle \geq -\ep\, .$$
Replacing now $u$ by $-u$ we find
$$\langle \psi '(z),u\rangle \leq \ep\, ,$$
which concludes our proof. \qed

\medskip
We give in what follows a variant of Ekeland's Theorem, whose proof
use in an essential manner the fact that the dimension of the
space is finite.

\medskip
\begin{teo}\label{hiri} Let $\psi :\RR^N\ri (-\infty ,+\infty ]$
be a lower semicontinuous function, $\psi\not\equiv +\infty$, bounded from below.
Let $x_\ep\in\RR^N$ be such that
\neweq{inff}
\inf \psi\leq\psi (x_\ep )\leq\inf\psi +\ep\, .
\endeq
Then, for every $\lambda >0$, there exists $z_\ep \in\RR^N$
such that

(i) $\psi (z_\ep )\leq\psi (x_\ep )$;

(ii) $\| z_\ep -x_\ep\|\leq\lambda$;

(iii) $\psi (z_\ep )\leq\psi (x)+{\ep\over\lambda}\,
\| z_\ep -x\|$, for every $x\in\RR^N$.
\end{teo}

\medskip
{\bf Proof}. Given $x_\ep$ satisfying \eq{inff}, let us consider
the function
$\varphi :\RR^N\ri (-\infty ,+\infty ]$ defined by
$$\varphi (x)=\psi (x)+{\ep\over\lambda}\, \| x-x_\ep\|\, .$$
By our hypotheses on $\psi$ it follows that $\varphi$ is
lower semicontinuous and bounded from below. Moreover,
$\varphi (x)\ri +\infty$ as $\| x\|\ri +\infty$. Therefore there
exists $z_\ep\in\RR^N$ which minimizes $\varphi$ on $\RR^N$, that is, for every
$x\in\RR^N$,
\neweq{egi}
\psi (z_\ep )+{\ep\over\lambda}\, \| z_\ep -x_\ep\|
\leq \psi (x )+{\ep\over\lambda}\, \| x -x_\ep\|\, .
\endeq
By letting $x=x_\ep$ we find
$$\psi (z_\ep )+{\ep\over\lambda}\, \| z_\ep -x_\ep\|
\leq\psi (x_\ep )\, ,$$
and (i) follows. Now, since
$\psi (x_\ep )\leq\inf\psi +\ep$, we clearly deduce from the above that
$\| z_\ep -x_\ep\|\leq\lambda$.

We infer from \eq{egi} that, for every $x\in\RR^N$,
$$\psi (z_\ep )\leq \psi (x)+{\ep\over\lambda}\, \left(
\| x-x_\ep\| -\| z_\ep -x_\ep\|\right)\leq \psi
(x)+{\ep\over\lambda}\, \| x-z_\ep\|\,$$ which is exactly the
desired inequality (iii). \qed

\medskip
The above result shows that, the closer to $x_\ep$ we desire
$z_\ep$ to be, the larger the perturbation of $\psi$ that must
be accepted. In practise, a good compromise is to take $\lambda =\sqrt\ep$.

\medskip
It is striking to remark that the Ekeland Variational Principle
characterizes the completeness of a metric space in a certain sense.
More precisely we have

\begin{teo}\label{reci}
Let $(M,d)$ be a metric space. Then $M$ is complete if and only
if the following holds: for every
application
 $\psi :M\ri (-\infty ,+\infty ]$, $\psi\not\equiv +\infty$,
which is bounded from below and for every $\ep >0$, there exists
$z_\ep\in M$ such that

\noindent (i) $\psi (z_\ep )\leq\inf_M\psi +\ep\, $

\noindent and

\noindent (ii) $\psi (z)>\psi (z_\ep )-\ep\, d(z,z_\ep )$, for any
$z\in M\setminus\{ z_\ep\}$.
\end{teo}

\medskip
{\bf Proof}. The ``only if" part  follows directly from
Corollary \ref{ek1}.

For the converse, let us assume that $M$ is an arbitrary metric
space satisfying the hypotheses. Let $(v_n)\subset M$ be an
arbitrary Cauchy sequence and consider the function $f:M\ri\RR$
defined by
$$f(u)=\lim_{n\ri\infty}d(u,v_n)\, .$$
the function $f$ is continuous and, since $(v_n)$ is a Cauchy
sequence, then $\inf f=0$. In order to justify the completeness
of $M$ it is enough to find $v\in M$ such that $f(v)=0$.
For this aim, choose an arbitrary $\ep\in (0,1)$. Now, from our
hypotheses (i) and (ii), there exists $v\in M$ with $f(v)\leq\ep$ and
$$f(w)+\ep\, d(w,v)>f(v)\, ,\qquad\mbox{for any}\ \, w\in M\setminus
\{ v\}\, .$$ From the definition of $f$ and the fact that $(v_n)$
is a Cauchy sequence we can take $w=v_k$ for $k$ large enough such
that $f(w)$ is arbitrarily small and $d(w,v)\leq\ep +\eta$, for
any $\eta >0$, because $f(v)\leq\ep$. Using (ii) we obtain that,
in fact, $f(v)\leq\ep^k$. Repeating the argument we may conclude
that $f(v)\leq\ep^n$, for all $n\geq 1$ and so $f(v)=0$, as
required. \qed

\section{Mountain Pass and Saddle Point type theorems}

Let $f:X\rightarrow\RR$ be a locally Lipschitz functional.
 Consider
$K$ a compact metric space and $K^*$ a closed nonempty subset of
 $K$. If $p^*:K^*\rightarrow X$ is a continuous mapping, set
$${\cal P}=\{p\in C(K,X);\ p=p^*\ \ \hbox{on}\ \ K^*\}\, .$$
By a celebrated theorem of Dugundji \cite{D}, the set ${\cal P}$ is nonempty.

Define
\neweq{1.2}
c=\inf_{p\in{\cal P}}\max_{t\in K}f(p(t))\, .
\endeq

Obviously, $\displaystyle c\geq\max_{t\in K^*}f(p^*(t))$.

\medskip
\begin{teo}\label
{1.2.1} \sl Assume that
\neweq{1.3}
c>\max_{t\in K^*}f(p^*(t))\, .
\endeq

Then there exists a sequence $(x_n)$ in $X$ such that

i)$\ \quad \displaystyle\lim_{n\rightarrow\infty}f(x_n)=c$;

ii)$\quad\displaystyle\lim_{n\rightarrow\infty}\lambda
(x_n)=0$.
\end{teo}

\medskip
For the proof of this theorem we need the following auxiliary result:

\medskip
\begin{lemma}\label
{1.2.2} \sl Let $M$ be a compact metric space and let
$\varphi :M\rightarrow 2^{\scriptstyle X^*}$ be an upper semicontinuous mapping
such that, for every
 $t\in M$, the set $\varphi (t)$ is convex and
$\sigma (X^*,X)$-compact. For $t\in M$, denote
$$\gamma (t)=\inf\{\ ||x^*||;\ \ x^*\in \varphi (t)\}$$
and
$$\gamma =\inf_{t\in M}\gamma (t)\, .$$

Then, for every fixed $\varepsilon >0$, there exists a continuous mapping
$v:M\rightarrow X$ such that for every $t\in M$ and
$x^*\in\varphi (t)$,
$$||v(t)||\leq 1\quad\hbox{and}\quad \langle x^*,v(t)\rangle
\geq\gamma -\varepsilon\, .$$\rm
\end{lemma}

\medskip
{\bf Proof of Lemma.} Assume, without loss of generality, that
 $\gamma >0$ and $0<\varepsilon <\gamma$. Denoting by $B_r$ the open ball in
  $X^*$ centered at the origin and with radius $r$, then, for every
$t\in M$ we have
$$B_{\gamma -{\varepsilon\over 2}}\cap\varphi (t)=\emptyset$$

Since $\varphi (t)$ and $B_{\gamma -{\epsilon\over 2}}$ are convex, disjoint and
 $\sigma (X^*,X)$-compact sets, it follows from the separation theorem
 in locally convex spaces  (Theorem 3.4 in \cite{Ru}), applied to the space
$(X^*,\sigma (X^*,X))$, and from the fact that the dual of this space
is  $X$, that: for every $t\in M$, there exists $v_t\in X$ such that
$$||v_t||=1      \quad\hbox{and}\quad\langle \xi,v_t\rangle
\leq\langle x^*,v_t\rangle\, ,$$
for any $\xi\in B_{\gamma -{\scriptstyle\epsilon\over 2}}$ and for every
$x^*\in\varphi (t)$.

Hence, for each $x^*\in\varphi (t)$,
$$\langle x^*,v_t\rangle\geq\sup_{\xi\in B_{\gamma -{\varepsilon\over 2}}}
\langle\xi ,v_t\rangle =\gamma -{\varepsilon\over 2}\, .$$

Since $\varphi$ is upper semicontinuous, there exists an open  neighbourhood
$V(t)$ of $t$ such that
for every $t'\in V(t)$ and $x^*\in \varphi (t')$,
$$\langle x^*,v_t\rangle >\gamma -\varepsilon\, .$$

Therefore, since $M$ is compact and $M=\bigcup_{t\in M}V(t)$,
there exists an open covering $\{V_1,...,V_n\}$ of $M$. Let $v_1,...,v_n$ be on the
unit sphere of $X$ such that
$$\langle x^*,v_i\rangle >\gamma -\varepsilon ,$$
for every $1\leq i\leq n$, $t\in V_i$ and $x^*\in \varphi (t)$.

If $\rho_i(t)=\hbox{dist}(t,\partial V_i)$, define
$$\zeta_i(t)={{\rho_i(t)}\over{\sum_{j=1}^n\rho_j(t)}} \quad\hbox{and}\quad
v(t)=\sum_{i=1}^n\zeta_i(t)v_i$$

It follows easily that  $v$ satisfies our conditions. \qed

\medskip
{\bf Proof of Theorem \ref{1.2.1}} We apply Ekeland's Principle to
$$\psi (p)=\max_{t\in K}f(p(t)),$$
defined on ${\cal P}$, which is a complete metric space if it is endowed with the
usual
metric. The mapping $\psi$ is continuous and bounded from below because, for every
 $p\in{\cal P}$,
$$\psi (p)\geq\max_{t\in K^*}f(p^*(t))$$
Since
$$c=\inf_{p\in{\cal P}}\psi (p),$$
it follows that for every $\varepsilon >0$, there is some
 $p\in {\cal P}$ such that
\neweq{1.4}
\psi (q)-\psi (p)+\varepsilon d(p,q)\geq 0,\quad\hbox{for all}
\ \ q\in{\cal P}\, ;
\endeq
$$c\leq \psi (p)\leq c+\varepsilon\, .
$$

Set
$$B(p)=\{t\in K;\ f(p(t))=\psi (p)\}\, .$$

For concluding the proof it is sufficient to show that there exists
 $t'\in B(p)$ such that
$$
\lambda (p(t'))\leq 2\varepsilon\, .
$$

Indeed the conclusion of the theorem follows then easily by choosing
$\varepsilon ={1\over\scriptstyle  n}$ and $x_n=p(t')$.

Applying Lemma \ref{1.2.2} for $M=B(p)$ and $\varphi (t)=\partial f(p(t))$,
we obtain a continuous map $v:B(p)\rightarrow X$ such that, for every
 $t\in B(p)$ and $x^*\in\partial f(p(t))$, we have
$$||v(t)||\leq 1\quad\hbox{and}\quad
\langle x^*,v(t)\rangle \geq\gamma -\varepsilon\, ,$$
where
$$\gamma =\inf_{t\in B(p)}\lambda (p(t))\, .$$
It follows that for every $t\in B(p)$,
$$f^0(p(t),-v(t))=\max\{\langle x^*,-v(t)\rangle;\ x^*\in\partial f(p(t))\}=$$
$$=-\min\{\langle x^*,v(t);\ x^*\in\partial f(p(t))\}\leq
 -\gamma +\varepsilon\, .$$

By \eq{1.3} we have $\ B(p)\cap K^*=\emptyset$. So, there exists a continuous
extension $w:K\rightarrow X$ of $v$ such that $\ w=0$ pe $K^*$ and, for every
 $t\in K$,
$$||w(t)||\leq 1\, .$$

Choose in the place of $q$ in \eq{1.4}
 small perturbations of the path $p$:
$$q_h(t)=p(t)-hw(t),$$
where $h>0$ is small enough.

It follows from \eq{1.4} that, for every $h>0$,
\neweq{1.8}
-\varepsilon\leq -\varepsilon
 ||w||_{\infty}\leq{{\psi (q_h)-\psi (p)}\over h}\, .
 \endeq

In what follows, $\varepsilon >0$ will be fixed, while $\ h\rightarrow 0$.
Let $t_h\in K$ be such that $f(q_h(t_h))=\psi (q_h)$. We may also assume that the
sequence
 $(t_{h_{\scriptstyle n}})$ converges to some
 $t_0$, which, obviously, is in $B(p)$. Observe that
$${{\psi (q_h)-\psi (p)}\over{h}}={{\psi (p-hw)-\psi (p)}\over{h}}
\leq{{\funu -f(p(t_h))}\over{h}}\, .$$

It follows from this relation and from \eq{1.8} that
$$-\varepsilon\leq{{\funu -f(p(t_h))}\over{h}}\leq$$
$$\leq{{f(p(t_h)-hw(t_0))-f(p(t_h))}\over{h}}+
{{\funu -f(p(t_h)-hw(t_0))}\over{h}}\, .$$

Using the fact that $f$ is a locally Lipschitz map and
$t_{h_{\scriptstyle n}}\rightarrow t_0$, we find that
$$\lim_{n\rightarrow \infty}{{\fdoi -f(p(t_{h_{\scriptstyle n}})
-h_nw(t_0))}\over{h_n}}=0\, .$$
Therefore

$$-\varepsilon\leq\limsup_{n\rightarrow\infty}{{f(p(t_0)+z_n-h_nw(t_0))-
f(p(t_0)+z_n)}\over{h_n}},$$
where $\ z_n=p(t_{h_{\scriptstyle n}})-p(t_0)$.
Consequently,
$$-\varepsilon\leq f^0(p(t_0),
-w(t_0))=f^0(p(t_0),-v(t_0))\leq -\gamma +\varepsilon\, .$$
It follows that
$$\gamma =\inf\{||x^*||;\
x^*\in\partial f(p(t)), t\in B(p)\}\leq 2\varepsilon\, .$$

Taking into account the lower semicontinuity of $\lambda$, we find the existence of
some $t'\in B(p)$ such that
$$\lambda (p(t'))=\inf\{||x^*||;\ x^*\in\partial
f(p(t'))\}\leq 2\varepsilon\, .\mfin$$

\medskip
\begin{cor}\label{1.2.3}
 If $f$ satisfies the condition $({\rm PS})_c$
and the hypotheses of Theorem \ref{1.2.1},
then $c$ is a critical value of $f$
corresponding to a critical point which is not
in  $p^*(K^*)$.
\end{cor}

\medskip
The proof of this result follows easily by applying Theorem
 \ref{1.2.1} and
the fact that $\lambda$ is lower semicontinuous. \qed

\medskip
The following result generalizes the classical
Ambrosetti-Rabinowitz Theorem.

\medskip
\begin{cor}\label{1.2.4}
 Let $f:X\rightarrow \RR$ be a locally Lipschitz
functional
such that $f(0)=0$ and there exists $v\in X\setminus\{0\}$ so that
$f(v)\leq 0$. Set
$${\cal P}=\{p\in C([0,1],X);\ p(0)=0\ \ \hbox{and}\ \ p(1)=v\}$$
$$c=\inf_{p\in{\cal P}}\max_{t\in [0,1]}f(p(t))\, .$$

If $c>0$ and $f$ satisfies the condition $({\rm PS})_c$, then $c$ is a critical
value of $f$.
\end{cor}

\medskip
For the proof of this result it is sufficient to apply Corollary
\ref{1.2.3} for $K=[0,1]$, $K^*=\{0,1\},\ p^*(0)=0$ and
$p^*(1)=v$. \qed

\medskip
\begin{cor}\label{1.2.5}
 Let $f:X\rightarrow \RR$ be a locally Lipschitz mapping.
Assume that there exists a subset $S$ of $X$ such that,
for every $p\in {\cal P}$,
$$p(K)\cap S\not =\emptyset\, .$$
If
$$\inf_{x\in S}f(x)>\max_{t\in K^*}f(p^*(t)),$$
then the conclusion of Theorem \ref{1.2.1} holds.
\end{cor}

\medskip
{\bf Proof.} It suffices to observe that
$$\inf_{p\in{\cal P}}\max_{t\in K}f(p(t))\geq
\inf_{x\in S}f(x)> \max_{t\in K^*}f(p^*(t))\, .
$$
\mfin

\medskip
Using now Theorem \ref{1.2.1} we may prove the following result, which is originally due
to Brezis-Coron-Nirenberg (see Theorem 2 in \cite{BCN}):

\medskip
\begin{cor}\label{1.2.6}
 Let $f:X\rightarrow\RR$ be a G\^ateaux differentiable
functional such that $f':(X,||\cdot ||)\rightarrow
(X^*,\sigma (X^*,X))$ is continuous. If $f$ satisfies \eq{1.3}, then there exists a
sequence $(x_n)$ in $X$ such that

i)$\ \ \displaystyle\lim_{n\rightarrow\infty}f(x_n)=c\, ;$

ii)$\ \displaystyle\lim_{n\rightarrow\infty}\| f'(x_n)\| =0\, .$

Moreover, if $f$ satisfies the condition $({\rm PS})_c$, then there exists
$x\in X$ such that $f(x)=c$ and $f'(x)=0$.
\end{cor}

\medskip
{\bf Proof.} Observe first that $f'$ is locally bounded. Indeed, if
 $(x_n)$ is a  sequence converging to $x_0$, then
$$\sup_{n}|\langle f'(x_n),v\rangle |<\infty ,$$
for every $v\in X$. Thus, by the Banach-Steinhaus Theorem,
$$\limsup_{n\rightarrow\infty}||f'(x_n)||<\infty$$

For $\lambda >0$ small enough and $h\in X$ sufficiently small we have
\neweq{1.9}
|f(x_0+h+\lambda v)-f(x_0+h)|=|\lambda\langle f'(x_0+h+\lambda\theta v),
v\rangle |\leq C||\lambda v||\, ,
\endeq
where $\theta\in (0,1)$. Therefore $f\in \hbox{Lip}_{loc}(X,\RR )$ and
$f^0(x_0,v)=\langle f'(x_0),v\rangle$, by the continuity assumption on $f'$.
In  \ref{1.9} the existence of $C$ follows
from the local boundedness property of
 $f'$.

Since $f^0$ is linear in $v$, we get
$$\partial f(x)=\{f'(x)\}$$
and, for concluding the proof, it remains to apply just Theorem
\ref{1.2.1} and Corollary \ref{1.2.3}. \qed

\medskip
A very useful result in applications is the following variant of the
 ``Saddle Point" Theorem of P.H. Rabinowitz.

\medskip
\begin{cor}\label{1.2.7}
 Assume that $X$ admits a decomposition of the form
$X=X_1\oplus X_2$, where $X_2$ is a finite dimensional subspace of
$X$. For some fixed $z\in X_2$, suppose that there exists $R>||z||$ such that
$$\inf_{x\in X_1}f(x+z)>\max_{x\in K^*}f(x)\, ,$$
where
$$K^*=\{x\in X_2;\ ||x|| = R\}\, .$$
Set
$$K=\{x\in X_2;\ ||x||\leq R\}$$
$${\cal P}=\{p\in C(K,X);\ p(x)=x\ \ \hbox{if}\ \ ||x||=R\}\, .$$

If $c$ is chosen as in \eq{1.2} and $f$ satisfies the condition
 $({\rm PS})_c$,
then $c$ is a critical value of $f$.
\end{cor}

\medskip
{\bf Proof.} Applying Corollary \ref{1.2.5} for $S=z+X_1$, we observe that it is
sufficient to prove that, for every $p\in {\cal P}$,
$$p(K)\cap (z+X_1)\not =\emptyset\, .$$

If $P:X\rightarrow X_2$ is the canonical projection, then the above condition is
equivalent to the fact that, for every $p\in {\cal P}$,
there exists $x\in K$ such that
$$P(p(x)-z)=P(p(x))-z=0\, .$$

For proving this claim, we shall make use of
an argument based on the topological
degree theory. Indeed, for every fixed $p\in{\cal P}$ we have
$$P\circ p=\hbox{Id}\ \ \hbox{on} \ K^{*}=\partial K\, .$$
Hence
$$d(P\circ p-z,\hbox{Int}\ K,0)=d(P\circ p,\hbox{Int}\ K,z)=$$
$$=d(\hbox{Id},\hbox{Int}\ K,z)=1\, .$$
Now, by the existence property of the Brouwer degree, we may find
$x\in\hbox{Int}\ K$ such that
$$(P\circ p)(x)-z=0,$$
which concludes our proof. \qed

\medskip
It is natural to ask us what happens if the condition
\eq{1.3} fails to be valid, more precisely, if
$$c=\max_{t\in K^*}f(p^*(t))\, .$$

The following example shows that in this case the conclusion of Theorem
\ref{1.2.1} does not hold.

\medskip
\begin{ex}\label
{1.2.8} Let $X=\RR^2,\ K=[0,1]\times\{0\}$,
$K^*=\{(0,0),(1,0)\}$ and let $p^*$ be the identic map of $K^*$.
As locally Lipschitz functional we choose
$$f:X\rightarrow\RR,\quad f(x,y)=x+|y|\, .$$
\end{ex}

In this case,
$$c=\max_{t\in K^*}f(p^*(t))=1\, .$$
An elementary computation shows that
$$\partial f(x,y)=\left(1\atop
1\right) ,\  \hbox{if}\ \ y>0\, ;$$
$$\partial f(x,y)=\left(1\atop
-1\right) ,\ \hbox{if}\ \ y<0\, ;$$
$$\partial f(x,0)=\left\{\left(1\atop
a\right) ;\ \hbox{if}\ \ -1\leq a\leq 1\right\}\, .$$ It follows
that $f$ satisfies the condition Palais-Smale. However $f$ has no
critical point.

The following result gives a sufficient condition for that
Theorem \ref{1.2.1}
holds provided that the condition \eq{1.3} fails.

\medskip
\begin{teo}\label
{1.2.9} \sl Assume that for every $p\in{\cal P}$ there exists
 $t\in K\setminus K^*$ such that $f(p(t))\geq c$.

Then there exists a sequence $(x_n)$ in $X$ such that

i)$\ \ \displaystyle\lim_{n\rightarrow\infty}f(x_n)=c\, ;$

ii)$\ \displaystyle\lim_{n\rightarrow\infty}\lambda (x_n)=0\, .$

Moreover, if $f$ satisfies the condition $({\rm PS})_c$, then
 $c$ is a critical value of $f$.
Furthermore, if $(p_n)$ is an arbitrary sequence in ${\cal P}$ satisfying
$$\lim_{n\rightarrow\infty}\max_{t\in K}f(p_n(t))=c,$$
then there exists a sequence $(t_n)$ in $K$ such that
$$\lim_{n\rightarrow\infty}f(p_n(t_n))=c\quad\hbox{and}
\quad\lim_{n\rightarrow\infty}\lambda (p_n(t_n))=0\, .$$\rm
\end{teo}

\medskip
{\bf Proof.} For every $\varepsilon >0$ we  apply
 Ekeland's Principle to the perturbed functional
 $$\psi_{\varepsilon}:{\cal P}\rightarrow\RR,\quad
 \psi_{\varepsilon}(p)=
\max_{t\in K}(f(p(t))+\varepsilon d(t))\, ,$$
where
$$d(t)=\min\{\hbox{dist} (t,K^*),1\}\, .$$
If
$$c_{\varepsilon}=\inf_{p\in{\cal P}}\psi_{\varepsilon}(p),$$
then
$$c\leq c_\varepsilon\leq c+\varepsilon\, .$$
Thus, by Ekeland's Principle, there exists a path $p\in{\cal P}$ such that,
for every $q\in{\cal P}$,
\neweq{1.10}
\psi_\varepsilon (q)-\psi_{\varepsilon}(p)+\varepsilon
d(p,q)\geq 0
\endeq
$$c\leq c_{\varepsilon}\leq\psi_{\varepsilon}
(p)\leq c_{\varepsilon}+\varepsilon \leq
c+2\varepsilon\, .$$
Denoting
$$ B_{\varepsilon}(p)=\{t\in K;\ f(p(t))+\varepsilon d(t)
=\psi_{\varepsilon}(p)\},$$
it remains to show that there is some $t'\in B_{\varepsilon}(p)$
such that $\lambda (p(t'))\leq 2\varepsilon$. Indeed, the conclusion of the first
part of the theorem follows easily, by choosing
 $\varepsilon
={1\over n}$ and $x_n=p(t')$.

Now, by Lemma \ref{1.2.2} applied for $M=B_{\varepsilon}(p)$ and $\varphi (t)=
\partial f(p(t))$, we find the existence of a continuous mapping
 $v:B_{\varepsilon}(p)\rightarrow X$
such that, for every $t\in B_{\varepsilon}(p)$ and
$x^*\in\partial f(p(t))$,
$$||v(t)||\leq 1\quad\hbox{and}\quad
\langle x^*,v(t)\rangle\geq\gamma_{\varepsilon}-\varepsilon,$$
where
$$\gamma_{\varepsilon}=\inf_{t\in B_{\varepsilon}(p)}\lambda (p(t))
\, .$$
On the other hand, it follows by our hypothesis that
$$\psi_{\varepsilon}(p)>\max_{t\in K^*}f(p(t))\, .$$
Hence
$$B_{\varepsilon}(p)\cap K^*=\emptyset\, .$$
So, there exists a continuous extension
 $w$ of $v$, defined on $K$
 and such that
$$w=0\ \ \hbox{ on}\ \ K^*\quad\ \hbox{ and}\ \
||w(t)||\leq 1, \quad\hbox{for any}\ \ t\in K\, .$$
Choose as $q$ in relation \eq{1.10} small variations of the path $p$:
$$q_h(t)=p(t)-hw(t),$$
for $h>0$ sufficiently small.

In what follows $\varepsilon >0$ will be fixed, while
$h\rightarrow 0$.

Let $t_h\in B_{\varepsilon}(p)$ be such that
$$f(q(t_h))+\varepsilon d(t_h)=\psi_{\varepsilon}(q_h)\, .$$

There exists a sequence $(h_n)$ converging to 0 and such that the corresponding
sequence $(t_{h_{\scriptstyle n}})$ converges to some $t_0$, which,
obviously, lies in $B_{\varepsilon }(p)$. It follows that
$$-\varepsilon\leq -\varepsilon ||w||_{\infty}\leq{{\psi_{\varepsilon}(q_h)-\psi_{\varepsilon}(p)}\over h}=
{{f(q_h(t_h))+\varepsilon d(t_h)-\psi_{\varepsilon}(p)}\over h}\leq$$
$$\leq {{f(q_h(t_h))-f(p(t_h))}\over h}={{f(p(t_h)-hw(t_h))
-f(p(t_h))}\over h}\, .$$

With the same arguments as in the proof of
 Theorem \ref{1.2.1} we obtain the existence of some
$t'\in B_{\varepsilon}(p)$ such that
$$\lambda (p(t'))\leq 2\varepsilon\, .$$

Furthermore, if $f$ satisfies (PS)$_c$ then $c$ is a critical value of
 $f$, since $\lambda$ is lower semicontinuous.

For the second part of the proof, applying again Ekeland's
Variational Principle, we deduce the existence of a sequence of
paths $(q_n)$ in ${\cal P}$ such that, for every $q\in {\cal P}$,
$$\psi_{\varepsilon_n^2}(q)-\psi_{\varepsilon_n^2}
(q_n)+\varepsilon_nd(q,q_n)\geq 0\, ;$$
$$\psi_{\varepsilon_n^2}(q_n)\leq\psi_{\varepsilon_n^2}(p_n)
-\varepsilon_nd(p_n,q_n)\, ,$$
where $(\varepsilon_n)$ is a sequence of positive numbers converging to
 0
and $(p_n)$ are paths in ${\cal P}$ such that
$$\psi_{\varepsilon_n^2}(p_n)\leq c+2\varepsilon_n^2\, .$$

Applying the same argument for $q_n$, instead of $p$,
we find $t_n\in K$ such that
$$c-\varepsilon_n^2\leq f(q_n(t_n))\leq c+2\varepsilon_n^2\, ;$$
$$\lambda (q_n(t_n))\leq 2\varepsilon_n\, .$$

We shall prove that this is the desired sequence $(t_n)$.
Indeed, by the Palais-Smale condition
(PS)$_c$, there exists a subsequence of $(q_n(t_n))$
which converges to a critical point. The corresponding subsequence of
$(p_n(t_n))$ converges to the same limit.
A standard argument, based on the continuity of $f$ and the lower
semicontinuity of
 $\lambda$ shows that for all the sequence we have
$$\limn f(p_n(t_n))=c$$
and
$$\limn \lambda (p_n(t_n))=0\, .$$
$\ $\qed

\medskip
\begin{cor}\label{1.2.10}  Let \fxr be a locally Lipschitz functional which
satisfies the Palais-Smale condition.

If $f$ has two different minimum points, then $f$ possesses
a third critical point.
\end{cor}

\medskip
{\bf Proof.} Let $x_0$ and $x_1$ be two different minimum points
of  $f$.

\underbar{Case 1.} $\ \ f(x_0)=f(x_1)=a$.
Choose $0<R<{1\over 2}\,\| x_1-x_0\|$
such that $f(x)\geq a$, for all $x\in B(x_0,R)\cup B(x_1,R)$.

Set $A=\overline{B}(x_0,{R\over 2})\cup\overline{B}(x_1,{R\over 2})$.

\underbar{Case 2.} $\ \ f(x_0)>f(x_1)$. Choose $0<R<\| x_1-x_0\|$
such that $f(x)\geq f(x_0)$, for every $x\in B(x_0,R)$. Put
$A=\overline{B}(x_0,{R\over 2})\cup\{ x_1\}$.

In both cases, fix $p^*\in C([0,1],X)$ such that
$p^*(0)=x_0$ and $p^*(1)=x_1$. If $K^*=(p^*)^{-1}(A)$ then, by
Theorem \ref{1.2.9}, we obtain the existence of a critical point of
 $f$, which is different from $x_0$ and $x_1$, as we can easily deduce by
 examining the  proof of Theorem \ref{1.2.9}.
\qed

\medskip
With the same proof as of Corollary \ref{1.2.6} one can show

\medskip
\begin{cor}\label{1.2.11}
 Let $X$ be a Banach space and let
 \fxr be a
 G\ab teaux-differentiable functional such that the operator
$f':(X,\|\cdot\| )\rightarrow (X^*,\sigma (X^*,X))$ is continuous.
Assume that for every $p\in{\cal P}$ there exists
$t\in K\setminus K^*$ such that $f(p(t))\geq c$.

Then there exists a sequence $(x_n)$ in $X$ so that

i) $\ \ \displaystyle\limn f(x_n)=c\, ;$

ii) $\ \ \displaystyle\limn \| f'(x_n)\| =0\, .$

If, furthermore, $f$ satisfies \pslc , then there exists
$x\in X$ such that $f(x)=c$ and $f'(x)=0$.
\end{cor}

\medskip
The following result is a strengthened variant of
Theorems \ref{1.2.1} and \ref{1.2.9}.

\medskip
\begin{teo}\label{1.2.12}  Let \fxr$\ $ be a locally Lipschitz functional
and let
$F$ be a closed subset of $X$, with no common point with $p^*(K^*)$.
Assume that
\neweq{1.13}
f(x)\geq c,\quad\hbox{for every}\ \ x\in F
\endeq
and
\neweq{1.14}
p(K)\cap F\not= \emptyset \ ,\quad\hbox{for all}\, .
p\in{\cal P}
\endeq

Then there exists a sequence $(x_n)$ in $X$ such that

i)$\ \ \ \displaystyle\limn {\rm dist}\ (x_n,F)=0$;

ii)$\ \ \ \displaystyle \limn f(x_n)=c$;

iii)$\ \ \displaystyle\limn\lambda (x_n)=0\, .$
\end{teo}

\medskip
{\bf Proof.} Fix $\varepsilon >0$ such that
$$\varepsilon <\min\{ 1;\ {\rm dist}\ (p^*(K^*),F)\}\, .$$

Choose $p\in{\cal P}$ so that
$$\max_{t\in K}f(p(t))\leq c+{{\varepsilon^2}\over 4}\, .$$
The set
$$K_0=\{ t\in K;\ {\rm dist}\ (p(t),F)\geq\varepsilon\}$$
is bounded and contains $K^*$. Define
$${\cal P}_0=\{ q\in C(K,X);\ q=p\ \ \hbox{on}\ \ K_0\}\, .$$
Set
$$\eta :X\rightarrow \RR ,\quad \eta (x)=\max \{ 0;\
\varepsilon^2-
\varepsilon\  {\rm dist}\ (x,F)\}\, .$$
Define $\psi :{\cal P}_0\rightarrow\RR$ by
$$\psi (q)=\max_{t\in K}\ (f+\eta )(q(t))\, .$$
The functional $\psi$ is continuous
and bounded from below.
By  Ekeland's Principle, there exists $p_0\in{\cal P}_0$ such that,
for every $q\in{\cal P}_0$,
$$\psi (p_0)\leq \psi (q)\, ,$$
\neweq{1.16}
d(p_0,q)\leq{{\varepsilon}\over 2}\, ,
\endeq
\neweq{1.17}
\psi (p_0)\leq \psi (q)+{\varepsilon\over 2}\ d(q,p_0)
\, .
\endeq
The set
$$B(p_0)=\{ t\in K;\ \ (f+\eta )(p_0(t))=\psi (p_0)\}$$
is closed.
For concluding the proof, it is sufficient to show that there
exists
$t\in B(p_0)$ such that
\neweq{1.18}
\hbox{dist}\, (p_0(t),F)\leq{{3\varepsilon}\over 2}\, ,
\endeq
\neweq{1.19}
c\leq f(p_0(t))\leq c+{{5\varepsilon^2}\over 4}\, ,
\endeq
\neweq{1.20}
\lambda (p_0(t))\leq{{5\varepsilon}\over 2}.
\endeq
Indeed, it is enough to choose then
 $\varepsilon ={1\over n}$ and
$x_n=p_0(t)$.

\underbar{Proof of \eq{1.18}}: It follows by the definition of
${\cal P}_0$ and \eq{1.14} that, for every
$q\in{\cal P}_0$, we have
$$q(K\setminus K_0)\cap F\not= \emptyset\, .$$
Therefore, for any $q\in{\cal P}_0$,
$$\psi (q)\geq c+\varepsilon^2\, .$$
On the other hand,
$$\psi (p)\leq c+{{\varepsilon^2}\over 4}
+\varepsilon^2=c+{{5\varepsilon^2}\over 4}.$$
Hence
\neweq{1.21}
c+\varepsilon^2\leq\psi (p_0)\leq\psi (p)\leq c+
{{5\varepsilon^2}\over 4}\, .
\endeq
So, for each $t\in B(p_0)$,
$$c+\varepsilon^2\leq\psi (p_0)=(f+\eta )(p_0(t)).$$
Moreover, if $t\in K_0$, then
$$(f+\eta )(p_0(t))=(f+\eta )(p(t))=f(p(t))\leq c+{{\varepsilon^2}\over 4}.$$
This implies that
$$B(p_0)\subset K\setminus K_0.$$
By the definition of $K_0$ it follows that, for every $t\in B(p_0)$ we have
$$\hbox{dist}\ (p(t),F)\leq\varepsilon .$$
Now, the relation \eq{1.16} yields
$$\hbox{dist}\ (p_0(t),F)\leq{\varepsilon\over 2}.$$

\underbar{Proof of \eq{1.19}}: For every $t\in B(p_0)$ we have
$$\psi (p_0)=(f+\eta )(p_0(t)).$$
Using \eq{1.21} and taking into account that
$$0\leq\eta\leq\varepsilon^2,$$
it follows that
$$c\leq f(p_0(t))\leq c+{{5\varepsilon^2}\over 4}.$$

\underbar{Proof of \eq{1.20}}: Applying Lemma \ref{1.2.2} for
$\varphi (t)=\partial f(p_0(t))$, we find a continuous mapping
$v:B(p_0)\rightarrow X$ such that, for every $t\in B(p_0)$,
$$\|v(t)\|\leq 1.$$
Moreover, for every $t\in B(p_0)$ and $x^*\in \partial f(p_0(t))$,
$$\langle x^*,v(t)\rangle\geq\gamma -\varepsilon ,$$
where
$$\gamma =\inf_{t\in B(p_0)}\lambda (p_0(t)).$$
Hence for every $t\in B(p_0)$,
$$f^0(p_0(t),-v(t))=\max\{\langle x^*,-v(t)\rangle ;
 x^*\in\partial f(p_0(t))\}=$$
$$=-\min\{\langle x^*,v(t)\rangle ;
 x^*\in\partial f(p_0(t))\}\leq -\gamma +\varepsilon .$$

Since $B(p_0)\cap K_0=\emptyset$, there exists a continuous extension $w$
of $v$ to the set $K$ such that $w=0$ on $K_0$ and
$\|w(t)\|\leq 1$, for all $t\in K$.

Now, by \eq{1.17} it follows that for every $\lambda >0$,
\neweq{1.22}
-{\varepsilon\over 2}\leq -{\varepsilon\over 2}\ \| w\|_{\infty}
\leq{{\psi (p_0-\lambda w)-\psi (p_0)}\over{\lambda}}\, .
\endeq
For every $n$, there exists $t_n\in K$ such that
$$\psi (p_0-{1\over n}w)=(f+\eta )(p_0(t_n)-{1\over n}w(t_n)).$$
Passing eventually to a subsequence, we may suppose that
 $(t_n)$ converges to
$t_0$, which, clearly, lies in $B(p_0)$. On the other hand, for every
 $t\in K$ and $\lambda >0$ we have
$$f(p_0(t)-\lambda w(t))\leq f(p_0(t))+\lambda\varepsilon .$$
Hence
$$n[\psi (p_0-\lambda w)-\psi (p_0)]\leq n[f(p_0(t_n)-{1\over n}w(t_n))+
{\varepsilon\over n}-f(p_0(t_n))] .$$
Therefore, by
\eq{1.22} it follows that
$$-{{3\varepsilon}\over 2}\leq n[\psi (p_0(t_n)-{1\over n}w(t_n))-f(p_0(t_n))]\leq$$
$$\leq n[\psi (p_0(t_n)-{1\over n}w(t_0))-f(p_0(t_n))]+$$
$$+n[f(p_0(t_n)-{1\over n}w(t_n))-f(p_0(t_n)-{1\over n}w(t_0)).$$

Using the fact that $f$ is locally Lipschitz and
$t_n\rightarrow t_0$ we find
$$\limsup_{n\rightarrow\infty}n\,
 [f(p_0(t_n)-{1\over n}w(t_n))-f(p_0(t_n)-{1\over n}
w(t_0))]=0.$$
Therefore
$$-{{3\varepsilon}\over 2}\leq\limsup_{n\rightarrow\infty}n\,
[f(p_0(t_0)+z_n-{1\over n}w(t_0))-
f(p_0(t_0)+z_n)],$$
where $z_n=p_0(t_n)-p_0(t_0)$. Hence
$$-{{3\varepsilon}\over 2}\leq f^0(p_0(t_0),-w(t_0))\leq
 -\gamma +\varepsilon .$$
So
$$\gamma =\inf\{\|x^*\| ; x^*\in\partial f(p_0(t)),t\in B(p_0)\}
\leq {{5\varepsilon}\over 2}.$$

Now, by the lower semicontinuity of $\lambda$, we find
 $t\in B(p_0)$ such that
$$\lambda (p_0(t))=\inf_{x^*\in\partial f(p_0(t))}\| x^*\|
 \leq {{5\varepsilon}\over 2},$$
which ends our proof. \qed

\medskip
\begin{cor}\label{1.2.13}
 Under hypotheses of Theorem \ref{1.2.12}, if
  $f$ satisfies,
$({\rm PS})_c$, then $c$ is a critical value of
 $f$.
 \end{cor}

\medskip
\begin{rem}\label{remi} If
$$\inf_{x\in X_1}f(x+z)=\max_{x\in K^*}f(x),$$
then the conclusion of Corollary
 \ref{1.2.7} remains valid, with an argument based on Theorem
 \ref{1.2.11}.
 \end{rem}

\medskip
\begin{cor}\label{1.2.14} (Ghoussoub-Preiss Theorem).
Let $f:X\rightarrow\RR$
be a lo\-cally Lip\-s\-chitz G\^ateaux-dif\-fer\-en\-ti\-a\-ble functional such that
$f':(X,\|\cdot\| )\rightarrow (X^*,\sigma (X^*,X))$ is continuous. Let $a$
and $b$ be in $X$ and define
$$c=\inf_{p\in{\cal P}}\max_{t\in [0,1]}f(p(t)),$$
where ${\cal P}$ is the set of continuous paths $p:[0,1]\rightarrow X$ such
that $p(0)=a$ and $p(1)=b$.
 Let $F$ be a closed subset of $X$ which does not contain $a$ and
$b$ and  $f(x)\geq c$, for all $x\in F$. Suppose, in addition that, for every
 $p\in{\cal P}$,
$$p([0,1])\cap{\cal P}\not =\emptyset .$$

Then there exists a sequence $(x_n)$ \ii n $X$ so that

i)$\ \ \displaystyle\lim_{n\rightarrow\infty}{\rm dist}(x_n,F)=0\, ;$

ii)$\ \displaystyle\lim_{n\rightarrow\infty}f(x_n)=c\, ;$

iii)$\ \displaystyle\lim_{n\rightarrow\infty}\| f'(x_n)\| =0\, .$

Moreover, if $f$ satisfies  $({\rm PS})_c$, then there exists
$x\in F$ such that $f(x)=c$ and $f'(x)=0$.
\end{cor}

\medskip
{\bf Proof.} With the same arguments as in the proof of
Corollary
\ref{1.2.6} we deduce that the functional $f$ is locally Lipschitz
and
$$\partial f(x)=\{f'(x)\} .$$

Applying Theorem \ref{1.2.12} for $K=[0,1]$, $K^*=\{0,1\}$, $p^*(0)=a$
 $p^*(1)=b$, our conclusion follows.

The last part of the theorem follows from Corollary \ref{1.2.13}.
\qed

\section{Applications of the Mountain-Pass Theorem}
\subsubsection{The simplest model}
Let $1<p<\frac{N+2}{N-2}$, if $N\geq 3$, and $1<p<+\infty$, provided that
$N=1,2$. Consider the problem
\neweq{model}
\itab
$-\Delta u=u^p\, ,$ & \qquad $\mbox{in}\ \, \Omega $\\
$ u>0\, ,$ & \qquad $\mbox{in}\ \, \Omega $\\
$ u=0\, ,$ & \qquad $\mbox{on}\, \ \partial\Omega\, .$\\
\ttab
\endeq

 Our aim is to prove in what follows the following
 \begin{teo}\label{teomodel}
There exists a solution of the problem \eq{model}, which is
not necessarily unique. Furthermore, this solution is unstable.
\end{teo}

\begin{rem}\label{comm}
If $p=\frac{N+2}{N-2}$ then the  energy functional associated to the
problem \eq{model} does not have the Palais-Smale property. The case
$p\geq\frac{N+2}{N-2}$ is difficult; for instant, there is no solution
even in the simplest case where $\Omega =B(0,1)$. If $p=1$ then the
existence of a solution depends on the geometry of the domain: if 1 is
not an eigenvalue of $(-\Delta )$ in $H_0^1(\Omega )$ then there is no
solution to our problem \eq{model}. If $0<p<1$ then there exists a unique solution
(since the mapping $u\longmapsto f(u)/u=u^{p-1}$ is decreasing) and,
moreover, this solution is stable. The arguments may be done in this case
by using the method of sub and super solutions.
\end{rem}

{\bf Proof}. We first argue the instability of the solution. So, in order
to justify that $\lambda_1(-\Delta -pu^{p-1})<0$, let $\varphi$ be
an eigenfunction corresponding to $\lambda_1$. We have
$$-\Delta\varphi -pu^{p-1}\varphi =\lambda_1\varphi\, ,
\qquad\mbox{in}\ \, \Omega\, .$$
Integrating by parts this equality we find
$$(1-p)\intom u^p\varphi =\lambda_1\intom\varphi u\, ,$$
which implies $\lambda_1<0$, since $u>0$ in $\Omega$ and $p>1$.

We will prove the existence of a solution by using two different methods:

\medskip
{\bf 1. A variational proof}. Let
$$m=\inf\{ \intom \mid\nabla v\mid^2;\, v\in H^1_0(\Omega )\
\mbox{and}\ \| v\|_{L^{p+1}}=1\}\, .$$
{\it First step: $m$ is achieved}. Let $(u_n)\subset H^1_0(\Omega )$
be a minimizing sequence.
Since $p<\frac{N+2}{N-2}$ then $H^1_0(\Omega )$ is
compactly embedded in $L^{p+1}(\Omega $. It follows that
$$\intom \mid\nabla u_n\mid^2\ri m\, ,$$
$$\| u_n\|_{L^{p+1}}=1\, .$$
So, up to a subsequence,
$$u_n\rightharpoonup u\, ,\qquad\mbox{weakly in}\ \, H^1_0(\Omega )$$
and
$$u_n\ri u\, ,\qquad\mbox{strongly in}\ \, L^{p+1}(\Omega )\, .$$
By the  lower semicontinuity of the functional
$\|\cdot\|_{L^2}$ we find that
$$\intom\mid\nabla u\mid^2\leq\liminf_{n\ri\infty}\intom\mid\nabla u_n\mid^2=m$$
which implies $\intom \mid\nabla u\mid^2=m$. Since $\| u\|_{L^{p+1}}=1$,
it follows that $m$ is achieved by $u$.

We remark that we have even $u_n\ri u$, strongly in $H^1_0(\Omega )$. This follows by
the weak convergence of $(u_n)$ in $H^1_0(\Omega )$ and by
$\| u_n\|_{H^1_0}\ri\| u\|_{H^1_0}$.

\medskip
{\it 2. $u\geq 0$, a.e. in $\Omega$}. We may assume that $u\geq 0$, a.e. in
$\Omega$. Indeed, if not, we may replace $u$ by $\mid u\mid$. This is possible
since $\mid u\mid\in H^1_0(\Omega )$ and so, by Stampacchia's theorem,
$$\nabla \mid u\mid =(\mbox{sign}\, u)\, \nabla u\, ,\qquad\mbox{if}\ u\not=0
\, .$$
Moreover, on the level set $[u=0]$ we have $\nabla u=0$, so
$$\mid\nabla\mid u\mid \mid =\mid\nabla u\mid\, ,\qquad\mbox{a.e. in}\
\, \Omega\, .$$

\medskip
{\it 3. $u$ verifies $-\Delta u=u^p$ in weak sense}. We have to prove that, for every
$w\in H^1_0(\Omega )$,
$$\intom\nabla u\nabla w =m\intom u^pw\, .$$
Put $v=u+\ep w$ in the definition of $m$. It follows that
$$\intom\mid\nabla v\mid^2=\intom\mid\nabla u\mid^2+2\ep\, \intom \nabla u
\nabla w+\ep^2\,\intom\mid\nabla w\mid^2$$
and
$$
\begin{array}{lll}
\intom\mid u+\ep w\mid^{p+1} = & \intom\mid u\mid^{p+1}+\ep (p+1)\, \intom
u^pw+o(\ep )=\\
&1+\ep (p+1)\, \intom
u^pw+o(\ep )\, .
\end{array}
$$
Therefore
$$\| v\|^2_{L^{p+1}}=\left(1+\ep (p+1)\, \intom
u^pw+o(\ep )\right)^{2/(p+1)}=1+2\ep\, \intom u^pw+o(\ep )\, .$$
Hence
$$m=\intom\mid u\mid^2\leq
\frac{m+2\ep\, \intom\nabla u\nabla w+o(\ep )}{1+2\ep\, \intom u^pw+o(\ep )}
=m+2\ep\, \left(\intom\nabla u\nabla w-m\intom u^pw\right)+o(\ep )\, ,$$
which implies
$$\intom\nabla u\nabla w=m\intom u^pw\, ,\qquad\mbox{for every}\ \, w\in
H^1_0(\Omega )\, .$$
Consequently, the function $u_1=m^\alpha u$ ($\alpha =\frac{1}{p-1}$) is
a weak solution of our problem \eq{model}, that is $u=u_1m^{-\alpha}$ is weak
solution to \eq{model}.

\medskip
{\it 4. Regularity of $u$}. We know until now that $u\in
H^1_0(\Omega )\subset L^{2^\star}(\Omega )$. In a general
framework, assuming that $u\in L^q$, it follows that $u^p\in
L^{q/p}$, that is, by Schauder regularity and Sobolev embeddings,
$u\in W^{2,q/p}\subset L^s$, where
$\frac{1}{s}=\frac{p}{q}-\frac{2}{N}$. So, assuming that
$q_1>(p-1)\, \frac{N}{2}$, we have $u\in L^{q_2}$, where
$\frac{1}{q_2}=\frac{p}{q_1}-\frac{2}{N}$. In particular,
$q_2>q_1$. Let $(q_n)$ be the increasing sequence we may construct
in this manner and set $q_\infty =\lim_{n\ri\infty}q_n$. Assuming,
by contradiction, that $q_n<\frac{Np}{2}$ we obtain, passing at
the limit as $n\ri\infty$, that $q_\infty =\frac{N(p-1)}{2}<q_1$,
contradiction. This shows that there exists $r>\frac{N}{2}$ such
that $u\in L^r(\Omega )$ which implies $u\in W^{2,r}(\Omega
)\subset L^\infty (\Omega )$. Therefore $u\in W^{2,r}(\Omega
)\subset C^k(\overline\Omega )$, where $k$ denotes the integer
part of $2-\frac{N}{r}$. Now, by H\"older continuity, $u\in
C^2(\overline\Omega )$.

\bigskip
{\bf 2. Second proof (Mountain-Pass Lemma)}. Set
$$F(u)=\frac{1}{2}\, \intom\mid\nabla u\mid^2-\frac{1}{p+1}\,
\intom (u^+)^{p+1}\, ,\qquad u\in H^1_0(\Omega )\, .$$
Standard arguments show that $F$ is a $C^1$ functional and $u$ is a critical
point of $F$ if and only if $u$ is a solution to the problem \eq{model}.
We observe that $F'(u)=-\Delta u-(u^+)^p\in H^{-1}(\Omega )$. So, if $u$ is
a critical point of $F$ then
$-\Delta u=(u^+)^p\geq 0$ in $\Omega$ and hence, by the Maximum Principle,
$u\geq 0$ in $\Omega$.

We verify the hypotheses of the Mountain-Pass Lemma. Obviously,
$F(0)=0$. On the other hand,
$$
\intom (u^+)^{p+1}\leq\intom \mid u\mid^{p+1}=\| u\|^{p+1}_{L^{p+1}}\leq
C\| u\|^{p+1}_{H^1_0}\, .
$$
Therefore
$$F(u)\geq\frac{1}{2}\, \| u\|^2_{H^1_0} -\frac{C}{p+1}\, \| u\|^{p+1}_{H^1_0}\geq
\rho >0\, ,$$
provided that $\| u\|_{H^1_0}=R$, small enough.

Let us now prove the existence of some $v_0$ such that $\| v_0\| >R$ and
$F(v_0)\leq 0$. For this aim, choose an arbitrary $w_0\geq 0$, $w_0\not\equiv 0$. We have
$$F(tw_0)=\frac{t^2}{2}\, \intom\mid\nabla w_0\mid^2-\frac{t^{p+1}}{p+1}\,
\intom (w_0^+)^{p+1}\leq 0\, ,$$ for $t>0$ large enough. \qed

\subsubsection{A bifurcation problem}

Let us consider a $C^1$ convex function $f:\RR\ri\RR$ such that $f(0)>0$ and
$f'(0)>0$. We also assume that there exists $1<p<\frac{N+2}{N-2}$ such that
$$\mid f(u)\mid\leq C(1+\mid u\mid^p)$$
and there exist $\mu >2$ and $A>0$ such that
$$\mu \int_0^u f(t)dt\leq uf(u)\, ,\qquad \mbox{for every}\ \, u\geq A\, .$$
A standard example of function satisfying these conditions is $f(u)=(1+u)^p$.

Consider the problem
\neweq{bifi}
\itab
$\di -\Delta u=\lambda f(u)\, ,$ & \qquad $\mbox{in}\ \, \Omega $\\
$ \di u>0\, ,$ & \qquad $\mbox{in}\ \, \Omega $\\
$\di u=0\, ,$ & \qquad $\mbox{on}\, \ \partial\Omega\, .$\\
\ttab
\endeq

We already know that there exists $\lambda^\star >0$ such that, for every
$\lambda <\lambda^\star$,
there exists a minimal and stable solution $\underline u$ to the problem \eq{bifi}.

\medskip
\begin{teo}\label{msol}
Under the above hypotheses on $f$, for every $\lambda \in (0,\lambda^\star )$, there exists
a second solution $u\geq\underline u$ and, furthermore, $u$ is unstable.
\end{teo}

\medskip
{\bf Proof}. We find a solution $u$ of the form $u=\underline u+v$ with $v\geq 0$. It
follows that $v$ satisfies
\neweq{v}
\itab
$\di -\Delta v=\lambda (f(\underline u+v)-f(\underline u)\, ,$ &
\qquad $\mbox{in}\ \, \Omega $\\
$\di v>0\, ,$ &\qquad $\mbox{in}\ \, \Omega $\\
$\di v=0\, ,$ & \qquad $\mbox{on}\, \ \partial\Omega\, .$\\
\ttab
\endeq
Hence $v$ fulfills an equation of the form
$$-\Delta v+a(x)v=g(x,v)\, ,\qquad\mbox{in}\ \,\Omega\, ,$$
where $a(x)=-\lambda f'(\underline u)$ and
$$g(x,v)=\lambda\left( f(\underline u(x)+v)-f(\underline u(x))\right)-
\lambda f'(\underline u(x))v\, .$$
We verify easily the following properties:

\noindent (i) $g(x,0)=g_v(x,0)=0$;

\noindent (ii) $\mid g(x,v)\mid\leq C(1+\mid v\mid^p)$;

\noindent (iii) $\mu\int_0^vg(x,t)dt\leq vg(x,v)$, for every $v\geq A$ large enough;

\noindent (iv) the operator $-\Delta -\lambda f'(\underline u)$ is coercive, since
$\lambda_1(-\Delta -\lambda f'(\underline u))>0$, for every $\lambda <\lambda_1$.

So, by the Mountain-Pass Lemma, the problem \eq{bifi} has a
solution which is, {\it a fortiori}, unstable. \qed

\section{Critical points and coerciveness of locally Lipschitz
functionals with the strong Palais-Smale property}

 The Palais-Smale property for $C^1$ functionals appears as the
most natural
compactness condition. In order to obtain corresponding results for
non-differentiable functionals the Palais-Smale condition introduced in
Definition \ref{1.1.5} is not always  an efficient tool, because of the nonlinearity
of the Clarke subdifferential. For this aim, we shall define a stronger
Palais-Smale type condition, which will be very useful in applications. In many
cases, our compactness condition will be a local one, similar to
 $(\hbox{PS})_c$ in Definition \ref{1.1.5}. The most efficient tool in our reasonings
 will be, as in the preceding paragraph, the Ekeland variational principle. As
 we shall remark the new Palais-Smale condition is in closed link with
 coerciveness properties of locally Lipschitz functionals.

As above, $X$ will denote a real Banach space.

\medskip
\begin{defin}\label
{1.3.1} \sl The locally Lipschitz functional
$f:X\rightarrow \RR$ is said to satisfy the strong Palais-Smale condition at
the point
 $c$ (notation: $(\hbox{ s-PS})_c$) provided that, for every sequence
 $(x_n)$ in $X$ satisfying
\neweq{1.23}
\lim_{n\rightarrow\infty}f(x_n)=c
\endeq
and
\neweq{1.24}
f^0(x_n,v)\geq -{1\over n}\, \| v\|, \ \ \hbox{ for every}\
\ v\in X,
\endeq
contains a convergent subsequence.

If this property holds for any real number $c$ we shall say that $f$ satisfies
the strong Palais-Smale condition
 $({\rm s-PS})$).\rm
 \end{defin}

\medskip
\begin{rem}\label
{1.3.2} It follows from the continuity of $f$ and the upper
semicontinuity of $f^0(\cdot ,\cdot )$ that if $f$ satisfies the
condition \spsc and there exist sequences $(x_n)$ and
$(\varepsilon_n)$ such that the conditions \eq{1.23} and \eq{1.24}
 are fulfilled, then $c$ is a critical value of
  $f$. Indeed, up to a subsequence,
 $x_n\rightarrow x$. It follows that
 $f(x)=c$ and, for every $v\in X$,
$$f^0(x,v)\geq\limsup_{n\rightarrow\infty}f^0(x_n,v)\geq 0,$$
that is  $0\in\partial f(x)$.
\end{rem}

\medskip
\begin{defin}\label
{1.3.3} \sl A mapping $f:X\rightarrow\RR$ is said to be
coercive provided that
$$\lim_{\| x\|\rightarrow\infty}f(x)=+\infty .$$\rm
\end{defin}

\medskip
For each $a\in\RR$ and $f:X\rightarrow\RR$, we shall denote from now on
$$[f=a]=\{x\in X; f(x)=a\}\, ;$$
$$[f\leq a]=\{x\in X; f(x)\leq a\}\, ;$$
$$[f\geq a]=\{x\in X; f(x)\geq a\}\, .$$

\medskip
\begin{prop}\label
{1.3.4} \sl Let $f:X\rightarrow\RR$ be a locally Lipschitz bounded
from below functional.  If $a=\displaystyle\inf_{X}f$ and $f$
satisfies the condition \pslie, then there exists
 $\alpha >0$ such that the set
 $[f\leq a+\alpha ]$ is bounded.\rm
 \end{prop}

\medskip
{\bf Proof.} We assume, by contradiction, that for every
$\alpha >0$, the set $[f\leq a+\alpha ]$ is unbounded.
So, there exists a sequence
 $(z_n)$ in $X$ such that, for every $n\geq 1$,
$$a\leq f(z_n)\leq a+{1\over{n^2}}\, ,$$
$$\| z_n\|\geq n.$$
Using Ekeland's Principle, for every
 $n\geq 1$ there is some $x_n\in X$
such that, for any $x\in X$,
$$a\leq f(x_n)\leq f(z_n)\, ,$$
$$f(x)-f(x_n)+{1\over n}\ \| x-x_n\|\geq 0\, ,$$
$$\| x_n-z_n\|\leq{1\over n}\, .$$
Therefore
$$\| x_n\|\geq n-{1\over n}\longrightarrow\infty\, ,$$
$$f(x_n)\longrightarrow a\, ,$$
and, for each $v\in X$,
$$f^0(x_n,v)\geq -{1\over n}\ \| v\|\, .$$
Now, by
 \pslie , it follows that the unbounded sequence
$(x_n)$ contains a convergent subsequence, contradiction. \qed

\medskip
The following is an immediate consequence of the above result

\medskip
\begin{cor}\label
{1.3.5} \sl If $f$ is a locally Lipschitz bounded from below
functional satisfying the strong Palais-Smale condition, then
 $f$ is coercive.\rm
 \end{cor}

\medskip
This result was proved by S.J. Li \cite{Li} for $C^1$ functionals. He used in his
proof the Deformation Lemma. Corollary \ref{1.3.5}
 also follows from

\medskip
\begin{prop}\label
{1.3.6} \sl Let $f:X\rightarrow\RR$ be a locally Lipschitz
functional satisfying
$$a=\liminf_{\| x\|\rightarrow\infty}f(x)<+\infty .$$
Then there exists a sequence $(x_n)$ in $X$ such that
$$\| x_n\|\rightarrow\infty,\quad f(x_n)\longrightarrow a$$
and, for every $v\in X$,
$$f^0(x_n,v)\geq -{1\over n}\, \| v\| .$$\rm
\end{prop}

\medskip
{\bf Proof.} For every $r>0$ we define
$$m(r)=\inf_{\| x\| \geq r}f(x).$$
Obviously, the mapping $m$ is non-decreasing
and $\displaystyle\lim_{r \rightarrow\infty}m(r)=a$.
For any integer $n\geq 1$ there exists $r_n>0$ such that, for every
$r\geq r_n$,
$$m(r)\geq a-{1\over{n^2}}.$$

Remark that we can choose $r_n$ so that $r_n\geq{ n\over\scriptstyle 2}
+{1\over\scriptstyle n}$.
Choose $z_n\in X$ such that $\| z_n\|\geq 2r_n$ and
\neweq{1.25}
f(z_n)\leq m(r_{2n})+{1\over{n^2}}\leq a+{1\over{n^2}}.
\endeq
Applying Ekeland's Principle to the functional $f$ restricted to
the set $\{x\in X;\ \| x\|\geq r_n\}$ and for
 $\varepsilon ={1\over n}$, $z=z_n$,
we get $x_n\in X$ \ai $\| x_n\|\geq r_n$ and, for every $x\in X$ with
$\| x\|\geq r_n$,
\neweq{1.26}
f(x)\geq f(x_n)-{1\over n}\ \| x-x_n\|\, ,
\endeq
\neweq{1.27}
a-{1\over{n^2}}\leq m(r_n)\leq f(x_n)\leq f(z_n)-
{1\over n}\ \| z_n-x_n\|\, .
\endeq

It follows from \eq{1.25} and \eq{1.27}
 that $\| x_n-z_n\|\leq{2\over\scriptstyle n}$, which implies
$$\| x_n\| \geq 2r_n-{2\over n}\longrightarrow +\infty .$$
On the other hand, $f(x_n)\longrightarrow a$.
For every $v\in X$ and $\lambda >0$, putting $x=x_n+\lambda v$ in
\eq{1.26}, we find
$$f^0(x_n,v)\geq\limsup_{\lambda\searrow 0}
{{f(x_n+\lambda v)-f(x_n)}\over\lambda}\geq -{1\over n}\ \| v\|\,
,$$ which concludes our proof. \qed

\medskip
\begin{prop}\label
{1.3.7} \sl Let \fxr be a locally Lipschitz bounded from below
functional. Assume there exists $c\in\RR$ such that $f$ satisfies
 \spsc and, for every $a<c$, the set $[f\leq a]$ is bounded.

Then there exists $\alpha >0$ such that the set
 $[f\leq c+\alpha ]$ is bounded.\rm
 \end{prop}

\medskip
{\bf Proof.} Arguing by contradiction, we assume that the set
$[f\leq c+\alpha ]$ is unbounded, for every $\alpha >0$. It
follows by our hypothesis that, for every
 $n\geq 1$, there exists $r_n\geq n$ such that
$$[f\leq c-{1\over{n^2}}]\subset B(0,r_n).$$
Set
$$c_n=\inf_{X\setminus B(0,r_n)}f\geq c-{1\over{n^2}}.$$
Since the set $[f\leq c+{1\over{n^2}}]$ is unbounded, we obtain the existence
of a sequence
$(z_n)$ in $X$ such that
$$\| z_n\|\geq r_n+1+{1\over n}$$
and
$$f(z_n)\leq c+{1\over{n^2}}.$$
So, $z_n\in X\setminus B(0,r_n)$ and
$$f(z_n)\leq c_n+{2\over{n^2}}.$$
Applying Ekeland's Principle to the functional $f$ restricted to
$X\setminus B(0,r_n)$, we find
$x_n\in X\setminus B(0,r_n)$ such that, for every $x\in X$ with
 $\| x\|\geq r_n$, we have
$$c_n\leq f(x_n)\leq f(z_n)\, ,$$
$$f(x)\geq f(x_n)-{2\over n}\ \| x-x_n\|\, ,$$
$$\| x_n-z_n\|\leq{1\over n}\, .$$
Hence
$$\| x_n\|\geq \| z_n\| -\| x_n-z_n\|\geq
r_n+1\longrightarrow +\infty\, ,$$
$$f(x_n)\longrightarrow c$$
and, for every $v\in X$,
$$f^0(x_n,v)\geq -{2\over n}\, \| v\| \, .$$
Now, by \spsc , we obtain that the unbounded sequence $(x_n)$
contains a convergent subsequence, contradiction. \qed

\medskip
\begin{prop}\label
{1.3.8} \sl Let \fxr be a locally Lipschitz bounded from below
functional. Assume that $f$ is not coercive. If
$$c=\sup\{ a\in\RR ;\ [f\leq a]\ \ \hbox{is bounded}\} ,$$
then $f$ does not satisfy the condition \spsc .\rm
\end{prop}

\medskip
{\bf Proof.} Denote
$$A=\{ a\in\RR ;\ [f\leq a]\ \ \hbox{is bounded}\} .$$
It follows from the lower boundedness of
 $f$ that the set $A$ is nonempty.
Since  $f$ is not coercive, it follows that
$$c=\sup A<+\infty .$$

Assume, by contradiction, that $f$  satisfies \spsc\, .
 Then, by Proposition
 \ref{1.3.7}, there exists $\alpha >0$ such that the set
$[f\leq a+\alpha ]$ is bounded, which contradicts the maximality
of $c$. \qed

\medskip
\begin{rem}\label{1.3.9}
 The real number $c$ defined in Proposition
\ref{1.3.8} may  also be characterized by
$$c=\inf\{ b\in\RR ;\,\,[f\leq b]\,\,\hbox{is unbounded}\} \, .$$
\end{rem}

\medskip
\begin{prop}\label{1.3.10}  Let \fxr be a locally Lipschitz functional
satisfying (s-PS). Assume there exists $a\in\RR$ such that the set
 $[f\leq a]$ is bounded.

Then the functional $f$ is coercive.\rm
\end{prop}

\medskip
{\bf Proof.} Without loss of generality, we may assume that
 $a=0$. It follows now from our hypothesis that there exists an integer
$n_0$ such that $f(x)>0$,for any $x\in X$ with $\| x\|\geq n_0$. We assume, by
contradiction, that
$$
0\leq c = \liminf_{\| x\|\rightarrow\infty} f(x)<+\infty .
$$
Applying Proposition \ref{1.3.6}, we find a sequence $(x_n)$ in $X$ such that
$\| x_n\|\rightarrow\infty $, $ f(x_n)\longrightarrow c$ and, for every
 $v\in X$,
$$f^{0}(x_n,v)\geq{}-{1\over n}\ \| v\| .$$
Using now the condition \psg , we obtain that the unbounded sequence
 $(x_n)$
contains a convergent subsequence, contradiction. So, the
functional $f$ is coercive. \qed

\medskip
\begin{cor}\label{1.3.11} Let \fxr be a locally Lipschitz bounded from
below functional which satisfies
\psg .

Then every minimizing sequence of $f$ contains a convergent subsequence.
\end{cor}

\medskip
{\bf Proof.} Let $(x_n)$ be a minimizing sequence of $f$.
Passing eventually at a subsequence we have
$$f(x_n)\leq\mbox{$\displaystyle\inf_{X} {f+{1\over{n^2}}}$} .$$
By Ekeland's Principle, there exists $z_n\in X$ such that, for every
 $x\in X$,
$$f(x)\geq f(z_n)-{1\over n}{}\| x-z_n\|\, ,$$
$$  f(z_n)\leq f(x_n)-{1\over n}{}\| x_n-z_n\| \, .$$
With an argument similar to that used in the proof of Proposition
 \ref{1.3.6}
we find
\neweq{1.28}
\| x_n-z_n\|\leq{2\over n}\, ,
\endeq
$$f(z_n)\leq\inf_{X}f+{1\over{n^2}}\, ,
$$
and, for every $v\in X$,
$$f^0(z_n,v)\geq -{1\over n}\ \| v\| \, .$$
Using now the condition \psg ,we find that the sequence
 $(z_n)$ is relatively compact.
By \eq{1.28} it follows that the corresponding subsequence of
 $(x_n)$ is convergent, too.
\qed

\medskip
Define the map
$$M:[0,+\infty )\rightarrow\RR,\quad\ \ M(r)=\inf_{\| x\| =r}f(x).$$
We shall prove in what follows some elementary properties of this functional.

\medskip
\begin{prop}\label{1.3.12}  Let \fxr be a locally Lipschitz bounded from
below functional which satisfies the condition
 (s-PS). Assume there exists $R>0$ such that all the critical points of
  $f$ are in the closed ball of radius  $R$.

Then the functional $M$ is increasing and continuous at the right on the set
 $(R,+\infty )$.\end{prop}

\medskip For the proof of this result an auxiliary one. First we introduce a
weaker variant of the condition
 \psg for functionals defined on a circular crown.

\medskip
\begin{defin}\label{1.3.13}  Let $0<a<b$ and let  $f$ be a locally Lipschitz
map defined on
$$A=\{ x\in X;\ a\leq\| x\|\leq b\} .$$

We say that $f$ satisfies the Palais-Smale type condition
 $({\rm PS})_A$ provided that
every sequence $(x_n)$ satisfying
$$a+\delta\leq\| x_n\|\leq b-\delta,\quad\hbox{for some}\
\ \delta >0\, ,$$
$$\sup_{n}|f(x_n)|<+\infty\, ,$$
$$f^0(x_n,v)\geq -{1\over n}\,
\| v\| ,\quad\hbox{for some}\ \ v\in X,$$
contains a convergent subsequence.
\end{defin}

\medskip
\begin{lemma}\label{1.3.14} \sl Let $A$ be as in Definition
\ref{1.3.13} and let
 $f$ be a locally Lipschitz bounded from below functional defined on
 $A$.
If $f$ satisfies $({\rm PS})_A$ and $f$ does not have critical points which are
interior point of
 $A$ then, for every $a<r_1<r<r_2<b$,
\neweq{1.29}
M(r)>\min\{ M(r_1),M(r_2)\}\, .
\endeq\rm
\end{lemma}

\medskip
{\bf Proof of Lemma.} Without loss of generality, let us assume that
 $f$ takes only positive values. Arguing by contradiction, let $r_1<r<r_2$ be
such that the inequality \eq{1.29} is not fulfilled.
There exists a sequence $(x_n)$
such that $\| x_n\| =r$ and
$$f(x_n)<M(r)+{1\over{n^2}}.$$
Applying now Ekeland's Principle to $f$ restricted to the set
$$B=\{ x\in X;\ r_1\leq \| x\|\leq r_2\} ,$$
we find $z_n\in B$ such that, for every $x\in B$,
$$f(x)\geq f(z_n)-{1\over n}{}\| x-z_n\|\, ,$$
$$f(z_n)\leq f(x_n)-{1\over n}{}\| x_n-z_n\| \, .$$
Moreover, $\ r_1<\| z_n\| <r_2$, for $n$ large enough.
Indeed, if it would exist
 $n\geq 1$ such that $\|z_n\| =r_1$, then
$$M(r_1)\leq f(x_n)-{1\over n}{}\| x_n-z_n\|
\leq M(r)+{1\over{n^2}}-{1\over n}{}\| x_n-z_n\|\leq$$
$$\leq M(r)+{1\over{n^2}}-{1\over n}(r-r_1)\leq M(r_1)+
{1\over{n^2}}-{1\over n}(r-r_1).$$
It follows that
 $\ r-r_1\leq {1\over n}$, which is not possible if  $n$ is sufficiently large.
Therefore
$$\sup_{n}|f(z_n)|=\sup_{n}f(z_n)\leq M(r)$$
and, for every $v\in X$,
$$f^0(z_n,v)\geq -{1\over n}\, \| v\| \, .$$

Using now $({\rm PS})_A$, the sequence $(z_n)$ contains a subsequence which
converges to a critical point of
 $f$ belonging to $B$. This contradicts one of the hypotheses imposed to $f$.
\qed

\medskip
{\bf Proof of Proposition \ref{1.3.12}} If $M$ is not  increasing, there exists
$r_1<r_2$ such that $M(r_2)\leq M(r_1)$. On the other hand, by Corollary
\ref{1.3.5} we have
$$\lim_{r\rightarrow\infty}M(r)=+\infty .$$
Choosing now $r>r_2$ so that $M(r)\geq M(r_1)$, we find that $r_1<r_2<r$ and
$$M(r_2)\leq M(r_1)=\min\{M(r_1),M(r)\} ,$$
which contradicts Lemma \ref{1.3.14}. So, $M$ is an increasing map.

The continuity at the right of $M$ follows from its upper
semicontinuity. \qed

\medskip
There exists a local variant of Proposition
 \ref{1.3.12} for locally Lipschitz functionals defined on the set
  $\ \{ x\in X;\ \| x\|\leq R_0\}$, for some $R_0>0$.

Assume $f$ satisfies the condition \psg in the following sense:
every sequence $(x_n)$ with the properties
$$\| x_n\|\leq R<R_0\, ,$$
$$\sup_{n}|f(x_n)|<+\infty$$
and
$$f^0(x_n,v)\geq -{1\over n}\, \| v\| ,
\quad\hbox{for every}\ \ v\in X\ \ \hbox{and}\ \ n\geq 1,$$
is relatively compact.

\medskip
\begin{prop}\label{1.3.15} \sl Let $f$ be a locally Lipschitz functional defined
on $\| x\|\leq R_0$ and satisfying the condition \psg .
Assume $f(0)=0$, $f(x)>0$ provided $0<\| x\| <R_0$ and
$f$ does not have critical point in the set $\{ x\in X;\ 0<\| x\| <R_0\}$.

Then there exists $0<r_0\leq R_0$ such that $M$ is increasing on $[0,r_0)$
and decreasing on $[r_0,R_0)$.\rm
\end{prop}

\medskip
{\bf Proof.} Let $(R_n)$ an increasing sequence of positive numbers
which converges to $R_0$.
It follows by the upper semicontinuity of $M$ restricted to
 $[0,R_n]$ that there exists
$r_n\in (0,R_n]$ such that $M$ achieves its maximum in $r_n$. Let
$r_0\in (0,R_0]$ be the limit of the increasing sequence $(r_n)$.
Our conclusion follows now easily by applying Lemma \ref{1.3.14}.
\qed

\chapter{Critical point theorems of
 Lusternik-Schnirelmann type}

\section{Basic results on the notions of genus and
 Lusternik-Schnirelmann category}

\hspace*{6mm} One of the most interesting problems related to the extremum problems is
how to find estimates of the eigenvalues and eigenfunctions of a given operator. In this
field the Lusternik-Schnirelmann theory plays a very important role. The starting point
of this theory is the eigenvalue problem
$$Ax=\lambda x,\quad \lambda\in\RR ,\ \ x\in\RR^n ,$$
where $A\in M_{n}(\RR )$ is a symmetric matrix.
This problem may be written, equivalently,
$$F'(x)=\lambda x,\quad \lambda\in\RR ,\ \ x\in\RR^n ,$$
where
$$F(x)={1\over 2}\sum_{i,j=1}^n a_{ij}x_ix_j ,$$
provided that
$A=(a_{ij})$ and $x=(x_1,...,x_n)$.

The eigenvalues of the operator
 $A$ are, by  Courant's Principle,
$$\lambda_k =\max_{M\in{\cal V}_k}\min_{x\in M}{{\langle Ax,x\rangle}\over{\| x\| ^2}}=
2\max_{A\in{\cal V}_k}\min_{x\in M}F(x)\, ,$$
for $1\leq k\leq n$, where ${\cal V}_k$ denotes the set of all vector subspaces of $\RR^n$
with the dimension
 $k$.

The first result in the
 Lusternik-Schnirelmann theory was proved in 1930
 and is the following:

\medskip
{\bf  Lusternik-Schnirelmann Theorem.} {\sl Let
 $f:\RR^n\rightarrow\RR$ be a $C^1$ even functional.
 Then $f'$ has at least
$2n$ distinct eigenfunctions on the sphere
 $S^{n-1}$}.

\medskip
The ulterior achievements in mathematics showed the signifiance of this
theorem. We only point out that the variational arguments play at this moment a
very strong instrument in the study of potential operators. Hence it is not a
coincidence the detail that the solutions of such a problem are found by analysing
the extrema of a suitable functional.

L. Lusternik and L. Schnirelmann developed their theory using the notion of
{\it Lusternik-Schnirelmann category} of a set. A
 simpler notion is that of {\it
 genus}, which is due to
  Coffman \cite{CO}, but equivalent to that introduced by
Krasnoselski.

Let $X$ be a real Banach space and denote by ${\cal F}$ the family of all closed
and symmetric with respect to the origin subsets of
 $X\setminus\{ 0\}$.

\medskip
\begin{defin}\label{2.1.1}  A nonempty subset $A$ of ${\cal F}$ has the genus $k$
provided $k$ is the least integer with the property that there exists
a continuous odd mapping  $h:A\rightarrow \RR^{k}\setminus\{ 0\}$.

We shall denote from now on by $\gamma (A)$ the genus of the set $A\in{\cal F}$.

By definition, $\gamma (\emptyset )=0$ and
$\gamma (A)=+\infty$,
if $\gamma (A)\not =k$, for every integer $k$.
\end{defin}

\medskip
\begin{lemma}\label{2.1.2} Let $D\subset \RR^n$ be a bounded open and symmetric set
which contains the origin. Let
 $f:\overline{D}\rightarrow\RR^n$ be a continuous function which does not vanish on
 the boundary of
 $D$.

Then $f(D)$ contains a neighbourhood of the origin.
\end{lemma}

\medskip
{\bf Proof.} By Borsuk's Theorem, $d[f;D,0]$ is an odd number, so different from 0.
Now, by the existence theorem for the topological degree, it follows that
$0\in f(D)$. The property of continuity of the topological degree implies the
existence of some $\varepsilon >0$ such that
 $a\in f(D)$, for all $a\in\RR^n$ with
$\| a\| <\varepsilon $. \qed

\medskip
\begin{lemma}\label{2.1.3}  Let $D$ be as in Lemma
\ref{2.1.2} and $g:\partial D\rightarrow\RR^n$
a continuous odd function, such that the set $g(\partial D)$ is contained in a proper
subspace of
 $\RR^n$.

Then there exists $z\in\partial D$ such that $g(z)=0$.
\end{lemma}

\medskip
{\bf Proof.} We may suppose that $g(\partial D)\subset\RR^{n-1}$.
If $g$ does not vanish on $\partial D$, then, by Tietze's Theorem,
there exists an extension $h$ of $g$ at the set $\overline{D}$. By
Lemma \ref{2.1.2}, the set $h(\overline{D})$ contains a
neighbourhood of the origin in $\RR^n$, which is not possible,
because $h(\overline{D})\subset \RR^{n-1}$. Thus, there is some
$z\in\partial D$ such that $g(z)=0$. \qed

\medskip
\begin{lemma}\label{2.1.4}  Let $A\in{\cal F}$ a set which is homeomorphic with $S^{n-1}$
by an odd homeomorphism.

Then $\gamma (A)=n$.
\end{lemma}

\medskip
{\bf Proof.} Obviously, $\gamma (A)\leq n$.

If $\gamma (A)=k<n$, then there exists $h:A\rightarrow\RR^{k}\setminus\{ 0\}$
continuous and odd.

Let $f:A\rightarrow S^{n-1}$ the homeomorphism given in the hypothesis. Then
$h\circ f^{-1}:S^{n-1}\rightarrow\RR^{k}\setminus\{ 0\}$ is continuous and odd, which
contradicts Lemma
 \ref{2.1.3}. Therefore $\gamma (A)=n$.
\qed

\medskip

The main properties of the notion of genus of a closed and symmetric set a listed in
what follows:

\medskip
\begin{lemma}\label{2.1.5}  Let $A,B\in{\cal F}$.

i) If there exists $f:A\rightarrow B$ continuous and odd then
$\gamma (A)\leq\gamma (B)$.

ii) If $A\subset B$, then $\gamma (A)\leq\gamma (B)$.

iii) If the sets $A$ and $B$ are homeomorphic, then
$\gamma (A)=\gamma (B)$.

iv) $\ \ \gamma (A\cup B)\leq\gamma (A)+\gamma (B)\, .$

v) If $\gamma (B)<+\infty$, then
$$\gamma (A)-\gamma (B)\leq \gamma
(\overline{A\setminus B})\, .$$

vi) If $A$ is compact, then $\gamma (A)<+\infty$.

vii) If $A$ is compact, then there exists $\varepsilon >0$ such that
$$\gamma (V_{\varepsilon}(A))=\gamma (A),$$
where
$$V_{\varepsilon}(A)=\{ x\in X;\
{\rm dist}(x,A)\leq\varepsilon\}\, .$$
\end{lemma}

\medskip
{\bf Proof.} i) If $\gamma (B)=n$, let $h:B\rightarrow\RR^n\setminus\{ 0\}$
continuous and odd. Then the mapping $h\circ f:A\rightarrow\RR^n\setminus\{ 0\}$
is also continuous and odd, that is $\gamma (A)\leq n$.

If $\gamma (B)=+\infty$, the result is trivial.

ii) We choose $f=$ Id in  the preceding proof.

iii) It follows from i), by interventing the sets $A$ and $B$.

iv) Let $\gamma (A)=m,\ \gamma (B)=n$ and $f:A\rightarrow \RR^m\setminus\{ 0\},\
g:B\rightarrow\RR^n\setminus\{ 0\}$ be continuous and odd. By
Tietze's Theorem let $F:X\rightarrow\RR^m$ and $G:X\rightarrow\RR^n$ be continuous
extensions of
 $f$ and $g$. Moreover, let us assume that $F$ and $G$ are odd. If not, we replace the
 function
 $F$ with
$$x\longmapsto{{F(x)-F(-x)}\over 2}\, .$$
Let  $$h=(F,G):A\cup B\rightarrow\RR^{m+n}
\setminus\{ 0\}\, .$$

Clearly, $h$ is continuous and odd, that is $\gamma (A\cup B)\leq
\gamma (A)+\gamma (B)$.

v) follows from ii), iv) and the fact that $A\subset (\overline{A\setminus B})\cup B$.

vi) If $x\not =0$ and $r<\| x\|$, then $B_r(x)\cap B_r(-x)=\emptyset$. So,
$$\gamma (B_r(x)\cup B_r(-x))=1\, .$$

By compactness arguments, we can cover the set $A$ with a finite number of open balls,
that is $\gamma (A)<+\infty$.

vii) Let $\gamma (A)=n$ and $f:A\rightarrow\RR^n\setminus\{ 0\}$ be continuous and odd.
 With the same arguments as  in iv), let $F:X\rightarrow\RR^n$  be a continuous and odd
 extension of
 $f$.

Since $f$ does not vanish on the compact set $A$, there is some
 $\varepsilon >0$ such that
$F$ does not vanish in $V_{\varepsilon}(A)$. Thus
 $\gamma (V_{\varepsilon}(A))\leq n=\gamma (A)$.

The reversed inequality follows from ii). \qed

\medskip
We give in what follows the notion of {\it  Lusternik-Schnirelmann category} of a set.
For further details and proof we refer to
Mawhin-Willem \cite{MW1} and Palais \cite{Pa}.

A topological space $X$ is said to be {\it contractible} provided that the identic map
is homotopic with a constant map, that is, there exist
 $u\in X$ and a continuous function
$F:[0,1]\times X\rightarrow X$ such that, for every $x\in X$,
$$F(0,x)=x\quad\hbox{and}\ \ \ F(1,x)=u\, .$$

A subset $M$ of $X$ is said to be {\it contractible in} $X$
is there exist $u\in X$ and a continuous function $F:[0,1]\times M\rightarrow X$
such that, for every $x\in M$,
$$F(0,x)=x\quad\hbox{and}\ \ F(1,x)=u\, .$$

If $A$ is a subset of $X$, define the {\it category of } $A$ {\it in} $X$,
denoted by $\catx (A)$, as follows:

$\catx (A)=0$, if $A=\emptyset$;

$\catx (A)=n$, if $n$ is the smallest integer such that
$A$ may be covered with $n$ closed sets which are
contractible in $X$;

$\catx (A)=\infty$, if contrary.

\medskip
\begin{lemma}\label{2.1.6}  Let $A$ and $B$ be subsets of $X$.

i) If $A\subset B$, then $\catx (A)\leq\catx (B)$.

ii) $\ \ \catx (A\cup B)\leq\catx (A)+\catx (B)\, .$

iii) Let $h:[0,1]\times A\rightarrow X$ be a continuous mapping
such that $h(0,x)=x$, for every $x\in A$.

If $A$ is closed and $B=h(1,A)$, then
$\catx (A)\leq\catx (B)$.
\end{lemma}

\section{A finite dimensional version of the
Lusternik-Schnirelmann theorem}

\hspace*{6mm}The first version of the celebrated Lusternik-Schnirelmann Theorem,
published
 in \cite{LS2}, was generalized in several ways. We shall prove in what follows a finite
 dimensional variant, by using the notion of
 {\it genus} of a set. Other variants of the
Lusternik-Schnirelmann Theorem may be found in
Krasnoselski \cite{Kr}, Palais \cite{Pa},
Rabinowitz \cite{Ra1}, Struwe \cite{St}.

Let $f,g\in C^1(\RR^n,\RR )$ and let $a>0$ be a fixed real number.

\medskip
\begin{defin}\label{2.2.1}  We say that the functional $f$ has a critical point
with respect to $g$ and $a$ if there exist $x\in\RR^n$ and $\lambda\in\RR$ such that
\neweq{2.1}
\left\{ f'(x)=\lambda g'(x)\atop g(x)=a\, .\right.
\endeq

In this case, $x$ is said to be a critical point of $f$ (with respect to the
mapping
$g$ and the number $a$), while $f(x)$ is called a critical value of $f$.

We say that the real number $c$ is a critical value of
 $f$ if the problem (2.1) admits a solution
$x\in\RR^n$ such that $f(x)=c$.
\end{defin}

\medskip
\begin{lemma}\label{2.2.2}  Let $g:\RR^n\rightarrow\RR$ be an even map which is
 Fr\'echet differentiable and such that

1) {\rm Ker}$\ g=\{ 0\}\, ;$

2) $\ \langle g'(x),x\rangle >0$, for every $x\not =0\, ;$

3) $\displaystyle\ \lim_{\| x\|\rightarrow\infty}g(x)=
\infty\, .$

Then the sets $[g=a]$ and $S^{n-1}$ are homeomorphic.
\end{lemma}

\medskip
{\bf Proof.} Let
$$h:[g=a]\rightarrow S^{n-1},\quad\ h(x)={x\over{\| x\|}}
\, .$$

Evidently, $h$ is well defined and continuous. We shall prove in what follows that
  $h$ is one-to-one and onto.

Let $y\in S^{n-1}$. Consider the mapping
$$f:[0,\infty )\rightarrow\RR ,\quad f(t)=g(ty)\, .$$

Then $f$ is differentiable and, for every $t>0$,
$$f'(t)=\langle g'(ty),y\rangle >0\, .$$

Since $f(0)>0$ and $\displaystyle \lim_{t\rightarrow\infty}f(t)=\infty$,
it follows that there exists a unique $t_0>0$ such that $f(t_0)=a$, that is
$g(t_0y)=a$. Thus, $t_0y\in [g=a]$ and $h(t_0y)=y$. Therefore
$h$ is surjective.

Let now $x,y\in [g=a]$ be such that $h(x)=h(y)$, that is
$${x\over{\| x\|}}={y\over{\| y\|}}\, .  $$

If $x\not= y$, then there is some $t_0>0$, $t\not= 1$ such that $y=t_0x$.

Consider the mapping
$$\psi :[0,\infty )\rightarrow\RR ,\quad \psi (t)=g(tx)\, .$$

It follows that $\psi (1)=\psi (t_0)$. But, for every $t>0$,
$$\psi '(t)=\langle g'(tx),x\rangle >0\, ,$$
which implies that the equality $\psi (1)=\psi (t)$ is not possible provided
$t\not= 1$. Thus, $x=y$, that is $h$ is one-to-one.

The condition 3) from our hypotheses implies the continuity of
$h^{-1}$. Moreover, $h^{-1}$ is odd, because $g$ is even. Thus,
$h$ is the desired homeomorphism. \qed

\medskip
For every $1\leq k\leq n$ define the set
$${\cal V}_k=\{ A;\ \ A\subset [g=a],\ A\ \hbox{compact,
symmetric and}\ \gamma (A)\geq k\}\, .$$

Let $F$ be a vector subspace of $\RR^n$ and of dimension $k$. If
$S_k=F\cap S^{n-1}$, it follows by Lemma \ref{2.2.2} that the set
$A=h^{-1}(S_k)$  lies in ${\cal V}_k$, that is
${\cal V}_k\not= \emptyset$, for every $1\leq k\leq n$. Moreover,
$${\cal V}_n\subset{\cal V}_{n-1}\subset{\cal V}_1\, .$$

\medskip
\begin{teo}\label{2.2.3}  Let $f,g:\RR^n\rightarrow\RR$ be two even functionals of
class
 $C^1$ and $a>0$ fixed. Assume that $g$
satisfies the following assumptions:

1) $\ \ ${\rm Ker} $g=\{ 0\}\, ;$

2) $\ \ \langle g'(x),x\rangle >0$,
for every $x\in\RR^n\setminus\{ 0\}\, ;$

3) $\ \ \displaystyle\lim_{\|
x\|\rightarrow\infty}g(x)=\infty\, .$

Under these assumptions, $f$ admits at least $2n$ critical points with respect to
the application $g$ and the number $a$.
\end{teo}

\medskip
{\bf Proof.} Observe first that the critical points appear in pairs, because of the
evenness of the mappings
 $f$ and
$g$.

\underbar{Step 1.} {\it The characterization of the critical values of } $f$.

Let, for every $1\leq k\leq n$,
$$c_k=\sup_{A\in{\cal V}_k}\min_{x\in A}f(x)\, .$$

We propose to show that $c_k$ are critical values of $f$.
This is not enough for concluding the proof, since it is possible that the numbers
$c_k$ are not distinct.

If $c$ is a real number, let
$$A_c=\{ x\in [g=a];\ f(x)\geq c\}\, .$$

We shall prove that, for every $1\leq k\leq n$,
$$c_k=\sup\{ r\in\RR ;\ \gamma (A_r)\geq k\}\, .$$

Set
$$x_k=\sup\{ r\in\RR ;\ \gamma (A_r)\geq k\}\, .$$

>From $\gamma (A_r)\geq k$ it follows that $\inf \{ f(x);\ x\in A_r\}\leq c_k$,
that is $x_k\leq c_k$.

If $x_k<c_k$, then there exists $A\in{\cal V}_k$ such that
$$c_k>\inf_{x\in A}f(x)=\alpha >x_k\, .$$

Thus, $A\subset A_{\alpha}$ and $k\leq \gamma (A)\leq\gamma (A_{\alpha})$, which
contradicts the definition of $x_k$. Consequently, $c_k=x_k$.

Using now the fact that, for every
$\varepsilon >0$, it follows that
$$\gamma (A_{c_k-\varepsilon})\geq k\, .$$

Let $K_c$ be the set of critical values of $f$ corresponding to the critical value
 $c$. By the Deformation Lemma (Theorem A.4 in \cite{Ra1}),
if $V$ is a neighbourhood of $K_c$, there exist $\varepsilon >0$
and $\eta\in C([0,1]\times\RR^n,\RR^n)$ such that, for every fixed
$t\in [0,1]$, the mapping
$$x\longmapsto \eta (t,x)$$
is odd and
$$\eta (1,A_{c-\varepsilon}\setminus V)\subset
 A_{c+\varepsilon}\, .$$

Now, putting for every $x\in \RR^n$,
$$s(x)=\eta (1,x),$$
we get
\neweq{2.2}
s(A_{c-\varepsilon}\setminus V)\subset A_{c+\varepsilon}
\, .\endeq

In particular, if $K_c=\emptyset$, then
$$s(A_{c-\varepsilon})\subset A_{c+\varepsilon}\, .$$

\underbar{Step 2.} {\it For every} $1\leq k\leq n$, {\it the number }
$ c_k$ {\it is a critical value of } $f$.

Indeed, if not, using the preceding result, there exists
 $\varepsilon >0$ such that
$$s(A_{c_k-\varepsilon})\subset A_{c_k+\varepsilon}
\, .$$

>From $\gamma (A_{c_{\scriptstyle k}-\varepsilon})\geq k$.
By Lemma \ref{2.1.5} ii),
it follows that
$$\gamma (s(A_{c_k-\varepsilon}))\geq k\, .$$

The definition of $c_k$ yields
\neweq{2.4}
c_k\geq\inf_{x\in s(A_{c_k-\varepsilon})}f(x)\, .\endeq

By \eq{2.2} and \eq{2.4} it follows that $c_k\geq c_k+\varepsilon$, contradiction.

\underbar{Step 3.} {\it A multiplicity argument}.

We study in what follows the case of multiple critical values.
  Let us assume that
$$c_{k+1}=...=c_{k+p}=c,\ \ \ p>1\, .$$

In this case we shall prove that $\gamma (K_c)\geq p$.

If, by contradiction,  $\gamma (K_c)\leq p-1$, then,
 by Lemma \ref{2.1.5} vii),
there is some $\varepsilon >0$ such that
$$\gamma (V_{\varepsilon}(K_c))\leq p-1\, .$$

Let $V=$ Int $V_\varepsilon (K_c)$. By \eq{2.2} it follows that
$$s(A_{c-\varepsilon}\setminus V)\subset
 A_{c+\varepsilon}\, .$$

Observe that $$B=\overline{A_{c-\varepsilon}-V_\varepsilon (K_c)}=A_{c-\varepsilon}
\setminus \hbox{Int}\ V_\varepsilon (K_c)\, .$$
But $\gamma (A_{c-\varepsilon})\geq k+p$. Using now Lemma
\ref{2.1.5} v), we have
$$\gamma (B)\geq\gamma (A_{c-\varepsilon})-\gamma
(V_\varepsilon (K_c))\geq k+1\, .$$
By Lemma \ref{2.1.5} i), it follows that
$$\gamma (s(B))\geq k+1\, .$$
The definition of $c=c_{k+1}$ shows that
\neweq{2.5}
\inf_{x\in s(B)}f(x)\leq c\, .\endeq
The inclusion $s(B)\subset A_{c+\varepsilon}$ implies
$$\inf_{x\in s(B)}f(x)\geq c+\varepsilon ,$$
which contradicts  \eq{2.5}. \qed

\section {Critical points of locally Lipschitz $Z$-periodic functionals}

\hspace*{6mm}Let $X$ be a Banach space and let $Z$ be a discrete subgroup of it. Therefore
$$\inf_{z\in Z\setminus\{ 0\}}\| z\| >0\, .$$

\begin{defin}\label{2.3.1}  A function \fxr is said to be $Z$-periodic
provided that $f(x+z)=f(x)$, for every $x\in X$ and
 $z\in Z$.
 \end{defin}

\medskip
If the locally Lipschitz functional \fxr is $Z$-periodic,
then, for every $v\in X$, the mapping $x\longmapsto f^0(x,v)$
is $Z$-periodic and $\partial f$ is $Z$-periodic, that is,
for every $x\in X$ and $z\in Z$,
$$\partial f(x+z)=\partial f(x)\, .$$

Thus the functional $\lambda$ inherits the property of $Z$-periodicity.

If $\pi :X\rightarrow X/Z$ is the canonical surjection and $x$ is a critical point of
 $f$, then the set $\pi^{-1}(\pi (x))$ contains only critical points.
Such a set is said to be   a {\it critical orbit} of $f$.
We also remark that $X/Z$ becomes a complete metric space if it is
endowed with the metric
$$ d(\pi (x),\pi (y))=\inf_{z\in Z}\| x-y-z\|\, .$$

A locally Lipschitz functional which is $Z$-periodic \fxr satisfies
the Palais-Smale \psz condition provided that,
 for every sequence $(x_n)$ in $X$ such that $(f(x_n))$ is bounded and
$\lambda (x_n)\rightarrow 0$, there exists a convergent subsequence of $(\pi (x_n))$.
Equivalently, this means that, up to a subsequence, there exists
 $z_n\in Z$ such that the sequence $(x_n-z_n)$ is convergent. If $c$ is a real number,
  then $f$ satisfies the local condition of type Palais-Smale \pszc
if, for every sequence $(x_n)$ in $X$ such that $f(x_n)\rightarrow c$
and $\lambda (x_n)\rightarrow 0$, there exists a convergent subsequence of
$(\pi (x_n))$.

\medskip
\begin{teo}\label{2.3.2}  Let \fxr be a locally Lipschitz functional which is
 $Z$-periodic and satisfies the assumption
 \eq{1.3}.

If $f$ satisfies the condition \pszc , then $c$ is a critical value of
 $f$, corresponding to a critical point which is not in $\pi^{-1}(\pi (p^*(K^*)))$.

\end{teo}

\medskip
{\bf Proof.} With the same arguments as in the proof of Theorem
\ref{1.2.1} we find a sequence $(x_n)$ in $X$ such that
$$\lim_{n\rightarrow\infty}f(x_n)=c\quad\hbox{and}
\quad\lim_{n\rightarrow\infty}\lambda (x_n)=0\, .$$

The Palais-Smale condition \pszc implies the existence of some
 $x$ such that , up to a subsequence, $\pi (x_n)\longrightarrow\pi (x)$. Passing now to the equivalence
 class
mod $Z$, we may assume that $x_n\longrightarrow x$ in $X$.
Moreover, $x$ is a critical point of $f$, because
$$f(x)=\lim_{n\rightarrow\infty}f(x_n)=c$$
and
$$\lambda (x)\leq\liminf_{n\rightarrow\infty}
\lambda (x_n)=0\, .$$ \qed

\medskip
\begin{lemma}\label{2.3.3}  If $n$ is the dimension of the vector space
spanned by the discrete subgroup  $Z$ of $X$,
then, for every $1\leq i\leq n+1$, the set
$${\cal A}_i=\{ A\subset X;\ A\ \hbox{is compact and}\ \ \catpi\pi (A)\geq i\}$$
is nonempty. Moreover
$${\cal A}_1\supset{\cal A}_2\supset ...
\supset{\cal A}_{n+1}\, .$$

\end{lemma}

\medskip
The proof of this result may be found in Mawhin-Willem \cite{MW1}.

\medskip
\begin{lemma}\label{2.3.4}  For every $1\leq i\leq n+1$, the set ${\cal A}_i$ becomes a
complete metric space if it is endowed with the Hausdorff metric
$$\delta (A,B)=\max\{ \sup_{a\in A}{\rm dist}(a,B),
\sup_{b\in B}{\rm dist}(b,A)\}\, .$$
\end{lemma}

\medskip
The proof of this result may be found  in  Kuratovski \cite{Ku}.

\medskip
\begin{lemma}\label{2.3.5}  If \fxr is continuous, then, for every $1\leq i\leq n+1$,
the mapping $\eta :{\cal A}_i\rightarrow\RR$ defined by
$$\eta (A)=\max_{x\in A}f(x)$$
is lower semicontinuous.
\end{lemma}

\medskip
{\bf Proof.} For any fixed $i$, let $(A_n)$ be a sequence in ${\cal A}_i$ and
$A\in{\cal A}_i$ such that $\delta (A_n,A)\longrightarrow 0$.

For every $x\in A$ there exists a sequence $(x_n)$ in $X$ such that $x_n\in A_n$ and
$x_n\rightarrow x$. Thus,
$$f(x)=\lim_{n\rightarrow \infty}f(x_n)
\leq\liminf_{n\rightarrow\infty}\eta (A_n)\, .$$

Since $x\in A$ is arbitrary, it follows that
$$\eta (A)\leq\liminf_{n\rightarrow\infty}\eta (A_n)\, .$$
\qed

\medskip
In what follows, $f:X\rightarrow\RR$ will be a locally Lipschitz
functional which is $Z$-periodic and satisfies the condition
 \psz . Moreover, we shall assume that $f$ is bounded from below.
Let \crfc be the set of critical points of $f$ having the real number
 $c$ as corresponding critical value. Thus,
$$\crfc =\{ x\in X;\ f(x)=c\ \ \hbox{and}\ \
\lambda (x)=0\}\, .$$

If $n$ is the dimension of the vector space spanned
by the discrete group $Z$, then, for every
$1\leq i\leq n+1$, let
$$c_i=\inf_{A\in{\cal A}_i}\eta (A)\, .$$

It follows by Lemma \ref{2.3.3} and the boundedness from below
of $f$  that
$$-\infty <c_1\leq c_2\leq ...\leq c_{n+1}<+\infty\, .$$

\medskip

\begin{teo}\label{2.3.6}  Under the above hypotheses, the functional $f$
has at least
$n+1$ distinct critical orbits.
\end{teo}

\medskip
{\bf Proof.} It is enough to show that if $1\leq i\leq j\leq n+1$
and $c_i=c_j=c$, then the set \crfc contains at least $j-i+1$ distinct critical orbits.
 Arguing by contradiction, let us assume that there exists $i\leq j$ such that
the set \crfc has $k\leq j-i$ distinct critical orbits, generated by $x_1,...,x_k$.
We first choose an open neighbourhood of  \crfc defined by
$$V_r=\bigcup_{l=1}^k\bigcup_{z\in Z}\ B(x_l+z,r)\, .$$

Moreover we may assume that $r>0$ is chosen so that $\pi$ restricted to the set
$\overline{B}\, (x_l,2r)$ is one-to-one. This contradiction shows that,
for every $1\leq l\leq k$,
$$\catpi \ \pi (\overline{B}(x_l,2r))=1\, .$$

In the above arguments, $V_r=\emptyset$ if $k=0$.

\underbar{Step 1.} We shall prove that there exists $\ 0<\varepsilon <\min \{ {1\over 4},r\}$
such that, for every $x\in [c-\varepsilon\leq f\leq c+\varepsilon ]\setminus V_r$, we have
\neweq{2.6}
\lambda (x)>\sqrt{\varepsilon}\, .\endeq
Indeed, if not, there exists a sequence $(x_m)$ in
$X\setminus V_r$ such that, for every $m\geq 1$,
$$c-{1\over m}\leq f(x_m)\leq c+{1\over m}$$
and
$$\lambda (x_m)\leq{1\over{\sqrt{m}}}\, .$$

Since $f$ satisfies the condition \psz , passing eventually to a
subsequence, $\pi (x_m)\longrightarrow \pi (x)$, for some $x\in
V\setminus V_r$. By the $Z$-periodicity property of $f$ and
$\lambda$, we may assume that
 $x_m\longrightarrow x$. The continuity of $f$ and the lower semicontinuity of
$\lambda$ imply $f(x)=0$ and $\lambda (x)=0$, contradiction, because
$x\in V\setminus V_r$.

\underbar{Step 2.} For $\varepsilon$ found above and taking into account the definition of
 $c_j$, there exists $A\in{\cal A}_j$ such that
$$\max_{x\in A}f(x)<c+{\varepsilon}^2\, .$$
Putting $B=A\setminus V_{2r}$ and applying Lemma
\ref{2.1.6} we find
$$j\leq \catpi \pi (A)\leq\catpi (\pi (B)\cup (\overline{V}_{2r}))\leq$$
$$\leq\catpi \pi (B)+\catpi \pi (\overline{V}_{2r})\leq\catpi \pi (B)+k\leq$$
$$\leq\catpi \pi (B)+j-i\, .$$
Thus
$$\catpi \pi (B)\geq i,$$
that is, $B\in{\cal A}_i$.

\underbar{Step 3.} For $\varepsilon$ and $B$ as above, we apply
Ekeland's Principle to the functional
 $\eta$ defined in Lemma \ref{2.3.5}. Thus, there exists
 $C\in {\cal A}_i$ such that, for every $D\in{\cal A}_i$, $D\not =C$,
$$\eta (C)\leq\eta (B)\leq\eta (A)\leq c+{\varepsilon}^2\, ,$$
$$\delta (B,C)\leq\varepsilon\, ,$$
\neweq{2.7}
\eta (D)>\eta (C)-\varepsilon\delta (C,D)\, .\endeq

Since $B\cap V_{2r}=\emptyset$ and $\delta (B,C)\leq\varepsilon <r$, we have
$C\cap V_r=\emptyset$. In particular, the set $F=[f\geq c-\varepsilon ]$ is contained
 in $[c-\varepsilon\leq f\leq c+\varepsilon ]$ and
$F\cap V_r=\emptyset$. Applying now Lemma \ref{1.2.2} for
 $\varphi =\partial f$ defined on
  $F$, we find a continuous map $v:F\rightarrow X$ such that, for every
$x\in F$ and $x^*\in \partial f(x)$,
$$\| v(x)\|\leq 1$$
and
$$\langle x^*,v(x)\rangle\geq\inf_{x\in F}\lambda (x)-\varepsilon \geq\inf_{x\in C}
\lambda (x)-\varepsilon \geq\sqrt{\varepsilon}-\varepsilon ,$$
the last inequality being justified by the relation
\eq{2.6}.
Thus, for every $x\in F$ and $x^*\in\partial f(x)$,
$$f^0(x,-v(x))=\max_{x^*\in\partial f(x)}\langle x^*,-v(x)\rangle
 =-\min_{x^*\in\partial f(x)}
\langle x^*,v(x)\rangle\leq$$
$$\leq \varepsilon -\sqrt{\varepsilon}<-\varepsilon ,$$
by the choice of $\varepsilon$.

By the upper semicontinuity of $f^0$ and the compactness of $F$
there exists $\delta >0$ such that, for every $x\in F,\ y\in X,\
\| y-x\|\leq\delta$ we have
\neweq{2.8}
f^0(y,-v(x))<-\varepsilon\, .\endeq

Since $C\cap\crfc =\emptyset$ and $C$ is compact and \crfc is closed,
 there exists a continuous extension  $w:X\rightarrow X$ of $v$
such that the restriction of $w$ to \crfc is the identic map
 and, for every $x\in X$,
$\| w(x)\|\leq 1$.

Let $\alpha :X\rightarrow [0,1]$ be a continuous function which is $Z$-periodic and
 such that
$\alpha =1$ on $[f\geq c]$ and $\alpha =0$ on $[f\leq c-\varepsilon ]$.
Let $h:[0,1]\times X\rightarrow X$ be the continuous map defined by
$$h(t,x)=x-t\delta \alpha (x)w(x)$$

If $D=h(1,C)$, it follows by Lemma \ref{2.1.6} that
$$\catpi \pi (D)\geq\catpi \pi(C)\geq i,$$
which shows that $D\in{\cal A}_i$, because $D$ is compact.

\underbar{Step 4.} By Lebourg's Mean Value Theorem we find that, for every
$x\in X$, there is some $\theta\in (0,1)$ such that
$$f(h(1,x))-f(h(0,x))\in\langle\partial f(h(\theta ,x)),-\delta \alpha
(x)w(x)\rangle\, .$$

Thus, there exists $x^*\in\partial f(h(\theta ,x))$ such that
$$f(h(1,x))-f(h(0,x))=\alpha (x)\langle
x^*,-\delta w(x)\rangle\, .$$

It follows now by \eq{2.8} that if $x\in F$, then
$$f(h(1,x))-f(h(0,x))=\delta\alpha (x)\langle x^*,-w(x)
\rangle \leq$$
\neweq{2.9}
\leq \delta\alpha (x)f^0(x-\theta\delta\alpha
(x)w(x),-v(x)\rangle \leq -\varepsilon\delta\alpha (x)
\, .\endeq
Thus, for every $x\in C$,
$$f(h(1,x))\leq f(x)\, .$$

Let $x_0\in C$ be such that $f(h(1,x_0))=\eta (D)$. Therefore
$$c\leq f(h(1,x_0))\leq f(x_0)\, .$$

By the definitions of $\alpha $ and $F$ it follows that
$\alpha (x_0)=1$ and $x_0\in F$. Thus, by \eq{2.9}, we have
$$f(h(1,x_0))-f(x_0)\leq -\varepsilon\delta\, .$$
Hence
\neweq{2.10}\eta (D) +\varepsilon\delta \leq f(x_0)\leq\eta
(C)\, .\endeq
Now, by the definition of $D$,
$$\delta (C,D)\leq\delta\, .$$
Thus
$$\eta (D)+\varepsilon\delta (C,D)\leq \eta (C),$$
that is, by \eq{2.7}, we find $C=D$, which contradicts the
relation \eq{2.10}. \qed

\chapter{Applications to the study of  multivalued elliptic problems}

\section{Multivalued variants of some results of
Brezis-Nirenberg
and Mawhin-Willem}

\hspace*{6mm} Let $\Omega$ be an open bounded set with the
boundary sufficiently smooth in $\RR^N$. Let $g$ be a measurable
function defined on $\Omega\times\RR$ and such that
\neweq{3.1}
|g(x,t)|\leq C(1+|t|^p),\quad a.e.\ \ (x,t)
\in\Omega\times\RR \, ,\endeq
where $C$ is a positive constant and $1\leq p <{{N+2}\over{N-2}}$ (if
$N\geq 3$) and $1\leq p<\infty$ (if $N=1,2$).

Define the functional $\psi :L^{p+1}(\Omega )\rightarrow\RR$ by
$$\psi (u)=\int_{\Omega}\int_0^{u(x)}g(x,t)dtdx\, .$$

We first prove that $\psi$ is a locally Lipschitz map.
Indeed the growth condition \eq{3.1} and  H\"older's Inequality yield
$$|\psi (u)-\psi (v)|\leq C'\ (|\Omega |^{p\over{p+1}}+\max_{w\in U}\|
 w\|^{p\over{p+1}}_{L^{p+1}(\Omega )})\cdot
\| u-v\|_{L^{p+1}(\Omega )},$$
where $U$ is an open ball containing $u$ and $v$.

Put
$$\ginf (x,t)=\limeps\, essinf\{ g(x,s);\
|t-s|<\varepsilon\}\, ,$$
$$\gsup (x,t)=\limeps\,
 esssup\{ g(x,s);\ |t-s|<\varepsilon\}\, .$$

\medskip
\begin{lemma}\label{3.1.1} \sl The mappings \ginf and
\gsup are measurable.
\end{lemma}

\medskip\rm
{\bf Proof.} Observe that
$$\gsup (x,t)=\limeps esssup\{ g(x,s);\ s\in
 [t-\varepsilon ,t+\varepsilon ]\}=$$
$$=\limn esssup\{ g(x,s);\ s\in [\tmn ,\tpn ]\}\, .$$

Replacing, locally, the map $g$ by $g+M$, for $M$ large enough, we may assume that
$g\geq 0$. It follows that
$$\gsup (x,t)=\limn \| g(x,\cdot )\|_{L^{\infty}([\tmn ,\tpn ])}=$$
$$=\limn \lim_{m\rightarrow\infty} \| g(x,\cdot )\|_{L^{m}([\tmn ,\tpn])}=$$
$$=\limn\lim_{m\rightarrow\infty} g_{m,n}(x,t),$$
where
$$g_{m,n}(x,t)=(\int_{\tmn}^{\tpn}(g(x,s))^mds)^{{1\over m}}$$
Thus, it is sufficient to prove that if
 $h\in L^{\infty}_{loc}(\Omega\times\RR )$ and
$a>0$, then the mapping
$$k:\Omega \times\RR\rightarrow\RR,\quad\ k(x,t)=\int_{t-a}^{t+a}h(x,s)ds$$
is measurable and to apply this result for
$h(x,s)=(g(x,s))^m$ and $a={1\over n}$.

Observe that, assuming we have already proved that
$$\Omega\times\RR\times\RR\ni (x,s,t) \buildrel l \over{\longmapsto} h(x,s+t)$$
is measurable, then the conclusion follows by
$$k(x,t)=\int_{-a}^a l(x,s,t)ds$$
and  Fubini's Theorem applied to the locally integrable function
 $l$ (the local integrability follows from its local boundedness).
So, it is sufficient to justify the measurability of the mapping
 $l$. In order to do this, it is enough to prove that reciprocal images of Borel
 sets (resp. negligible sets) through the
 function
$$\Omega\times\RR\times\RR\ni (x,t,s)\buildrel \omega \over{\longmapsto}
(x,s+t)\in\Omega\times\RR$$
are Borel sets (resp., negligible). The first condition is, obviously, fulfilled.
For the second, let  $A$ be a measurable set of null measure
 in $\Omega\times\RR$. Consider the map
$$\RR^N\times\RR\times\RR\ni (x,t,s)\buildrel \eta
\over{\longmapsto}(x,s+t,t)\in \RR^N\times\RR\times\RR\, .$$

We observe that $\eta$ is invertible and of class $C^1$. Moreover
$$\eta (\omega^{-1}(A))\subset A\times\RR ,$$
which is negligible. So, $\omega^{-1}(A)$ which is negligible,
too. \qed

\medskip
Let $G:\Omega\times\RR\rightarrow\RR$ be defined by
$$G(x,t)=\int_0^t g(x,s)ds\, .$$

\medskip
\begin{lemma}\label{3.1.2} \sl Let $g$ be a locally bounded measurable function,
 defined on
$\Omega\times\RR$ and $\ginf ,\gsup$ as above.

Then the Clarke subdifferential of $G$ with respect to $t$ is given by
$$\partial_{t}G(x,t)=[\ginf (x,t),\gsup (x,t)],\quad a.e.\ \ (x,t)
\in\Omega\times\RR\, .$$

The condition ``a.e." can be removed if, for every $x$, the mapping
$t\longmapsto g(x,t)$ is measurable and locally bounded.
\end{lemma}

\rm\medskip
{\bf Proof.} We show that we have equality on the set
$$\{ (x,t)\in\Omega\times\RR ;\ g(x,\cdot )\ \ \hbox{is locally bowhered and
measurable}
 \}\, .$$

In order to prove our result, it is enough to consider mappings
$g$ which do not depend on $x$. For this aim the equality that we have to prove is
equivalent to
\neweq{3.2}
G^{0}(t;1)=\gsup (t)\quad\hbox{and}\ \
G^0(t;-1)=\ginf (t)\, .\endeq

Examining the definitions of $G^0,\gsup$ and $\ginf$, it follows that
$$\ginf (t)=-(\overline{-g})(t)\quad\hbox{and}\ \
\ G^0(t,-1)=-(-G)^0(t,1)\, .$$

So, the second equality appearing in \eq{3.2} is equivalent to the first one.

The inequality $G^0(t,1)\leq\gsup (t)$ is proved in Chang
\cite{Ch}. For the reversed
inequality, we assume by contradiction that there exists
 $\varepsilon >0$ such that
$$G^0(t,1)=\gsup (t)-\varepsilon\, .$$

Let $\delta >0$ be such that
$${{G(\tau +\lambda)-G(\tau )}\over{\lambda}}<\gsup
(t)-{\varepsilon\over 2} \, ,$$
if $0<|\tau -t|<\delta$ and $0<\lambda <\delta$. Then
\neweq{3.3}
{1\over\lambda}\int_{\tau}^{\tau +\lambda} g(s)ds<\gsup
(t)-
{\varepsilon\over 2},\quad \hbox{if}\ \ |\tau -t
 |<\delta,\ \lambda >0\, .\endeq

We now justify the existence of some $\lambda_n\searrow 0$ such that
\neweq{3.4}
{1\over{\lambda_n}}\int_{\tau}^{\tau +\lambda_n}g(s)ds\longrightarrow
g(\tau ),\quad a.e.\ \ \tau\in (t-\delta , t+\delta )
\, .\endeq

Assume, for the moment, that \eq{3.4} has already been
proved.
 Then by \eq{3.3} and \eq{3.4} it follows that
for every $\tau\in (t-\delta ,t+\delta )$,
$$g(\tau )\leq \gsup (t)-{\varepsilon\over 2}\, .$$
Thus we get the contradiction
$$\gsup (t)\leq\ esssup\, \{g(s);\ s\in [t-\delta
,t+\delta ]\}\leq\gsup (t)-{\varepsilon\over 2}\, .$$

For concluding the proof it is sufficient to prove  \eq{3.4}. Observe that we can
 ``cut off" the mapping $g$, in order to have
  $g\in L^{\infty}\cap L^1$. Then \eq{3.4} is nothing else that the classical result
\neweq{3.5}
T_{\lambda}\longrightarrow \hbox{Id}_{L^{1}(\RR )},
\quad\hbox{if}\ \
\lambda\searrow 0\, .\endeq
in ${\cal L}(L^1(\RR ))$, where
$$T_{\lambda }\ u(t)={1\over\lambda}\int_t^{t+\lambda}
\ u(s)ds\, ,$$
for $\lambda >0,\ t\in\RR ,\ u\in L^1(\RR )$.

Indeed, we observe easily that $T_{\lambda}$ is linear and continuous in
in $L^1(\RR )$ and that
$$\lim_{\lambda\searrow 0}T_{\lambda }\ u=u\quad\ \hbox{in}\ \ {\cal D}(\RR ),$$
for $u\in{\cal D}(\RR )$. The relation \eq{3.5} follows now
 by a density argument.
\qed

\medskip
Returning to our problem, it follows by Theorem 2.1 in
Chang \cite{Ch} that
\neweq{3.6}
\partial\psi_{| H^1_0(\Omega )}\ (u)\subset
 \partial\psi (u)\, .\endeq

For obtaining further information concerning $\partial\psi$, we need the following
refinement of Theorem 2.1 in \cite{Ch}.

\medskip
\begin{teo}\label{3.1.3} \sl If $u\in\lpu$ then
$$\partial\psi (u)(x)\subset [\ginf (x,u(x)),\gsup (x,u(x))],\quad\ a.e.\
\ x\in\Omega ,$$
in the sense that if $w\in\partial \psi(u)$, then
\neweq{3.7}
\ginf (x,u(x))\leq w(x)\leq\gsup (x,u(x)),\quad\ a.e.\
 \ x\in\Omega\, .\endeq
 \rm\end{teo}

\medskip
{\bf Proof.} Let $h$ be a Borel function such that $h=g$ a.e. in
$\Omega\times\RR$. Then
$$A=\Omega\setminus\{ x\in\Omega ;\ h(x,t)=g(x,t)\ \ \ a.e.\ \ t\in\RR\}$$
is a negligible set. Thus
$$B=\{ x\in\Omega ;\ \ \hbox{there exists}\ \ t\in\RR\ \ \hbox{such that} \ \
\ginf (x,t)\not =\underline{h} (x,t)\}$$
is a negligible set. A similar reasoning may be done for
$\gsup$ and $\overline{h}$.

It follows that we can assume $g$ is a Borel function.

\medskip
\begin{lemma}\label{3.1.4} \sl Let $g:\Omega\times\RR
\rightarrow\RR$ be a locally bounded
Borel function. Then \gsup is a Borel function.\rm
\end{lemma}

\medskip
{\bf Proof of Lemma.} Since the restriction is local, we may assume that
$g$ is nonnegative and bounded by 1. Since
$$g=\limn\ \lim_{m\rightarrow\infty} g_{m,n}\, ,$$
where
$$g_{m,n}(x,t)=(\int_{\tmn}^{\tpn}\ |g^m(x,s)|ds)^{{1\over m}},$$
it is enough to show that $g_{m,n}$ is a Borel function.

Set
$${\cal M}=\{ g:\Omega\times\RR\rightarrow\RR ;\ |g|\leq 1\ \hbox{and}\ \
g\ \ \hbox{is borelian}\}\, ,$$
$${\cal N}=\{ g\in{\cal M};\ g_{m,n}\ \
\hbox{is borelian}\}\, .$$

Evidently, ${\cal N}\subset{\cal M}$. By a classical result from
Measure Theory (see Berberian \cite{Be}, p.178) and the Lebesgue
Dominated Convergence Theorem, we find ${\cal M}\subset{\cal N}$.
Consequently, ${\cal M}= {\cal N}$. \qed

\medskip
{\bf Proof of Theorem \ref{3.1.3} continued.} Let $v\in L^{\infty}(\Omega )$.
There exist the sequences $\lambda_{i}\searrow 0$ and $h_i\rightarrow 0$ in
\lpu such that
$$\psi^0(u,v)=\lim_{i\rightarrow\infty}\ {1\over{\lambda_{i}}}
\int_{\Omega}\ \ \int_{u(x)+h_i(x)}^{u(x)+
h_i(x)+\lambda_iv(x)}\ \ g(x,s)dsdx\, .$$

We may assume that $h_i\rightarrow 0$ a.e. So
$$\psi^0(u,v)=\lim_{i\rightarrow\infty}\ {1\over{\lambda_{i}}}
\int_{[v>0]}\ \ \int_{u(x)+h_i(x)}^{u(x)+h_i(x)+
\lambda_iv(x)}\ \ g(x,s)dsdx\leq$$
$$\leq\int_{[v>0]}(\limsup_{i\rightarrow\infty}
{1\over{\lambda_i}}\int_{u(x)+h_i(x)}^{u(x)+h_i(x)+\lambda_iv(x)}\ \
g(x,s)ds)dx\leq$$
$$\leq\int_{[v>0]}\gsup (x,u(x))v(x)dx\, .$$
So, for every $v\in L^{\infty}(\Omega )$,
\neweq{3.8}
\psi^0(u,v)\leq\int_{[v>0]}\ \gsup (x,u(x))v(x)dx\endeq

Let us now assume that \eq{3.7} is not true. So, there exist $\varepsilon >0$,
 a set $E$ with $|E|>0$ and $w\in\partial\psi (u)$ such that, for every $x\in E$,
\neweq{3.9}
w(x)\geq\gsup (x,u(x))+\varepsilon\, .\endeq
Putting $v=\chi _{E}$ in \eq{3.8} it follows that
$$\langle w,v\rangle =\int_{E}w\leq \psi ^0(u,v)\leq\int_{E}\gsup (x,u(x))dx,$$
which contradicts \eq{3.9}. \qed

\medskip
We assume in what follows that
\neweq{3.10}
g(x,0)=0\ \ \hbox{and}\ \ \limeps \, esssup\left\{
|{{g(x,t)}\over t}|;\ (x,t)\in\Omega\times
 [-\varepsilon ,\varepsilon ]\right\}<\lambda_1\, ,
 \endeq
where $\lambda_1$ is the first eigenvalue of the operator
 $(-\Delta )$ in \huo .
Furthermore, we shall assume that the following
``tehnical" condition is fulfilled:

 there exist $\mu >2$ and $r\geq 0$ such that
\neweq{3.11}
\mu G(x,t)\leq\left\{ t\ginf (x,t),\ a.e.\
x\in\Omega ,t\geq r\atop
t\gsup (x,t),\ a.e.\ x\in \Omega ,t\leq -r\right.\ \ \hbox{and}\
g(x,t)\geq 1\ a.e.\ x\in\Omega ,t\geq r\, .\endeq

It follows from \eq{3.1} and \eq{3.10} that there exist some constants
 $0<C_1<\lambda_1$
and $C_2>0$ such that
\neweq{3.12}
|g(x,t)|\leq C_1|t|+C_2|t|^p,\ \ a.e.\ (x,t)
\in\Omega\times\RR \, .\endeq

\medskip
\begin{teo}\label{3.1.5} \sl Let $\alpha ={{\lambda_1}\over{\lambda_1-C_1}}>1$
and
$a\in L^{\infty}(\Omega )$ be such that the operator
 $-\Delta +\alpha a$ is coercive. Assume, further, that
the conditions \eq{3.1}, \eq{3.10} and \eq{3.11}
 are fulfilled.

Then the multivalued elliptic problem
\neweq{3.13}
-\Delta u(x)+a(x)u(x)\in [\ginf (x,u(x)),\gsup (x,u(x))],\ \ a.e.
 \ x\in\Omega\endeq
has a solution in $\huo\cap W^{2,q}(\Omega )\setminus\{ 0\}$, where
$q$ is the conjugated exponent of $p+1$.\rm\end{teo}

\medskip
{\bf Remark}. The technical condition imposed to $a$ is automatically
fulfilled
if $a\geq 0$ or, more general, if $\| a^{-}\|_{L^{\infty}(\Omega )}
<\lambda_1\alpha^{-1}$.

\medskip
{\bf Proof.} Consider in the space \huo the locally Lipschitz functional
$$\varphi (u)={1\over 2}{}\|\nabla u\|^2_{L^2(\Omega )}+{1\over 2}
\int_{\Omega}
a(x)u^2(x)dx-\psi (u)\, .$$
In order to prove the theorem, it is enough to show that
  $\varphi$ has a critical point $u_0\in\huo$ corresponding to a positive
  critical value. Indeed, it is obvious that
$$ \partial\varphi (u)=-\Delta u+a(x)u-\partial\psi _{|\huo}(u),
\quad\hbox{in}\ \ H^{-1}(\Omega )\, .$$

If $u_0$ is a critical point of $\varphi$, it follows that there exists
 $w\in\partial\psi_{|\huo}(u_0)$ such that
$$-\Delta u_0+a(x)u_0=w,\quad\hbox{in}\ \ H^{-1}(\Omega )$$
But $w\in L^q(\Omega )$. By a standard  regularity result for elliptic
equations,
we find that $u_0\in W^{2,q}(\Omega )$ and that
 $u_0$ is a solution of \eq{3.13}.

In order to prove that $\varphi$ has such a critical point we shall apply
Corollary
 \ref{1.2.4}. More precisely, we shall prove that $\varphi$ satisfies the
Palais-Smale condition, as well as the following ``geometric" hypotheses:

\neweq{3.14}
\varphi (0)=0\ \ \hbox{and there exists} \ v\in\huo\setminus\{ 0\}
 \ \hbox{such that} \ \varphi (v)\leq 0\, ;\endeq
\neweq{3.15}
\hbox{there exist}\ \ c>0\ \hbox{and}\ 0<R<\| v\|\ \hbox{such that}
  \ \varphi \geq c\ \hbox{pe}\ \partial B(0,R)\, .\endeq

\underbar{Verification of \eq{3.14}.} Evidently, $\varphi (0)=0$.
For the other assertion appearing in \eq{3.15}
we need

\medskip
\begin{lemma}\label{3.1.6} \sl There exist positive constants
 $C_1$ and $C_2$ such that
\neweq{3.16}
g(x,t)\geq C_1t^{\mu -1}-C_2,\ \ a.e.\ \ (x,t)
\in\Omega\times [r,\infty )\, .\endeq\rm
\end{lemma}

\medskip
{\bf Proof of Lemma.} We shall show that
\neweq{3.17}
\ginf\leq g\leq\gsup,\ \ \ \hbox{ a.e. in}\ \ \Omega\times
[r,\infty )\ . \endeq

Assume for the moment that the relation \eq{3.17} was proved.
Then, by \eq{3.11},
\neweq{3.18}
\mu \underline{G}(x,t)\leq t\ginf (x,t),\ \ a.e.\ \ (x,t)
\in\Omega\times\RR\, ,\endeq
where
$$\underline{G}(x,t)=\int_0^t \ginf (x,s)ds\, .$$

Consequently it is sufficient to prove \eq{3.16} for $\ginf$ instead of $g$.

Since $\ginf (x,\cdot )\in L^{\infty}_{loc}(\RR )$, it follows that
$\underline{G} (x,\cdot )\in W^{1,\infty}_{loc}(\RR )$. if we choose
 $C>0$ sufficiently large so that
\neweq{3.19}
\mu\underline{G}(x,t)\leq C+t\ginf (x,t),\ \ a.e.\ \ (x,t)
\in\Omega\times [0,\infty )\, .\endeq
We observe that \eq{3.19} shows that
$$(0,+\infty )\ni t\longmapsto {{\underline{G}}(x,t)
\over{t^\mu}}-{C\over{\mu t^{\mu}}}$$
is increasing. So, there exist $R$ large enough and positive constants
  $K_1,K_2$
such that
$$\underline{G}(x,t)\geq K_1t^\mu -K^2,\ \ a.e.\ \
(x,t)\in\Omega\times [R,+\infty )\, .$$
Relation \eq{3.16} follows now from the above inequality
 and from \eq{3.18}.

For proving the second inequality appearing in \eq{3.17}, we observe that
\neweq{3.20}
\gsup =\limn g_n\quad\hbox{in}\ \ L^{\infty}_{loc}(\Omega )
\, ,\endeq
where
$$g_n(x,t)=esssup\ \{g(x,s);\ \ |t-s|\leq{1\over n}\}\, .$$

For fixed $x\in\Omega$, it is sufficient to show that for every
interval $I=[a,b]\subset\RR$ we have
$$\int_{I}\ \gsup (x,t)dt\geq\int_{I}\ g(x,t) dt$$
and to use then a standard argument from Measure Theory.

Taking into account \eq{3.20}, it is enough to show that
$$\liminf_{n\rightarrow\infty}\int_{I}g_n(x,t)dt\geq
\int_{I}g(x,t)dt\, .$$
We have
$$\int_{I}g_n(x,t)dt=\int_{I}esssup\ \{g(x,s);\ s\in [\tmn ,\tpn ]\}\geq$$
$$\geq\int_{I}{n\over 2}{}\int_{\tmn}^{\tpn}g(x,s)ds=\ ({\rm Fubini})
\ \int_{a-{1\over n}}^{b+{1\over n}}{n\over 2}
\int_{t_1(s)}^{t_2(s)}\ g(x,s)dtds=$$
$$=\int_{a-{1\over n}}^{b+{1\over n}}\ {n\over 2}
\ (t_2(s)-t_1(s))\ g(x,s)ds=$$
$$=\int_a^bg(x,s)ds+{n\over 2}\int_{a-{1\over n}}^{b+{1\over n}}
(s-{1\over n}-a)g(x,s)ds+$$
$$+{n\over 2}\int_{b-{1\over n}}^{b+{1\over n}}(b-s-{1\over n})
g(x,s)ds\longrightarrow \int_a^bg(x,s)ds\, .$$

We have chosen $n$ such that ${2\over n}\leq b-a$, and
$$t_1(s)=\itab
$a,$ & \quad $\hbox{if}\ \ a-{1\over n}\leq s\leq a+{1\over n}$\\
$s-{1\over n},$ & \quad $\hbox{if}\ \ a+{1\over n}\leq s\leq b+{1\over n}$\\
\ttab $$
$$t_2(s)=\itab
$s+{1\over n},$ & \quad $\hbox{if}\ \ a-{1\over n}\leq s\leq b-{1\over n}$\\
$b,$ & \quad $\hbox{if}\ \ b-{1\over n}\leq s\leq b+
{1\over n}\, .$\\
\ttab $$ \qed

\medskip
{\bf Proof of Theorem \ref{3.1.5} continued.} If $e_1>0$ denotes the first
eigenfunction of the operator
 $-\Delta$ in \huo , then, for $t$
large enough,
$$\varphi (te_1)\leq{{\lambda_1t^2}\over 2}\| e_1\|^2_{L^2(\Omega )}+
{{t^2}\over 2}\intom ae_1^2-\psi (te_1)\leq$$
$$\leq ({{\lambda_1}\over 2}\| e_1\|^2_{\ldo}+\intom ae_1^2)
t^2+C_2t\intom e_1-C_1't^{\mu}\intom e_1^{\mu}<0$$

So, in order to obtain \eq{3.14} it is enough to choose
 $v=te_1$, for $t$ found above.

\underbar{Verification of \eq{3.15}.} Applying Poincar\'e's Inequality,
the
Sobolev embedding theorem and  \eq{3.12}, we find that, for every $u\in\huo$,
$$\psi (u)\leq{{C_1}\over 2}\intom u^2+{{C_2}\over{p+1}}
\intom |u|^{p+1}\leq
{{C_1}\over{2\lambda_1}}\| \nabla u\|^2_{\ldo}+C'
\|\nabla u\|^{p+1}_{\ldo}\, .$$

Hence
$$\varphi (u)\geq{1\over 2}(1-{{C_1}\over{\lambda_1}})\| \nabla u
\|^2_{\ldo}+{1\over 2}\intom au^2-C'\|\nabla u\|^{p+1}_{\ldo}\geq$$
$$\geq{1\over 2}(1-{{C_1}\over{\lambda_1}}-{1\over{\alpha +\varepsilon}})\|\nabla u\|^2_{\ldo}-C'\|\nabla u\|^{p+1}_{\ldo}\geq C>0,$$
for $\varepsilon >0$ sufficiently small, if
$\|\nabla u\|_{\ldo}=R$ is close to 0.

\underbar{Verification of the Palais-Smale condition.} Let
 $(u_k)$ be a sequence in
\huo such that
$$\varphi (u_k)\ \ \ \hbox{is bounded}$$
and
$$\lim_{k\rightarrow\infty}\lambda (u_k)=0\, .$$

The definition of $\lambda$ and  \eq{3.6} imply the existence of a sequence
 $(w_k)$ such that
$$w_k\in\partial \psi _{|\huo}\subset L^q(\Omega )$$
and
$$-\Delta u_k+a(x)u_k-w_k\longrightarrow 0\ \ \hbox{in}\
\ H^{-1}(\Omega )\, .$$

Since, by \eq{3.1}, the mapping $G$ is locally bounded with respect to the
variable $t$ and uniformly with respect to
 $x$, the hypothesis \eq{3.11} yields
$$\mu G(x,u(x))\leq\left\{u(x)\ginf (x,u(x))+C,\ \ \ a.e.\ \hbox{in}\
[u\geq 0]\atop
u(x)\gsup (x,u(x))+C,\ \ \ a.e.\ \hbox{in}\
 [u\leq 0]\, \right.$$
where $u$ is a measurable function defined on $\Omega$, while
 $C$ is a positive constant not depending on
 $u$. It follows that, for every $u\in\huo$,
$$\psi (u)=\int_{\upo}G(x,u(x))dx+\int_{\umo}G(x,u(x))dx\leq$$
$$\leq{1\over\mu}\int_{\upo}u(x)\ginf (x,u(x))dx+{1\over\mu}
\int_{\umo}u(x)\gsup (x,u(x))dx+C|\Omega |$$

This inequality  and \eq{3.7} show that, for every $u\in\huo$ and
$w\in\partial\psi (u)$,
$$\psi (u)\leq{1\over\mu}\intom u(x)w(x)dx+C'\, .$$

We first prove that the sequence $(u_k)$ contains a subsequence which is
weakly convergent
in \huo . Indeed,
$$\varphi (u_k)={1\over 2}\intom |\nabla u_k|^2+{1\over 2}
\intom au_k^2-\psi (u_k)=$$
$$=({1\over 2}-{1\over\mu})\intom (|\nabla u_k|^2+au_k^2)+
{1\over\mu}\langle
-\Delta u_k+au_k-w_k,u_k\rangle +$$
$$+{1\over\mu}\langle w_k,u_k\rangle -\psi (u_k)\geq$$
$$\geq ({1\over 2}-{1\over\mu})\intom (|\nabla u_k|^2+au_k^2)+
{1\over\mu}\langle -\Delta u_k+au_k-w_k,u_k\rangle -C'\geq$$
$$\geq C''\intom |\nabla u_k|^2-{1\over\mu}\, o(1)\,
\sqrt{\intom |\nabla u_k|^2}-C'\, .$$

This implies easily that the sequence $(u_k)$ is bounded in \huo . So, up to
a subsequence,
 $(u_k)$ is weakly convergent to $u\in\huo$.
Since the embedding $\huo\subset L^{p+1}(\Omega )$ is compact,
it follows that, up to a subsequence,
 $(u_k)$ is strongly convergent in
$L^{p+1}(\Omega )$. We remark that $(u_k)$ is bounded in $L^q(\Omega )$.
Since $\psi$ is Lipschitz on the bounded subsets  of $L^{p+1}(\Omega )$,
it follows that $(w_k)$ is bounded in $L^q(\Omega )$. From
$$\|\nabla u_k\|^2_{\ldo}=\intom \nabla u_k\nabla u-\intom au_k(u_k-u)+$$
$$+\intom w_k(u_k-u)+\langle -\Delta u_k+au_k-w_k,u_k-u\rangle_{H^{-1},H_0^1(\Omega )}$$
it follows that
$$\|\nabla u_k\|_{\ldo}\longrightarrow
\|\nabla u\|_{\ldo}\, .$$

Since $\huo$ is a Hilbert space, and
$$u_k\rightharpoonup u,\ \|u_k\|_{\huo}
\rightarrow\|u\|_{\huo}\, ,$$ we deduce that $(u_k)$ converges to
$u$ in $\huo$. \qed

\medskip
A stronger variant of Theorem \ref{3.1.5} is

\medskip
\begin{teo}\label{3.1.7} \sl Under the same hypotheses as
 in Theorem \ref{3.1.5}, let
$b\in L^q(\Omega )$ such that there exists $\delta >0$ so that
$$\| b\|_{L^\infty(\Omega )}<\delta\, .$$

Then the problem
\neweq{3.24}
-\Delta u(x)+a(x)u(x)+b(x)\in [\ginf (x,u(x)),\gsup (x,u(x))],\
\ a.e.\ \ x\in\Omega\endeq
has a solution.\rm
\end{teo}

\medskip
{\bf Proof.} Define
$$\varphi (u)={1\over 2}\| \nabla u\|^2_{\ldo}+{1\over 2}
\intom au^2+\intom bu
-\intom\int_0^{u(x)}g(x,t)dtdx\, .$$

We have seen that if $b=0$, then the problem \eq{3.24} has a solution.
For $\| b\|_{L^\infty (\Omega )}$ sufficiently small, the verification of the
Palais-Smale condition and of the geometric hypotheses
 \eq{3.14} and \eq{3.15} follows the same lines as
 in the proof of Theorem \ref{3.1.5}.
\qed

\medskip
As a second application of the abstract theorems proved in the first two
chapters, we shall study the multivalued pendulum problem
\neweq{3.25}
\left\{x''+f\in [\ginf (x),\gsup (x)]\atop x(0)=
x(1)\, ,\right. \endeq
where
\neweq{3.26}
f\in L^p(0,1),\quad\ \hbox{for some}\ \ p>1\, ,
\endeq
\neweq{3.27}
g\in L^\infty (\RR )\quad\ \hbox{and there exists}\ \ T>0\ \ \hbox{such that} \
g(x+T)=g(x),\ a.e.\ x\in\RR\, ,\endeq
\neweq{3.28}
\int_0^T g(t)dt=\int_0^1 f(t)dt=0\, .\endeq

The smooth variant of the problem
 \eq{3.25} was studied in  Mawhin-Willem \cite{MW2}.

\medskip
\begin{teo}\label{3.1.8} \sl If $f$ and $g$ are as above then the problem
\eq{3.25} has at least two  solutions in the space
$$X:=H^1_p(0,1)=\{ x\in H^1(0,1);\ x(0)=x(1)\}\, .$$
Moreover, these solutions are distinct, in the sense that their difference is
not an integer multiple of
 $T$.\rm\end{teo}

\medskip
{\bf Proof.} As in the proof of Theorem \ref{3.1.5}, the critical points of the
functional $\varphi :X\rightarrow\RR$ defined by
$$\varphi (x)=-{1\over 2}\int_0^1x'^2+\int_0^1fx-\int_0^1G(x)$$
are solutions of the problem \eq{3.25}.
The details of the proof are, essentially, the same as above.

Since $\varphi (x+T)=\varphi (x)$, we may apply Theorem
 \ref{2.3.6}. All we have to do is to prove that $\varphi$
verifies the condition $\pszc$, for any $c$.

In order to do this, let $(x_n)$ be a sequence  in $X$ such that
\neweq{3.29}
\limn\varphi (x_n)=c\endeq
\neweq{3.30}
\limn\lambda (x_n)=0\, .\endeq

Let
\neweq{3.31}
w_n\in\partial\varphi (x_n)\subset L^\infty (0,1)
\endeq
be such that
$$\lambda (x_n)=x_n''+f-w_n\rightarrow 0\quad\hbox{in}\
\ H^{-1}(0,1)\, .$$

Observe that the last inclusion appearing in \eq{3.31}
 is justified by the fact
that
$$\ginf\circ x_n\leq w_n\leq\gsup\circ x_n,$$
and $\ginf ,\gsup\in L^{\infty}(\RR )$.

By \eq{3.30}, after multiplication with $x_n$ it follows that
$$\int_0^1 (x_n')^2-\int_0^1fx_n+\int_0^1w_nx_n=o(1)\,
\| x_n\|_{H^1_p}\, .$$

Then, by \eq{3.29},
$$-{1\over 2}\int_0^1(x_n')^2+\int_0^1fx_n-
\int_0^1G(x_n)\longrightarrow c\, .$$

So, there exist positive constants $C_1$ and $C_2$ such that
$$\int_0^1(x_n')^2\leq C_1+C_2\| x_n\|_{H^1_p}\, .$$

We observe that $G$ is also $T$-periodic, so bounded.

For every $n$, replacing $x_n$ with $x_n+k_nT$ for a convenable integer $k_n$,
we can assume that
$$x_n(0)\in [0,T]\, .$$

We have obtained that the sequence $(x_n)$ is bounded in $H^1_p(0,1)$.

Let $x\in H^1_p(0,1)$ be such that , up to a subsequence,
$$x_n\rightharpoonup x\quad\ \hbox{and}\quad\
x_n(0)\rightarrow x(0)\, .$$
Thus
$$\int_0^1(x_n')^2=\langle -x_n''-f+w_n,x_n-x\rangle +$$
$$+\int_0^1w_n(x_n-x)-\int_0^1f(x_n-x)+\int_0^1x_n'x'\rightarrow
\int_0^1x'^2,$$
since $x_n\rightarrow x$ in $L^{p'}$, where $p'$ is the conjugated exponent
of $p$.

It follows that $x_n\rightarrow x$ in $H^1_p$, so \pszc is
fulfilled. \qed

\section{Multivalued problems of
Landesman-Lazer type  with strong resonance at infinity}

\hspace*{6mm}In \cite{LL} Landesman and Lazer studied for the first time
problems with resonance and they found sufficient conditions for the
existence of a solution. We shall first recall the main definitions, in the
framework of the singlevalued problems.

Let $\Omega$ be an open bounded open set in $\RR^N$, and
$f:\RR\rightarrow\RR$  a $C^1$ map. We consider the problem
\neweq{3.312}
\left\{ -\Delta u=f(u),\quad\hbox{in}\ \ \Omega\atop
u\in\huo\, .\quad\right.\endeq

For obtaining information on the existence of solutions, as well
as estimates on the number of solutions, it is necessary
 to have further information $f$. In fact the solutions of
  \eq{3.312} depend in an essential manner on the asymptotic behaviour of
 $f$. Let us assume, for example, that
$f$ is asymptotic linear, that is
 $\displaystyle{{f(t)}\over t}$ has a finite limit as
 $|t|\rightarrow\infty$. Let
\neweq{3.32}
a=\lim_{|t|\rightarrow\infty}{{f(t)}\over t}\, .\endeq
We write
$$f(t)=at-g(t),$$
where
$$\lim_{|t|\rightarrow\infty}{{g(t)}\over t}=0\, .$$

Let $0<\lambda_1<\lambda_2\leq ...\leq\lambda_n\leq ...$ the eigenvalues of
the operator $(-\Delta )$ in \huo . The problem \eq{3.312}
is said to be with
{\it resonance at infinity} if the number $a$ from \eq{3.32} is one of the
eigenvalues of $\lambda_k$. With respect to the growth of $g$ at infinity
there are several ``degrees" of resonance. If $g$ has a  ``smaller" rate of
increase at infinity, then its resonance is
 ``stronger".

There are several situations:

a)$\displaystyle\ \ \lim_{t\rightarrow\pm\infty}g(t)=l_{\pm}\in\RR\ \
\hbox{and}\ \ (l_+,l_-)\not =(0,0)$.

b)$\displaystyle\ \ \lim_{|t|\rightarrow\infty}g(t)=0\ \ \hbox{and}\ \
 \lim_{|t|\rightarrow\infty}\int_0^tg(s)ds=\pm\infty$.

c)$\displaystyle\ \ \lim_{|t|\rightarrow\infty}g(t)=0\ \ \hbox{and}\ \
\lim_{|t|\rightarrow\infty}\int_0^tg(s)ds=\beta\in\RR$.

The case c) is called as the case of a  {\it strong resonance}.

For a treatment of these cases, we refer only to
 \cite{LL}, \cite{ALP}, \cite{BBF}, \cite{AM1},
\cite{AM2}, \cite{He}, \cite{Th}, \cite{CS}, \cite{MW2}.

We shall study in what follows a multivalued variant of the
Landesman-Lazer problem.

\neweq{3.33}
\left\{-\Delta u(x)-\lambda_1u(x)\in [\finf (u(x)),
\fsup (u(x))\quad {\rm a.e.}\ \ x\in\Omega\atop u\in\huo\setminus
\{ 0\}\, ,\right.
\endeq
where

i) $\Omega\subset\RR^N$ is an open bounded set with
sufficiently smooth boundary;

ii) $\lambda_1$ (respectively $e_1$) is the first eigenvalue (respectively
 eigenfunction) of the operator
 $(-\Delta )$ in $H^1_0(\Omega )$;

iii) $f\in L^{\infty }(\RR )$;

iv) $\finf (t)=\limeps {\rm essinf}\{ f(s);\ |t-s|<\varepsilon \},\ \
\fsup (t)=\limeps {\rm esssup}
\{ f(s);\ |t-s|<\varepsilon\}\, .$

Consider in \huo the functional $\varphi (u)=\varphi_1(u)-\varphi_2(u)$, where
$$\varphi_1(u)={1\over 2}\intom (|\nabla u|^2-\lambda_1u^2)\quad\ \hbox{and}\quad\
\varphi_2(u)=\intom F(u)\, .$$
Observe first that $\varphi$ is locally Lipschitz in \huo . Indeed,
it is enough to show that $\varphi_2$ is locally Lipschitz in \huo ,
which follows from
$$|\varphi_1(u)-\varphi_2(u)|=|\intom (\int_{u(x)}^{v(x)}f(t)dt)dx|\leq$$
$$\leq \| f\|_{L^\infty}\cdot \| u-v\|_{L^1}\leq C_1\| u-v\|_{L^2}\leq
C_2\| u-v\|_{H^1_0}\, .$$

We shall study the problem \eq{3.33} under the hypothesis
\neweq{f1}
f(+\infty ):=\lim_{t\rightarrow +\infty}f(t)=0,\quad
 F(+\infty ):=\lim_{t\rightarrow +\infty}F(t)=0\, .
 \endeq

Thus, by \cite{LL} and \cite{BBF}, the problem \eq{3.33} becomes a
Landesman-Lazer type problem, with strong resonance at $+\infty$.

As an application of Corollary \ref{1.2.7} we shall prove the following sufficient
condition for the existence of a solution to our problem.

\medskip
\begin{teo}\label{3.2.1} \sl Assume that $f$ satisfies the condition $(f1)$, as well as
\neweq{3.34}
-\infty\leq F(-\infty )\leq 0\, .
\endeq

If $F(-\infty )$ is finite, we assume further that
\neweq{3.35}
\hbox{there exists}\ \ \eta >0\ \ \hbox{such that}\ \ F\ \
\hbox{is nonnegative on}\ \
(0,\eta )\ \ \hbox{or}\ \ (-\eta ,0)\, .
\endeq

Under these hypotheses, the problem \eq{3.33}
 has at least one solution.\rm
 \end{teo}

\medskip
We shall make use in the proof of the following auxiliary results:

\medskip
\begin{lemma}\label{3.2.2}
\sl Assume that $f\in L^\infty (\RR )$ and there exists
$F(\pm \infty )\in\overline{\RR}$. Moreover, assume that

i) $f(+\infty )=0$ if $F(+\infty )$ is finite.

ii) $f(-\infty )=0$ if $F(-\infty )$ is finite.

Under these hypotheses,
$$\RR\subset\ \{\alpha |\Omega |\ ;\ \ \alpha =-F(\pm\infty )\}\ \subset\
\{c\in\RR;\ \varphi\ \hbox{satisface}\ \pslc\}\, .$$
\end{lemma}

\medskip
\begin{lemma}\label{3.2.3}
\sl Assume that $f$ satisfies the condition (f1).
Then $\varphi$ satisfies the condition \pslc\, ,
 for every $c\not =0$ such that
$c<-F(-\infty )\cdot |\Omega |$.\rm
\end{lemma}

\medskip
Assume, for the moment, that these results have been proved.

\medskip
{\bf Proof of Theorem \ref{3.2.1}} There are two distinct situations:

\underbar{Case 1.} $\ \ F(-\infty )$ is finite, that is
$-\infty <F(-\infty )\leq 0$. In this case, $\varphi$ is bounded from below,
 because
$$\varphi (u)={1\over 2}\intom (|\nabla u|^2-\lambda_1u^2)-\intom F(u)\geq -
\intom F(u)$$
and, in virtue of the hypothesis on $F(-\infty )$,
$$\sup_{u\in H^1_0(\Omega )}\intom F(u)<+\infty\, .$$
Hence
$$-\infty <a:=\inf_{H^1_0(\Omega )}\varphi \leq 0=\varphi (0)\, .$$

There exists a real number $c$, sufficiently small in absolute value and
 such that $F(ce_1)<0$ (we observe that $c$ may be chosen positive if $F>0$
in $(0,\eta )$ and negative, if $F<0$ in $(-\eta ,0)$ ). So,
 $\varphi (ce_1)<0$, that is
$a<0$. By Lemma \ref{3.2.3}, it follows that $\varphi$ satisfies the condition
 $({\rm PS})_a$.

\underbar{Case 2.} $\ \ F(-\infty )=-\infty$. Then, by Lemma
\ref{3.2.2}, $\varphi$
satisfies the condition \pslc , for every $c\not =0$.

Let $V$ be the orthogonal complement with respect to \huo of the space spanned by
 $e_1$, that is
$$\huo =\hbox{Sp}\ \{e_1\}\oplus V\, .$$

For some fixed $t_0>0$, let
$$V_0:=\{t_0e_1+v;\ v\in V\}$$
$$a_0:=\inf_{V_0}\varphi\, .$$

We remark that $\varphi$ is coercive on $V$. Indeed, for every $v\in V$,
\neweq{3.36}
\varphi (v)={1\over 2}\, \intom (|\nabla v|^2-\lambda_1v^2)-\intom F(v)\geq\leqno (3.36)$$
$$\geq {{\lambda_2-\lambda_1}\over 2}\| v\|^2_{H^1_0}-\intom F(v)\rightarrow +\infty,\ \ \hbox{if}\ \
\| v\|_{H^1_0}\rightarrow +\infty \, ,
\endeq
because the first term in the right hand side of
 \eq{3.36} has a quadratic growth at infinity
 ($t_0$ being fixed), while $\intom F(v)$ is uniformly bounded
(with respect to $v$), by the behaviour of $F$ near
$\pm\infty$. So, $a_0$ is attained, because of the coercivity of
 $\varphi$ in $V$. Taking into account the boundedness of
$\varphi$ in \huo , it follows that  $-\infty <a\leq0=\varphi (0)$ and
$a\leq a_0$.

At this stage, there are again two possibilities:

i) $\ \ a<0$. Thus, by Lemma \ref{3.2.3},
 $\varphi$ satisfies $(PS)_a$. So,
$a<0$ is a critical value of $\varphi$.

ii)$\ \ a=0\leq a_0$. If $a_0=0$, then, by a preceding remark,
$a_0$ is achieved. So, there exists $v\in V$ such that
$$0=a_0=\varphi (t_0e_1+v)\, .$$
Hence $u=t_0e_1+v\in\huo\setminus\{ 0\}$ is a critical point of
$\varphi$, that is a solution of the problem \eq{3.33}.

If $a_0>0$, we observe that $\varphi$ satisfies $(PS)_b$ for every
$b\not =0$. Since $\displaystyle\lim_{t\rightarrow +\infty}\varphi
(te_1)=0$, we can apply Corollary \ref{1.2.7} for $X=\huo ,\
X_1=V,\ X_2=\hbox{Sp}\ \{e_1\},\ f=\varphi,\ z=t_0e_1$. Therefore,
$\varphi$ has a critical value $c\geq a_0>0$.\qed

\medskip
{\bf Proof of Lemma \ref{3.2.2}} We shall assume, without loss of generality, that
 $F(-\infty )\notin\RR$
and $F(+\infty )\in\RR$. In this case, if $c$ is a critical value such that
$\varphi$ does not satisfy \pslc , then it is enough to prove that
 $c=-F(+\infty )\cdot |\Omega |$.
Since $\varphi$ does not satisfy the condition \pslc , there exist $t_n\in\RR$ and
$v_n\in V$ such that the sequence $(u_n)\subset\huo$, where $u_n=t_ne_1+v_n$, has no
convergent subsequence, while
\neweq{3.37}
\limn\varphi (u_n)=c\endeq
\neweq{3.38}
\limn \lambda (u_n)=0\, .
\endeq

\underbar{Step 1.} {\it The sequence} $(v_n)$ {\it is bounded in}  \huo .

By \eq{3.38} and
$$\partial\varphi (u)=-\Delta u-\lambda_1u-\partial\varphi_2(u)\, ,$$
it follows that there exists $w_n\in\partial\varphi_2(u_n)$ such that
$$-\Delta u_n-\lambda_1u_n-w_n\rightarrow 0\quad\hbox{in}\ \
H^{-1}(\Omega )\, .$$
So
$$\langle -\Delta u_n-\lambda_1u_n-w_n,v_n\rangle =\intom
|\nabla v_n|^2-\lambda_1\intom v_n^2-
\intom g_n(t_ne_1+v_n)=o(\| v_n\|_{H^1_0}),$$
as $n\rightarrow\infty$, where $\finf\leq g_n\leq \fsup$.
Since $f$ is bounded, it follows that
$$\| v_n\|^2_{H^1_0}-\lambda_1\| v_n\|^2_{L^2}=O(\| v_n\|_{H^1_0})\, .$$
So, there exists $C>0$ such that, for every $n\geq 1$, $\| v_n\|_{H^1_0}\leq C$.
Now, since $(u_n)$ has no convergent subsequence,
 it follows that the sequence $(u_n)$ has no convergent subsequence,
 too.

\underbar{Step 2.} $\ \ t_n\rightarrow +\infty$.

Since $\| v_n\|_{H^1_0}\leq C$ and the sequence $(t_ne_1+v_n)$ has no
convergent subsequence,
 it follows that
$|t_n|\rightarrow +\infty$.

On the other hand, by Lebourg's Mean Value Theorem, there exist
$\theta\in (0,1)$ and $x^*\in \partial F(te_1(x)+\theta v(x))$ such that
$$\varphi_2(te_1+v)-\varphi_2(te_1)=\intom \langle x^*,v(x)\rangle dx\leq$$
$$\leq\intom F^0(te_1(x)+v(x),v(x))dx=$$
$$=\displaystyle\intom
\limsup_{{y\rightarrow te_1(x)+v(x)}\atop
{\lambda\searrow 0}}
{{F(y+\lambda v(x))-F(y)}\over\lambda}dx\leq$$
$$\leq \| f\|_{L^\infty}\cdot\intom |v(x)|dx=\| f\|_{L^\infty}\cdot \| v\|_{L^1}\leq
C_1\| v\|_{H^1_0}\, .$$
A similar computation for $\varphi_2(te_1)-\varphi_2(te_1+v)$ together with the
above inequality shows that, for every $t\in\RR$ and $v\in V$,
$$|\varphi_2(te_1+v)-\varphi_2(te_1)|\leq C_2\| v\|_{H^1_0}\, .$$
So, taking into account the boundedness of $(v_n)$ in \huo , we find
$$|\varphi_2(t_ne_1+v_n)-\varphi_2(t_ne_1)|\leq C\, .$$
Therefore, since $F(-\infty )\notin\RR$ and
$$\varphi (u_n)=\varphi_1(v_n)-\varphi_2(t_ne_1+v_n)\rightarrow c,$$
it follows that  $t_n\rightarrow +\infty$. In this argument we have also used
the fact that
 $\varphi_1(v_n)$ is bounded.

\underbar {Step 3.} $\ \ \|v_n\|_{H^1_0}\rightarrow 0$ if $n\rightarrow\infty$.

By (f1) and Step 2 it follows that
$$\limn\intom f(t_ne_1+v_n)v_n=0\, .$$

Using now \eq{3.38} and Step 1 we find
$$\limn \| v_n\|_{H^1_0}=0\, .$$

\underbar{Step 4.}
\neweq{3.39}
\lim_{t\rightarrow +\infty}\varphi_2(te_1+v)=
F(+\infty )\cdot |\Omega | \, ,
\endeq
uniformly on the bounded subsets of  $V$.
Assume the contrary. So, there exist $r>0,\ t_n\rightarrow +\infty ,
\ v_n\in V$ with $\| v_n\|\leq r$,
such that (3.39) is not fulfilled. Thus, up to a subsequence, there exist
$v\in\huo$ and $h\in L^2(\Omega )$ such that
\neweq{3.40}
v_n\rightharpoonup v\quad\hbox{weakly in}\ \ \huo\, ,
\endeq
\neweq{3.41}
v_n\rightarrow v\quad\hbox{strongly in}\ \ L^2(\Omega )\, ,
\endeq
\neweq{3.42}
v_n(x)\rightarrow v(x)\quad\hbox{a.e.}\ \ x\in\Omega\, ,
\endeq
\neweq{3.43}
|v_n(x)|\leq h(x)\quad\hbox{a.e.}\ \ x\in\Omega \, .
\endeq

For any $n\geq 1$ we define
$$A_n=\{ x\in\Omega ;\ t_ne_1(x)+v_n(x)<0\}\, ,$$
$$h_n(x)=F(t_ne_1+v_n)\chi_{A_n}\, ,$$
where $\chi_{A}$ represents the characteristic function of the set
 $A$. By \eq{3.43} and the choice of $t_n$
it follows that  $|A_n|\rightarrow 0$ if $n\rightarrow\infty$.

Using \eq{3.42} we remark easily that
$$h_n(x)\rightarrow 0\quad\hbox{a.e.}\ \ x\in\Omega\, .$$
Therefore
$$|h_n(x)|=\chi_{A_n}(x)\cdot |\int_0^{t_ne_1(x)+v_n(x)}f(s)ds|\leq$$
$$\leq\chi_{A_n}(x)\cdot \| f\|_{L^\infty}\cdot |t_ne_1(x)+v_n(x)|\leq C|v_n(x)|
\leq Ch(x)\quad\hbox{a.e.}\ \ x\in\Omega\, .$$
So, by Lebesgue's Dominated Convergence Theorem,
$$\limn \int_{A_n}F(t_ne_1+v_n)=0\, .$$
On the other hand,
$$\limn \int_{\Omega\setminus A_n} F(t_ne_1+v_n)=F(+\infty )
\cdot |\Omega |\, .$$
So
$$\limn \varphi_2(t_ne_1+v_n)=\limn\intom F(t_ne_1+v_n)=F(+\infty )
\cdot |\Omega | ,$$
which contradicts our initial assumption.

\underbar{Step 5.} Taking into account the preceding step and the fact that
$\varphi (te_1+v)=\varphi_1(v)-\varphi_2(te_1+v)$, we obtain
$$\limn \varphi (t_ne_1+v_n)=\limn\varphi_1(v_n)-\limn\varphi_2(t_ne_1+v_n)=
-F(+\infty )\cdot |\Omega | ,$$ that is $c=-F(+\infty )\cdot
|\Omega |$, which concludes our proof. \qed

\medskip
{\bf Proof of Lemma \ref{3.2.3}} It is enough to show that for every
 $c\not= 0$ and $(u_n)\subset\huo$ such that
$$\varphi (u_n)\rightarrow c\, ,$$
\neweq{3.45}
\lambda (u_n)\rightarrow 0\, ,
\endeq
$$\| u_n\|\rightarrow\infty \, ,$$
we have $c\geq -F(-\infty )\cdot |\Omega |$.

Let $t_n\in\RR$ and $v_n\in V$ be such that, for every $n\geq 1$,
$$u_n=t_ne_1+v_n\, .$$
As we have already remarked,
$$\varphi (u_n)=\varphi_1(v_n)-\varphi_2(u_n)\, .$$
Moreover,
$$\varphi_1(v)={1\over 2}\intom (|\nabla v|^2-\lambda_1v^2)\geq
{1\over 2}(1-{{\lambda_1}\over{\lambda_2}})
\cdot\| v\|^2_{H^1_0}\rightarrow +\infty\quad\hbox{if}\ \ \| v\|_{H^1_0}
\rightarrow\infty\, .$$
So, $\varphi_1$ is positive and coercive on $V$. We also have that
$\varphi_2$ is bounded from below, by eq{f1}. So,
again by \eq{f1}, we conclude that the sequence $(v_n)$ is bounded
in \huo . Thus there exists $v\in\huo$ such that, up to a subsequence,
$$v_n\rightharpoonup v\quad\hbox{weakly in}\ \ \huo\, ,$$
$$v_n\rightarrow v \quad\hbox{strongly in}\ \ L^2(\Omega )\, ,$$
$$v_n(x)\rightarrow v(x)\quad\hbox{a.e.}\ \ x\in\Omega\, ,$$
$$|v_n(x)|\leq h(x)\quad\hbox{a.e.}\ \ x\in\Omega \, ,$$
for some $h\in L^2(\Omega )$.

Since $\| u_n\|_{H^1_0}\rightarrow\infty$ and $(v_n)$ is bounded
in \huo, it follows that  $|t_n|\rightarrow +\infty$.

Assume for the moment that we have already proved that
 $\| v_n\|_{H^1_0}\rightarrow 0$,
if $t_n\rightarrow +\infty$. So,
$$\varphi (u_n)=\varphi_1(v_n)-\varphi_2(u_n)\rightarrow
0\quad\hbox{if}\ \ n\rightarrow\infty\, .$$

 Here, for proving that $\varphi_2(u_n)\rightarrow 0$, we have used (f1).
The last relation yields a contradiction, since
 $\varphi (u_n)\rightarrow c\not= 0$.
So, $t_n\rightarrow -\infty$.

Moreover, since $\varphi (u)\geq -\varphi_2(u)$ and $F$ is bounded from below,
 it follows that
$$c=\liminf_{n\rightarrow\infty}\varphi (u_n)\geq\liminf_{n\rightarrow\infty}
(-\varphi_2(u_n))=$$
$$=-\limsup_{n\rightarrow\infty}\intom F(u_n)\geq -\intom\limsup_{n\rightarrow\infty}F(u_n)=
-F(-\infty )\cdot |\Omega | ,$$
which gives the desired contradiction.

So, for concluding the proof, we have to show that
$$\| v_n\|_{H^1_0}\rightarrow 0\quad\hbox{if}\ \
t_n\rightarrow +\infty\, .$$
Since
$$\partial\varphi (u)=-\Delta u-\lambda_1 u-\partial\varphi_2(u),$$
it follows from \eq{3.45} that there exists
 $w_n\in\partial\varphi_2(u_n)$ such that
$$-\Delta u_n-\lambda_1u_n-w_n\rightarrow 0\quad\hbox{in}\ \
H^{-1}(\Omega )\, .$$
Thus
$$\langle -\Delta u_n-\lambda_1u_n-w_n,v_n\rangle =\intom
|\nabla v_n|^2-\lambda_1\intom v_n^2-$$
$$-\intom g_n(t_ne_1+v_n)v_n=o(\| v_n\|)\quad\hbox{if}\
\ n\rightarrow\infty ,$$
where $\finf \leq g_n\leq \fsup$.

Now, for concluding the proof, it is sufficient to show that the last
 term tends to 0,
as $n\rightarrow\infty$.

Let $\varepsilon >0$. Because $f(+\infty )=0$, it follows that
there exists $T>0$ such that
$$|f(t)|\leq\varepsilon\quad\quad\hbox{a.e.}\ \ t\geq T\, .$$

Set
$$A_n=\{ x\in\Omega ; \ t_ne_1(x)+v_n(x)\geq T\}\quad\hbox{and}
\quad B_n=\Omega\setminus A_n\, .$$

We remark that for every $x\in B_n$,
$$|t_ne_1(x)+v_n(x)|\leq |v_n(x)|+T\, .$$
So, for every $x\in B_n$,
$$|g_n(t_ne_1(x)+v_n(x))v_n(x)|\cdot\chi_{B_n}(x)\leq\|
f\|_{L^\infty}\cdot h(x)\, .$$
By
$$\chi_{B_n}(x)\rightarrow 0\quad\hbox{a.e.}\ \ x\in\Omega $$
and the Lebesgue Dominated Convergence Theorem it follows that
\neweq{3.47}
\int_{B_n}g_n(t_ne_1+v_n)v_n\rightarrow 0\quad\hbox{if}\ \
n\rightarrow\infty\, .
\endeq

On the other hand, it is obvious that
\neweq{3.48}
|\int_{A_n}g_n(t_ne_1+v_n)v_n|\leq\varepsilon\int_{A_n}|v_n|
\leq\varepsilon \, \| h\|_{L^1}\, .
\endeq

By \eq{3.47} and \eq{3.48} it follows that
$$\limn\intom g_n(u_n)v_n=0,$$
which concludes our proof. \qed

\section{Multivalued problems of
 Landesman-Lazer type with mixed resonance}

\hspace*{6mm}Let $X$ be a Banach space. Assume that there exists
 $w\in X\setminus\{ 0\}$, which can be supposed to have the norm 1, and a linear
 subspace
 $Y$ of $X$ such that
\neweq{3.49}
X=\hbox{Sp}\{ w\}\oplus Y\, .
\endeq

\medskip
\begin{defin}\label{3.3.1}
\sl A subset $A$ of $X$ is called to be $w$-bounded
if there exists $r\in\RR$ such that
$$A\subset\{ x=tw+y;\ t<r,y\in Y\}\, .$$

A functional \fxr is said to be $w$-coercive (or, coercive with respect to the
decomposition \eq{3.49}) if
$$\lim_{t\rightarrow\infty}f(tw+y)=+\infty ,$$
uniformly with respect to $y\in Y$.

The functional $f$ is called $w$-bounded from below if there exists $a\in\RR$
such that the set $[f\leq a]$ is $w$-bounded.
\rm
\end{defin}

\medskip
All the results from 1.3 can be extended in this more general framework. For example,
we shall give the variant of Proposition
 \ref{1.3.7}:

\medskip
\begin{prop}\label{3.3.2}
\sl Let $f$ be a $w$-bounded from below locally Lipschitz
functional. Assume that
 there exists $c\in\RR$ such that $f$ satisfies the condition
(s-PS)$_c$ , and the set $[f\leq a]$ is $w$-bounded from below for every $a<c$.

Then there exists $\alpha >0$ such that the set
 $[f\leq c+\alpha ]$ is $w$-bounded.\rm
 \end{prop}

\medskip
{\bf Proof.} Assume the contrary. So, for every $\alpha >0$,
the set $[f\leq c+\alpha ]$  is not $w$-bounded. Thus, for every
$n\geq 1$, there exists $r_n\geq n$ such that
$$[f\leq c-{1\over{n^2}}]\subset A_n:=\{ x=tw+y;\ t<r_n,\ y\in Y\}\, .$$
Therefore
\neweq{3.50}
c_n:=\inf_{X\setminus A_n}f\geq c-{1\over{n^2}}\, .
\endeq

Since the set $[f\leq c+{1\over{n^{\scriptstyle 2}}}]$  is not
$w$-bounded, there exists a sequence $(z_n)$ in
$X$ such that $z_n=t_nw+y_n$ and
$$t_n\geq r_n+1+{1\over n}\, ,$$
\neweq{3.51}
f(z_n)\leq c+{1\over{n^2}}\, .
\endeq
It follows that $z_n\in X\setminus A_n$ and, by
\eq{3.50} and \eq{3.51}, we find
$$f(z_n)\leq c_n+{2\over{n^2}}\, .$$

We apply Ekeland's Variational Principle to $f$ restricted to the set
$X\setminus A_n$, provided $\varepsilon ={2\over{n^{\scriptstyle{2}}}}$ and
$\lambda ={1\over n}$. So, there exists $x_n=t'_nw+y'_n\in X\setminus A_n$ such that,
for every $x\in X\setminus A_n$,
$$c\leq f(x_n)\leq f(z_n)\, ,$$
\neweq{3.52}
f(x)\geq f(x_n)-{2\over n}\, \| x-x_n\|\, ,
\endeq
\neweq{3.53}
\| x_n-z_n\|\leq{1\over n}\, .
\endeq

Let $P$ ($\| P\| =1$) be the projection of $X$ on Sp $\{ w\}$. Using the continuity
of $P$
and the relation \eq{3.53} we obtain
$$|t_n-t'_n|\leq{1\over n}\, .$$
So,
$$|t'_n|\geq r_n+1\, .$$
Therefore $x_n$ is an interior point of the set $A_n$. By \eq{3.52}
it follows that, for every $v\in X$,
$$f^0(x_n,v)\geq -{2\over n}\, \| v\|\, .$$
This relation and the fact that $f$ satisfies the condition \spsc imply
 that the sequence
$(x_n)$ contains a  convergent subsequence, contradiction, because
this sequence is not $w$-bounded. \qed

\medskip
As an application of these results we shall study the following multivalued
 Landesman-Lazer problem with mixed resonance.

Let $\Omega\subset \RR^N$ be an open, bounded set with the boundary sufficiently
smooth. Consider
the problem
\neweq{3.54}
\left\{ -\Delta u(x)-\lambda_1u(x)\in [\finf (u(x)),\fsup (u(x))]\quad a.e.\
 \ x\in\Omega\atop
u\in\huo\, .\right.
\endeq

We shall study this problem under the following hypotheses:

(f1) $\quad f\in L^\infty (\RR )\, ;$

(f2) $\quad \hbox{if}\ \ F(t)=\int_0^1f(s)ds$, then
 $\displaystyle\lim_{t\rightarrow +\infty}F(t)=0\, ;$

(f3) $\quad\displaystyle\lim_{t\rightarrow -\infty}F(t)=+\infty\, ;$

(f4) $\quad\hbox{there exists}\ \ \alpha <{1\over 2}(\lambda_2-\lambda_1)$ such that,
for every $t\in\RR$, $F(t)\leq\alpha t^2$.

Define on the space \huo the functional
 $\varphi (u)=\varphi_1(u)-\varphi_2(u)$, where
$$\varphi_1(u)={1\over 2}\intom (|\nabla u|^2-\lambda_1u^2)dx\, ,$$
$$\varphi_2(u)=\intom F(u)dx\, .$$

As in the preceding section, we observe that the functional $\varphi$ is
locally Lipschitz.

Let $Y$ be the orthogonal complement of the space spanned by $e_1$, that is
$$\huo =\hbox{Sp}\{ e_1\} \oplus Y\, .$$

\medskip
\begin{lemma}\label{3.3.3} \sl With the above notations, The following hold:

i) there exists $r_0>0$ such that, for every $t\in\RR$ and $v\in Y$ with
 $\| v\|_{H^1_0}\geq r_0$,
$$\varphi (te_1+v)\geq -\varphi_2(te_1)\, ;$$

$$\lim_{t\rightarrow +\infty}\varphi (te_1+v)=\varphi_1(v),\leqno ii)$$
uniformly on the bounded subsets of  $Y$.\rm
\end{lemma}

\medskip
{\bf Proof.} i) We first observe that, for every $t\in\RR$ and
$v\in Y$,
$$\varphi_1(te_1+v)=\varphi_1(v)\, .$$
So,
\neweq{3.55}
\varphi (te_1+v)=\varphi_1(v)-\varphi_2(te_1+v)\, .
\endeq
On the other hand,
\neweq{3.56}
\varphi_2(te_1+v)-\varphi_2(te_1)=\intom (\int_{te_1(x)}^{te_1(x)+v(x)}
f(s)ds)dx=
\endeq
$$=\intom (F(te_1(x)+v(x))-F(te_1(x)))dx\, .$$

By the Lebourg Mean Value Theorem, there exist $\theta\in (0,1)$ and
$x^*\in\partial F(te_1(x)+\theta v(x))$ such that
$$F(te_1(x)+v(x))-F(te_1(x))=\langle x^*,v(x)\rangle\, .$$

The relation \eq{3.56} becomes
\neweq{3.57}
\left\{\begin{array}{lll}
 \varphi_2(te_1+v)-\varphi_2(te_1)=\intom \langle x^*,v(x)\rangle dx\leq
\intom F^0(te_1(x)+\theta v(x),v(x))dx= \\
=\displaystyle\intom\limsup_{{y\rightarrow te_1(x)+\theta v(x)}\atop{\lambda\searrow 0}}
{{F(y+\lambda v(x))-F(y)}\over\lambda}dx\leq  \\
 \leq\| f\|_{L^\infty}\cdot\intom |v(x)|dx=\|
f\|_{L^\infty}\cdot \| v\|_{L^1}\leq C\| v\|_{H^1_0} \, .
\end{array}\right.
\endeq

By \eq{3.55} and \eq{3.57} it follows that
$$\varphi (te_1+v)\geq\varphi_1(v)-C\| v\|_{H^1_0}-\varphi_2(te_1)\, .$$
So, for concluding the proof, it is enough to choose
$r_0>0$ such that
$$\varphi_1(v)-C\| v\|_{H^1_0}\geq 0,$$
for every $v\in Y$ cu $\| v\|_{H^1_0}\geq r_0$. This choice is possible if
we take into account the variational characterization of $\lambda_2$ and
the fact that $\lambda_1$ is a simple eigenvalue. Indeed,
$$\varphi_1(v)-C\| v\|_{H^1_0}={1\over 2}(\| \nabla v\|^2_{L^2}-
\lambda_1\| v\|^2_{L^2})-C\| \nabla v\|_{L^2}\geq$$
$$\geq ({1\over 2}-\varepsilon )\| \nabla v\|^2_{L^2}-{{\lambda_1}\over 2}\| v\|^2_{L^2}\geq$$
$$\geq ({1\over 2}-\varepsilon )\lambda_2\| v\|^2_{L^2}-{{\lambda_1}\over 2}\| v\|^2_{L^2}\geq 0,$$
for $\varepsilon >0$ sufficiently small
 and $\| v\|_{L^2}$ (hence, $\| v\|_{H^1_0})$
large enough.

ii) By \eq{3.55}, our statement is equivalent with
\neweq{3.58}
\lim_{t\rightarrow +\infty}\varphi_2(te_1+v)=0\, ,
\endeq
uniformly with respect to $v$, in every closed ball.

Assume, by contradiction, that there exist
 $R>0,\ t_n\rightarrow +\infty $ and
$v_n\in Y$ with $\| v_n\|\leq R$ such that \eq{3.58}
does not hold. So, up to a subsequence,
we may assume that
 there exist $v\in\huo$ and $h\in L^1(\Omega )$ such that
\neweq{3.59}
v_n\rightharpoonup v\quad\hbox{weakly in}\ \ \huo\, ,
\endeq
\neweq{3.60}
v_n\rightarrow v\quad\hbox{strongly in}\ \ L^2(\Omega )\, ,
\endeq
\neweq{3.61}
v_n(x)\rightarrow v(x)\quad\hbox{for a.e.}\ \ x\in\Omega\, ,
\endeq
\neweq{3.62}
|v_n(x)|\leq h(x)\quad\hbox{a.e.}\ \ x\in\Omega\, .
\endeq

For every $n\geq 1$, denote
$$A_n=\{ x\in\Omega ;\ t_ne_1(x)+v_n(x)<0\}\, ,$$
$$g_n=F(t_ne_1+v_n)\chi_{A_n}\, .$$

By \eq{3.62} and the choice of
 $t_n$ it follows that  $|A_n|\rightarrow 0$.

Using now \eq{3.61}, it is easy to observe that
$$g_n(x)\rightarrow 0\quad\hbox{a.e.}\ \ x\in\Omega\, .$$

By (f1) and \eq{3.62} it follows that
$$|g_n(x)|=\chi_{A_n}(x)\cdot |\int_0^{t_ne_1(x)+v_n(x)}f(s)ds|\leq
\chi_{A_n}(x)\cdot\| f\|_{L^\infty}\cdot |t_ne_1(x)+v_n(x)|\leq$$
$$\leq C|v_n(x)|\leq Ch(x)\quad\hbox{a.e.}\ \ x\in\Omega\, .$$

Thus, by Lebesgue's Dominated Convergence Theorem,
$$\limn\int_{A_n}F(t_ne_1+v_n)dx=0\, .$$

By (f2) it follows that  $F$ is bounded on $[0,\infty )$.
Using again (f2) we find
$$\limn\int_{\Omega\setminus A_n}F(t_ne_1+v_n)dx=0\, .$$

So
$$\limn \varphi_2(t_ne_1+v_n)=\limn\intom F(t_ne_1+v_n)dx=0,$$
which contradicts our initial assumption. \qed

\begin{rem}\label{3.3.4}
As a consequence of the above result,
$$\lim_{t\rightarrow +\infty}\varphi (te_1)=0\, ,$$
$$\liminf_{t\rightarrow +\infty}\inf_{v\in Y}\varphi
(te_1+v)\geq 0\, .$$

Thus, the set $[\varphi\leq a]$ is $e_1$-bounded if $a<0$
and  is not $e_1$-bounded for $a>0$. Moreover,
$\varphi$  is not bounded from below, because
$$\lim_{t\rightarrow -\infty}\varphi (te_1)=-\lim_{t\rightarrow -\infty}
\varphi_2(te_1)=-\lim_{t\rightarrow -\infty}\intom F(te_1)=-\infty\, .$$

Thus, by Proposition \ref{3.3.2}, it follows that
$\varphi$ does not satisfy the condition $(s - \hbox{PS})_0$.
\end{rem}

\medskip
\begin{teo}\label{3.3.5} \sl Assume that $f$ does not satisfy the conditions
(f1)-(f4). If \fla $\varphi$ has the strong Palais-Smale property
$(s-\hbox{PS})_a$, for every $a\not= 0$, then the multivalued problem
\eq{3.54} has at least a nontrivial solution.\rm
\end{teo}

\medskip
{\bf Proof.} It is sufficient to show that $\varphi$ has a critical point
 $u_0\in \huo \setminus\{ 0\}$. It is obvious that
$$\partial\varphi (u)=-\Delta u-\lambda_1u-\partial\varphi_2(u)\quad
\hbox{in}\ \ H^{-1}(\Omega )\, .$$

If $u_0$ is a critical point of $\varphi$, then there exists
$w\in\partial\varphi_2(u_0)$ such that
$$-\Delta u_0-\lambda_1u_0=w\quad\mbox{in} \ \ H^{-1}(\Omega )\, .$$

Set
$$Y_1=\{ e_1+v;\ v\in Y\}\, .$$

We first remark the coercivity of  $\varphi$ on $Y_1$:
$$\varphi (e_1+v)={1\over 2}(\| \nabla\|^2_{L^2}-\lambda_1\| v\|^2_{L^2})
-\lambda_1\intom e_1vdx-\intom G(e_1+v)dx\geq$$
$$\geq {1\over 2}(\| \nabla v\|^2_{L^2}-\lambda_1\| v\|^2_{L^2})-
\lambda_1\| v\|_{L^2}-
\alpha (1+\| v\|^2_{L^2}+2\intom e_1vdx)\geq$$
$$\geq ({{\lambda_2-\lambda_1}\over 2}-\alpha )\| v\|^2_{L^2}-(\lambda_1+
2|\alpha |)\| v\|^2_{L^2}-\alpha ,$$
which tends to $+\infty$ as $\| v\|_{L^2}\rightarrow\infty$, so,
 as $\| v\|_{H^1_0}\rightarrow\infty$.

Putting
$$m_1=\inf_{Y_1}\varphi >-\infty ,$$
it follows that  $m_1$ is attained, by the coercivity of
$\varphi$ pe $Y_1$. Thus, there exists $u_0=e_1+v_0\in Y_1$ such that
$\varphi (u_0)=m_1$. There are two possibilities:

1) $\ \ m_1>0$. It follows from (f1) and (f3) that
$$\lim_{t\rightarrow +\infty}\varphi (te_1)=-\lim_{t\rightarrow +\infty}
\varphi_2(te_1)=0\, ,$$
$$\lim_{t\rightarrow -\infty}\varphi (te_1)=-\lim_{t\rightarrow -\infty}
\varphi_2(te_1)=-\infty\, .$$

Since $m_1>0=\varphi (0)$ and $\varphi$ has the property
$(s-\hbox{PS})_a$ for every $a>0$, it follows from Corollary \ref{1.2.7} that
 $\varphi$ has a critical value $c\geq m_1>0$.

2) $\ \ m_1\leq 0$. Let
$$W=\{ te_1+v;\ t\geq 0, \ v\in Y\}\, ,$$
$$c=\inf_{W}\varphi\, .$$
So, $c\leq m_1\leq 0$.

If $c<0$, it follows by \spsc and the Ekeland Principle that $c$
is a critical value of $\varphi$.

If $c=m_1=0$, then $u_0$ is a local minimum point of $\varphi$, because
$$\varphi (u_0)=\varphi (e_1+v_0)=0\leq\varphi (te_1+v),$$
for every $t\geq 0$ and $v\in Y$.

So, $u_0\not= 0$ is a critical point of $\varphi$. \qed


\begin{thebibliography}{999}
{\footnotesize
\bibitem{ab} G. A. Afrouzi, K. J. Brown, Positive solutions for a semilinear
elliptic problem with sign-changing nonlinearity, {\it Nonlinear
Analysis, T.M.A.} {\bf 36} (1999), 507-510.

\bibitem{ALP} S. Ahmad, A. C. Lazer, J. L. Paul, Elementary critical point theory and perturbations of elliptic
boundary value problems at resonance, {\it Indiana Univ. Math. J.}
 {\bf 25} (1976), 933-944.

\bibitem{a} H. Amann, On the existence of positive solutions of nonlinear
elliptic boundary value problems, {\it Indiana Univ. Math.~J.}
 {\bf 21} (1971), 125-146.

\bibitem{Am} H. Amann, Multiple positive fixed points of assymptotically linear maps,
{\it J. Funct. Anal.} {\bf 14} (1973), 162-171.

\bibitem{AC} H. Amann, M. Crandall, On some existence theorems for semilinear elliptic
equations, {\it Indiana Univ. Math. J.} {\bf 27} (1978), 779-790.

\bibitem{A} A. Ambrosetti, Critical points and nonlinear
 variational problems, {\it Bull.  Soc.
Math. France} {\bf 120} (1992), 5-139.

\bibitem{AM1} A. Ambrosetti, G. Mancini, Theorems of existence and multiplicity for nonlinear
elliptic problems with noninvertible linear part, {\it Ann. Sc.
Norm. Sup. Pisa.} {\bf 5} (1978), 15-38.

\bibitem{AM2} A. Ambrosetti, G. Mancini, Existence and multiplicity results for nonlinear elliptic
problems with linear part at resonance, {\it J. Diff. Eq.} {\bf
28} (1978), 220-245.

\bibitem{AR} A. Ambrosetti, P. H. Rabinowitz, Dual variational methods in critical
point theory and applications, {\it J. Funct. Anal.} {\bf 14}
(1973), 349-381.

\bibitem{Ar} A. Arcoya, Periodic solutions of Hamiltonian systems with strong
resonance at infinity, {\it Diff.  Int. Equations} {\bf 3} (1990),
909-921.

\bibitem{ArC1} A. Arcoya, A. Canada, Critical point theorems and applications to
nonlinear boundary value problems, {\it J. Nonlinear Anal. TMA}
{\bf 14} (1990), 393-411.

\bibitem{ArC2} A. Arcoya, A. Canada, The dual variational principle and discontinuous elliptic
problems with strong resonance at infinity, {\it  Nonlinear Anal.
TMA} {\bf 15} (1990), 1145-1154.

\bibitem{AUC} J. P. Aubin, F. H. Clarke, Shadow prices and duality for a class of optimal control problems,
{\it SIAM J. Control and Optimization} {\bf 17} (1979), 567-586.

\bibitem{BBF} P. Bartolo, V. Benci, D. Fortunato, Abstract critical point theorems
and applications to some nonlinear problems with strong resonance
at infinity, {\it Nonlinear Analysis}, T. M. A. {\bf 7} (1983),
981-1012.

\bibitem{BR} V. Benci, P. H. Rabinowitz, Critical point theorems for indefinite functionals,
{\it Invent. Math.} {\bf 52} (1979), 241-273.

\bibitem{Be} S. Berberian, {\it Measure and Integration}, Mac Millan, 1967.

\bibitem{BP} M. S. Berger, E. Podolak, On the solutions of nonlinear Dirichlet problem,
{\it Indiana Univ. Math. J.} {\bf 24} (1975), 837-846.

\bibitem{Br} H. Brezis, {\it Analyse Fonctionnelle}, Masson, Paris, 1992.

\bibitem{BCN} H. Brezis, J. M. Coron, L. Nirenberg, Free vibrations for a nonlinear wave
equation and a theorem of Rabinowitz, {\it Comm. Pure Appl. Math.}
{\bf 33} (1980), 667-689.

\bibitem{BN2} H. Brezis, L. Nirenberg, Remarks on finding critical points,
{\it Comm. Pure Appl. Math.} {\bf 44} (1991), 939-964.

\bibitem{BN1} H. Brezis, L. Nirenberg, {\it Nonlinear Functional Analysis and Applications to Partial
Differential Equations}, in preparation.

\bibitem{cfg} N. Cac, A.  Fink, J. Gatica, Nonnegative solutions of the
radial Laplacian with nonlinearity which changes sign, {\it Proc.
Amer. Math. Soc.} {\bf 123} (1995), 1393-1398.

\bibitem{CLW} L. Caklovici, S. Li, M. Willem, A note on Palais-Smale condition versus
coercivity, {\it Diff. and Int. Eq. } {\bf 3} (1990), 799-800.

\bibitem{Ch} K. C. Chang, Variational methods for non-differentiable functionals and
applications to partial diff. equations, {\it J. Math. Anal.
Appl.} {\bf 80} (1981), 102-129.

\bibitem{CDR} M. Choulli, R. Deville, A. Rhandi, A general Mountain Pass principle for
non-differentiable functionals and applications, {\it Rev. Mat.
Apl.} {\bf 13} (1992), 45-58.

\bibitem{Cl1} F. H. Clarke, Generalized gradients of Lipschitz functionals, {\it Adv.
in Math.} {\bf 40} (1981), 52-67.

\bibitem{Cl2} F. H. Clarke, Generalized gradients and applications, {\it Trans. Amer. Math.
Soc.} {\bf 205} (1975), 247-262.

\bibitem{CO} C. Coffman, A min-max principle for a class of nonlinear integral
equations, {\it J. d'Analyse Math.} {\bf 22} (1969), 391-418.

\bibitem{CG} D. G. Costa, J. E. Gon\c calves, Critical point theory for nondifferentiable
functionals and applications, {\it J. Math. Anal. Appl.} {\bf 153}
(1990), 470-485.

\bibitem{CS} D. G. Costa, E. A. De B. Silva, The Palais-Smale condition versus coercivity,
{\it Nonlinear Analysis, T. M. A.} {\bf 16} (1991), 371-381.

\bibitem{cr} M. Crandall, P. H. Rabinowitz, Bifurcation from simple eigenvalues,
{\it J.~Functional Anal.} {\bf  8} (1971), 321-340.

\bibitem{CR} M. Crandall, P. H. Rabinowitz, Some continuation and variational methods
for positive solutions of nonlinear elliptic eigenvalues problems,
{\it Arch. Rational Mech. Anal.} {\bf 58} (1975), 207-218.

\bibitem{dibe}  E. DiBenedetto, {\it Partial Differential
Equations}, Birkh\"auser, Boston Basel Berlin, 1995.

\bibitem{D} J. Dugundji, An extension of Tietze's theorem, {\it Pacific J. Math.} {\bf
1} (1951), 353-367.

\bibitem{E1} I. Ekeland, On the variational principle, {\it J. Math. Anal. Appl.} {\bf
47} (1974), 324-353.

\bibitem{E2} I. Ekeland, Nonconvex minimization problems, {\it Bull. Amer. Math. Soc.} {\bf 1} (1979), 443-474.

\bibitem{evans}  L. C. Evans, {\it Partial Differential Equations},
Graduate Studies in Mathematics, Vol.~19, Amer. Math. Soc.,
Providence, RI, 1998.

\bibitem{FNSS} S. Fucik, J. Necas, J. Soucek, V. Soucek, {\it Spectral Analysis of
Nonlinear Operators}, Springer-Verlag, Heidelberg, 1973.

\bibitem{GaR} F. Gazzola, V. R\u adulescu,
A nonsmooth critical point theory approach to some nonlinear
elliptic equations in unbounded domains, {\it Differential and
Integral Equations} {\bf 13} (2000), 47-60.

\bibitem{G} I. M. Gelfand, Some problems in the theory of quasi-linear equations,
{\it Uspekhi Mat. Nauk} {\bf 14} (1959), 87-158.

\bibitem{GP} N. Ghoussoub, D. Preiss, A general mountain pass principle for locating
and classifying critical points, {\it Ann. Inst. H. Poincare} {\bf
6} (1989), 321- 330.

\bibitem{GT} D. Gilbarg, N. S. Trudinger, {\it Elliptic Partial Differential Equations
of Second Order}, Springer-Verlag, Berlin Heidelberg New York,
1983.

\bibitem{He} P. Hess, Nonlinear perturbations of linear elliptic and parabolic problems at resonance,
 {\it Ann. Sc. Norm. Sup. Pisa} {\bf 5} (1978), 527-537.

 \bibitem{HU} J.B. Hiriart-Urruty, A short proof of the variational principle
 for approximate solutions of a minimization problem, {\it American Math. Monthly} {\bf
 90} (1983), 206-207.

\bibitem{Ho} L. H\"ormander, {\it The Analysis of Linear Partial Differential Operators},
I, Springer-Verlag, Berlin Heildelberg New York, 1983.

\bibitem{If} V. Iftimie, {\it Partial Differential Equations},
 Bucharest University Press 1980 (in Romanian).

\bibitem{JL} D. D. Joseph, T. S. Lundgren, Quasilinear Dirichlet problems driven by
positive sources, {\it Arch. Rat. Mech. Anal.} {\bf 49} (1973),
241-269.

\bibitem{jost} {\sc J. Jost}, {\em Postmodern Analysis}, Springer-Verlag, Berlin Heidelberg New York,
1998.

\bibitem{Kr} M. A. Krasnoselski, {\it Topological Methods in the Theory of Nonlinear
Integral Equations}, Mac Millan, 1968.

\bibitem{Ku} K. Kuratovski, {\it Topologie}, 1934.

\bibitem{LL} E. A. Landesman, A. C. Lazer, Nonlinear perturbations of linear elliptic
boundary value problems at resonance, {\it J. Math. Mech.} {\bf
19} (1976), 609-623.

\bibitem{Le} G. Lebourg, Valeur moyenne pour gradient generalis\'e,
 {\it C. R. Acad. Sci. Paris} {\bf 281} (1975), 795-797.

\bibitem{LeR} C. Lefter, V. R\u{a}dulescu,
Minimization problems and corresponding renormalized energies, {
\it Diff. Integral Equations} {\bf 9} (1996), 903-918.

\bibitem{Li} S. J. Li, Some aspects of critical point theory, Preprint 1986.

\bibitem{LS2} L. Lusternik, L. Schnirelmann, {\it M\'ethodes Topologiques dans les
Probl\`{e}mes Variationnels}, Hermann, Paris, 1934.

\bibitem{MW2} J. Mawhin, M. Willem, Multiple solutions of the periodic boundary value
problem for some forced pendulum-type equations, {\it J. Diff.
Eq.} {\bf 52} (1984), 264-287.

\bibitem{MW1} J. Mawhin, M. Willem, {\it Critical Point Theory and Hamiltonian
Systems}, Springer-Verlag, Berlin Heidelberg New York, 1989.

\bibitem{MR1} P. Mironescu, V. R\u adulescu, A bifurcation problem associated to a convex
asymptotically linear function, {\it C. R. Acad. Sci. Paris} {\bf
316} (1993), 667-672.

\bibitem{MR2} P. Mironescu, V. R\u adulescu, The study of a bifurcation problem associated
to an asymptotically linear function, {\it Nonlinear Analysis,
T.M.A. } {\bf 26} (1996), 857-875.

\bibitem{MR3} P. Mironescu, V. R\u adulescu, A multiplicity theorem for
locally Lipschitz periodic functionals,
 {\it J. Math. Anal. Appl.} {\bf 195} (1995), 621-637.

\bibitem{monradul} E. Montefusco, V. R\u adulescu, Nonlinear eigenvalue
problems for quasilinear operators on unbounded domains, {\it
Nonlinear Differential Equations Applications (NoDEA)} {\bf 8}
(2001), 481-497.

\bibitem{MoR} D. Motreanu, V. R\u adulescu, Existence theorems for some
classes of boundary value problems involving the $p$-Laplacian",
{\it PanAmerican Math. Journal} {\bf 7} (1997), No. 2, 53-66.

\bibitem{NR} C. Niculescu, V. R\u adulescu,  A Saddle Point type theorem and applications to the study of some
problems with strong resonance at infinity, {\it Ann. Acad. Sci.
 Fennicae} {\bf 21} (1996), 117-131.

\bibitem{N} L. Nirenberg, Variational and topological methods in nonlinear problems,
{\it Bull. Amer. Math. Soc.} {\bf 4} (1981), 267-302.

\bibitem{Pa} R. Palais, Lusternik-Schnirelmann theory on Banach manifolds, {\it Topology} {\bf
5} (1966), 115-132.

\bibitem{PS1} P. Pucci, J. Serrin, Extensions of the Mountain Pass theorem, {\it J.
Funct. Anal.} {\bf 59} (1984), 185-210.

\bibitem{PS2} P. Pucci, J. Serrin, A Mountain Pass theorem, {\it J. Diff. Eq.} {\bf
60} (1985), 142-149.

\bibitem{Ra1} P. H. Rabinowitz, {\it Mini-max Methods in Critical Point Theory with
Applications to Differential Equations}, CBMS Reg. Conf.     Ser. in Math., No. 65, AMS, Providence
1986.

\bibitem{Ra2} P. H. Rabinowitz, Some critical point theorems and applications to
semilinear elliptic partial differential equations, {\it Ann. Sc.
Norm. Sup. Pisa} {\bf 2} (1978), 215-223.

\bibitem{Ra3} P. H. Rabinowitz, Free vibrations for a semilinear wave equation,
{\it Comm. Pure Appl. Math.} {\bf 31} (1978), 31-68.

\bibitem{R1} V. R\u adulescu, Sur la th\'eorie Lusternik-Schnirelmann en dimension finie, {\it Matarom,
Lab. d'Analyse Num. de l'Univ. de Paris 6}, 2(1992), 9-16.

\bibitem{R2} V. R\u adulescu, Mountain Pass theorems for non-differentiable functions
and applications, {\it Proc. Japan Acad.}  {\bf 69A} (1993),
193-198.

\bibitem{R3} V. R\u adulescu, Mountain Pass type theorems for nondifferentiable convex
functions, {\it Rev. Roum. Math. Pures Appl.} {\bf 39} (1994),
53-62.

\bibitem{R4} V. R\u adulescu, Locally Lipschitz functionals with the strong Palais-Smale
property, {\it Rev. Roum. Math. Pures Appl.} {\bf 40} (1995),
355-372.

\bibitem{R5} V. R\u adulescu, Nontrivial solutions for a multivalued
problem with strong resonance, {\it Glasgow Math. Journal} {\bf
38} (1996), 53-61.

\bibitem{R6} V. R\u adulescu, A Lusternik-Schnirelman type theorem for locally Lipschitz functionals
with applications to multivalued periodic problems, {\it Proc.
Japan Acad.} {\bf 71A} (1995), 164-167.

\bibitem{Ru} W. Rudin, {\it Functional Analysis}, Mc Graw-Hill 1973.

\bibitem{S1} M. Schechter, Nonlinear elliptic boundary value problems at strong
resonance, {\it Amer. J. Math.} {\bf 112} (1990), 439-460.

\bibitem{S2} M. Schechter, A variation of the Mountain Pass lemma and applications,
 {\it J. London Math. Soc.} {\bf 44} (1991), 491-502.

\bibitem{St} M. Struwe, {\it Variational Methods}, Springer-Verlag, Berlin Heidelberg New York, 1990.

\bibitem{SZ} A. Szulkin, {\it Critical Point Theory of Lusternik-Schnirelmann Type and
Applications to Partial Differential Equations}, Presses de l'Univ. Montreal, 1989.

\bibitem{Su} F. Sullivan, A characterization of complete metric spaces,
{\it Proc. Amer. Math. Soc.} {\bf 83} (1981), 345-346.

\bibitem{Th} K. Thews, Nontrivial solutions of elliptic equations at resonance,
{\it Proc. Royal Soc. Edinburgh}  {\bf 85A} (1980), 119-129.

\bibitem{willem} M. Willem, {\it Minimax Theorems}, Birkh\"auser,
Boston, 1996.

\bibitem{willem2}  M. Willem, {\it Analyse fonctionnelle \'el\'ementaire},
Cassini, 2003.

}
\end{thebibliography}
\end{document}